\documentclass[10pt]{amsart}  

\usepackage[utf8]{inputenc}

\usepackage{graphicx,amssymb,amsmath,amsthm,wrapfig}
\usepackage{caption}
\usepackage{subcaption}
\usepackage{indentfirst}
\usepackage{cancel}
\usepackage{epstopdf}
\usepackage{comment}
\usepackage{enumerate,color}
\usepackage{empheq}
\usepackage{fullpage}
\usepackage{hyperref}
\usepackage{xcolor}

\theoremstyle{plain}
\newtheorem{theorem}{Theorem}[section]

\newtheorem{lemma}[theorem]{Lemma}

\newtheorem{proposition}[theorem]{Proposition}
\newtheorem{remark}[theorem]{Remark}

\numberwithin{equation}{section}
\newtheorem*{theorem*}{Theorem}

\newcommand{\RNum}[1]{\uppercase\expandafter{\romannumeral #1\relax}}
\newcommand{\veps}{\varepsilon}
\newcommand{\wt}{\widetilde}
\newcommand{\UU}{\overline{U}}
\newcommand{\RR}{\mathbb{R}}
\newcommand{\UUb}{\overline{U}_\beta}
	
\newcommand{\ve}{\varepsilon}
\newcommand{\mo}{d_0}

\title{Stable \(C^{7/9}\) cusp formation for the Novikov equation} 
\author{Yunjoo Kim, Dowan Koo, Bongsuk Kwon, and Wanyong Shim}

\address{(Yunjoo Kim) Department of Mathematical Sciences, Ulsan National Institute of Science and Technology, Ulsan, 44919, Korea}
\email{gomuli3@unist.ac.kr}

\address{(Dowan Koo) Mathematical Institute, University of Oxford, Oxford OX2 6GG, United Kingdom}
\email{dowan.koo@maths.ox.ac.uk}

\address{(Bongsuk Kwon) Department of Mathematical Sciences, Ulsan National Institute of Science and Technology, Ulsan, 44919, Korea}
\email{bkwon@unist.ac.kr}

\address{(Wanyong Shim) Department of Mathematical Sciences, Korea Advanced Institute of Science and Technology, Daejeon, 34141, Korea}
\email{wyshim25@kaist.ac.kr}

\thanks{\textbf{Acknowledgment.} D.K. was supported by NRF grant no. RS-2025-02312778. B.K. was supported by the National Research Foundation of Korea(NRF) grant funded by the Korea government(MSIT)(00560003).}

\begin{document}

\begin{abstract}
We establish stable cusp formation for the Novikov equation, a cubically nonlinear Camassa--Holm-type equation. We identify an open set of smooth initial data for which the first gradient blow-up produces a cusp with sharp H\"older regularity \(C^{7/9}\). This result shows that, in nonlocal wave-breaking problems, the sharp regularity of the cusp is not determined by the nonlocal or nonlinear structure alone. While the conserved \(H^1\)-type quantity excludes the \(C^{1/3}\) cusp associated with Burgers-type gradient blow-up, the precise H\"older exponent is selected by the coupling between the nonlocal term and the algebraic structure of the nonlinearity. In the Novikov equation, this coupling yields the exponent \(7/9\), rather than the \(3/5\) exponent known for the Camassa--Holm and Hunter--Saxton equations. The main difficulty is that the naive high-frequency limit retains the cubic character of the equation and therefore does not exhibit a self-similar leading flow. We overcome this by introducing a Galilean-type change of variables around a nonzero background, which reveals a quadratic Hunter--Saxton-type leading equation. Its self-similar profiles determine the \(C^{7/9}\) cusp, while the nonlocal and cubic remainders are controlled perturbatively in modulated similarity variables.

\end{abstract}
	
	\maketitle
	

\section{Introduction}

We consider the initial value problem for the Novikov equation:
\begin{subequations} \label{Nov}
    \begin{align}
        & \label{Nov_1} v_t - v_{xxt} + 4v^2 v_x = 3v v_x v_{xx} + v^2 v_{xxx}, \quad  x \in \mathbb{R}, \ t > -\varepsilon, \\
        & \label{Nov_2} v(x,-\varepsilon) = v_0(x),
    \end{align}
\end{subequations}
where $t=-\varepsilon<0$ is the initial time. The Novikov equation \eqref{Nov_1} was first presented in \cite{HW} as a new integrable peakon equation found by Novikov, and later appeared in the classification of generalized Camassa--Holm equations in \cite{Nov}. Consistent with this origin, it shares several remarkable features with the Camassa--Holm equation, including a bi-Hamiltonian structure, infinitely many conserved quantities, and peaked solitary wave solutions known as peakons. A distinctive feature of the Novikov equation, however, is that its nonlinearity is cubic rather than quadratic. One of the main points of this paper is that, in the presence of the nonlocal Camassa--Holm-type structure, this algebraic difference changes the effective blow-up dynamics and leads to a new cusp regularity.

Over the past decade, the Novikov equation has been extensively studied from various perspectives, including local and global well-posedness, blow-up phenomena, and the stability of peakon solutions. For results on peakon stability, we refer the reader to \cite{CLWX, Laf, LLQ} and the references therein. In the one-dimensional setting, local well-posedness in Sobolev spaces $H^s(\mathbb{R})$ for $\textstyle s>\frac{3}{2}$ was established in \cite{HH, WY}, while ill-posedness at or below the critical regularity $\textstyle s=\frac{3}{2}$ has been investigated in both Sobolev and Besov settings; see \cite{HHK, NZ, WY2, YLZ, YLZ2}. Global existence and uniqueness of weak solutions have also been established under suitable structural assumptions \cite{Lai, WY2011}. Moreover, several sufficient conditions for the occurrence of gradient blow-up have been derived in the literature; see \cite{CGLQ, JN, ZY}.

Despite these advances, the fine structure of singularity formation for the Novikov equation has remained largely unexplored. In particular, the mechanism of the first gradient blow-up and the sharp H\"older regularity of the limiting singular profile were not previously known. This question is part of a broader problem in nonlinear wave breaking: beyond the fact that the \(C^1\)-norm becomes unbounded, what determines the geometry of the cusp formed at the first singularity?

The sharp H\"older regularity of the cusp is not merely a refinement of a \(C^1\) blow-up result. Rather, it identifies the critical threshold between H\"older norms that remain uniformly controlled up to the breaking time and those that fail to remain bounded. In this sense, the sharp H\"older exponent serves as a fingerprint of the underlying blow-up mechanism: different exponents correspond to different self-similar regimes and different structures of the equation.

For scalar conservation laws and Burgers-type models, this fingerprint is largely universal. Near a genuinely nonlinear breaking point, the detailed algebraic form of the flux can typically be normalized away, and the leading singularity is governed by the Burgers mechanism, producing a \(C^{1/3}\) cusp. This Burgers-type behavior is known to occur in a broad class of models, including the multidimensional Euler equations \cite{BSV} and related variants considered in \cite{BKK2, BKK, OP, Y}. Recent work has shown, however, that nonlocal Camassa--Holm-type equations can exhibit a different singularity class. In \cite{KKY}, solutions to the Camassa--Holm and Hunter--Saxton equations were shown to form cusps with sharp \(C^{3/5}\) H\"older regularity at the first blow-up of the \(C^1\)-norm. Subsequent work \cite{KKS} showed that the same \(C^{3/5}\) regime also arises in a Hamiltonian regularization of the Saint--Venant equations and of the Burgers equation. Taken together, these results show that the \(C^{3/5}\) regime is not specific to the Camassa--Holm equation and raise the question of which features select the cusp regularity.

A natural structural explanation for the difference between the Burgers $C^{1/3}$ regime and the Camassa--Holm $C^{3/5}$ regime is the presence of a uniformly bounded $H^1$-type quantity. Such a bound imposes a strict H\"older barrier on the solution profile. Indeed, for any $x,y$ in a compact set $K\subset\mathbb{R}$ and any $v\in C^1_c(\mathbb{R})$, one has
\begin{equation*}
    |v(x)-v(y)|
    \leq \left| \int_x^y v'(s)\,ds \right|
    \leq |x-y|^{1/2}\|v'\|_{L^2},
\end{equation*}
which shows that a uniform $\dot H^1$ bound yields local $C^{1/2}$ control. This simple estimate rules out a Burgers-type $C^{1/3}$ asymptotic blow-up profile. Since the Novikov equation also possesses a conserved $H^1$ quantity, one might expect the same sharp $C^{3/5}$ regularity as in the Camassa--Holm and Hunter--Saxton equations.

The main result of this paper shows that this expectation is false. We identify an open class of smooth initial data, in the relative topology determined by the modulation constraints, for which the first gradient blow-up of the Novikov equation produces a cusp of sharp $C^{7/9}$ H\"older regularity; see Theorem~\ref{mainthm}. Thus the conserved $H^1$ quantity gives only a coarse regularity barrier: it excludes the Burgers-type $C^{1/3}$ regime, but it does not determine the precise cusp exponent. The sharp exponent is instead selected by the way the nonlocal structure couples to the algebraic form of the nonlinearity. In this sense, the Novikov equation reveals a new exponent-selection mechanism in nonlocal wave-breaking problems: the resulting \(C^{7/9}\) cusp belongs to a singularity class distinct from those of both the Burgers-type \(C^{1/3}\) cusp and the Camassa--Holm-type \(C^{3/5}\) cusps.

\medskip

The cubic nonlinearity also creates a substantial analytical obstruction. In the Camassa--Holm equation, the high-frequency limit is the Hunter--Saxton equation \cite{DP, HZ}, and the two equations share the same sharp blow-up regularity \cite{KKY}. It is therefore natural to try to identify the leading equation for Novikov by the same high-frequency procedure. Introducing
\begin{equation*}
    x=\eta X, \quad t=\eta T, \quad w(X,T):=v(\eta X,\eta T),
\end{equation*}
and taking the limit $\eta\searrow0$, one obtains, after relabeling $(X,T)$ as $(x,t)$,
\begin{equation} \label{eq:high}
    w_{xxt}+3ww_xw_{xx}+w^2w_{xxx}=0.
\end{equation}
This equation captures the large-gradient scaling of the Novikov equation, but it is not the correct self-similar leading-order model. The reason is structural: equation \eqref{eq:high} retains the cubic dependence on \(w\) and does not possess the Galilean invariance needed to normalize the value of the solution at a prospective blow-up point. In particular, if \(w(x,t)\) is a solution, the shifted function \(w(x-ct,t)+c\) is not, in general, a solution. Therefore, at a point \((x_*,T_*)\) with \(w(x_*,T_*)\neq0\), one cannot reduce the height of the solution to zero, and no self-similar ansatz centered at \(x=x_*\) is available in the usual form
\begin{equation*}
    w(x,t)=(T_*-t)^\alpha
    \overline W\!\left(\frac{x-x_*}{(T_*-t)^{\alpha+1}}\right).
\end{equation*}
This obstruction prevents a direct application of the self-similar perturbation frameworks developed in \cite{KKS, KKY}.

To overcome this difficulty, we introduce a Galilean-type change of variables adapted to the cubic nonlinear structure; see \eqref{u}. This transformation is not merely a technical normalization. It is the device that exposes the hidden leading-order dynamics around a nonzero background. After the change of variables, the cubic structure decomposes into a quadratic Hunter--Saxton-type leading term and perturbative nonlocal and cubic remainders; see \eqref{Nov-u} and \eqref{leading}. The resulting leading equation has the self-similar structure needed for the blow-up analysis. The far-field behavior of its self-similar profiles determines the sharp $C^{7/9}$ regularity. This is the mechanism by which the coupling between the nonlocal term and the cubic nonlinearity produces a cusp exponent that is not visible from the conservation law alone, nor from the naive high-frequency limit.

The strategy developed here may provide a useful framework for the study of singularity formation in other evolution equations with higher-order nonlinearities. In particular, it suggests that when the naive high-frequency limit does not reveal a self-similar flow, one should look for a nonlinear change of variables adapted to the background value at the blow-up point. Such a transformation may isolate the correct self-similar leading dynamics and allow the remaining nonlocal and higher-order nonlinear terms to be treated perturbatively.

\subsection{Reformulation via a Galilean-type change of variables}\label{sec1.1}
To expose the leading-order structure in the gradient blow-up regime, we reformulate the Novikov equation. We first introduce \(p:=-v_{xt}-v^2v_{xx}-\frac12v(v_x)^2+v^3\) and rewrite \eqref{Nov_1} as the following nonlocal system:
	\begin{subequations} \label{Nov2}
		\begin{align}
			& \label{Nov2_1} v_t+v^2v_x+\frac{1}{2}v_x^3=-p_x, \\
			& \label{Nov2_2} p-p_{xx}=v^3+\frac{3}{2}vv_x^2+\frac{3}{2}v_x^2v_{xx}.
		\end{align}
	\end{subequations}
Differentiating \eqref{Nov2_1} with respect to $x$ and using \eqref{Nov2_2}, we obtain
	\begin{equation*}
		v_{xt}+v^2v_{xx}=-\frac{1}{2}vv_x^2+(v^3-p).
	\end{equation*}
Here we observe that, if
$v(x,t) \to 0 \quad \text{as } (x,t) \to (x_*,T_*)$
for some gradient blow-up point $(x_*,T_*) \in \mathbb{R} \times (-\varepsilon,\infty)$, then the coefficient $-\frac12 v$ of the quadratic term vanishes, so the usual Riccati-type blow-up mechanism for $v_x$ breaks down.\footnote{This phenomenon arises from the cubic nonlinearity of the Novikov equation and is therefore not observed in the Burgers equation or the Camassa--Holm equation (cf. \cite{KKY}).} Motivated by this observation, we focus on the nondegenerate blow-up regime, in which \(v\) remains uniformly separated from zero along the blow-up point; see Remark~\ref{center}, and see Appendix~\ref{sec1.2} for further discussion of the non-vanishing condition.

Moreover, since the system \eqref{Nov2} is invariant under spatial translations \(x \mapsto x + \delta\), we may fix the origin so that
	\begin{equation} \label{assinf}
	(\partial_x v_0)(0) = \inf_{x\in\mathbb{R}} (\partial_x v_0)(x),
	\end{equation}
	provided that $\partial_x v_0$ attains its infimum at some point in $\mathbb{R}$.	

Having fixed the origin by \eqref{assinf}, we impose the nondegeneracy condition \(v_0(0)\neq0\) and introduce a new variable \(u\) by
	\begin{equation}\label{u}
		u(x,t):=v\left( x+c^2(t+\varepsilon),t \right)-c, \quad c:=v_0(0)\neq 0.
	\end{equation}
Using the scaling invariance of \eqref{Nov2}, we set $\textstyle c = \frac{1}{2}$ without loss of generality. Under the change of variables \eqref{u}, the system \eqref{Nov2} and the initial data \eqref{Nov_2} can be rewritten as
	\begin{subequations} \label{Nov-u}
		\begin{align}
			& u_t +(u+1)uu_x+\frac{1}{2}u_x^3=-p_x, \quad  \label{Nov-u1}\\
			& p-p_{xx}=\left(u+\frac{1}{2}\right)^3+\frac{3}{2}u_x^2\left(\frac{1}{2}+u+u_{xx}\right), \label{Nov-u2}
		\end{align}
	\end{subequations}
and
\begin{equation} \label{initu}
u(x,-\veps) = u_0(x) := v_0(x) - c.
\end{equation}

This reformulation identifies, as we prove in this paper, the correct leading-order equation governing gradient blow-up. Differentiating \eqref{Nov-u1} in $x$, we have
	\begin{equation}\label{ux}
		u_{xt}+\frac{1}{2}\left(u+\frac{1}{2}\right)u_x^2+(u+1)uu_{xx}= - p + \left(u+\frac{1}{2}\right)^3.
	\end{equation}
	Formally neglecting lower-order terms in the large-gradient regime, we are led to the leading-order equation for \eqref{ux}:
	\begin{equation}\label{leading}
		w_{xt}+\frac{1}{4}w_x^2+ww_{xx}=0.
	\end{equation}
This equation is of Hunter--Saxton type, with the standard Hunter--Saxton equation recovered when the coefficient $\textstyle \frac{1}{4}$ is replaced by $\textstyle \frac{1}{2}$. We note that, unlike the formal high-frequency model \eqref{eq:high}, the equation \eqref{leading} possesses a self-similar structure. The corresponding self-similar solutions are described in the following subsection.

\subsection{Self-similar blow-up profiles}

We introduce a smooth function $\overline{U}(y): \mathbb{R} \to \mathbb{R}$ satisfying the ODE:
\begin{equation} \label{ovU}
\left( 1 + \frac{1}{4}\overline{U}'(y) \right) \overline{U}'(y) + \left( \frac{9}{2}y + \overline{U}(y) \right) \overline{U}''(y) = 0, \quad y \in \mathbb{R}.
\end{equation}
This ODE is obtained by inserting the self-similar ansatz into \eqref{leading}; see Section~\ref{sec2}.

In Proposition~\ref{Profile-construct}, we construct a one-parameter family $\{\overline{U}_\beta\}_{\beta>0}$ of smooth, odd, monotonically decreasing solutions to
\eqref{ovU}. Each profile satisfies
\begin{equation*}
\overline{U}_\beta'(0)=-4, \quad -4\leq \overline{U}_\beta'(y)<0 \quad \text{for all } y\in\mathbb R.
\end{equation*}
Moreover,
\begin{equation*}
\overline U_\beta(y) = -4y+\frac{2^{25}\beta}{6}y^3+\mathcal O(y^5) \quad \text{as } y\to 0
\end{equation*}
and, as $|y|\to\infty$,
\begin{equation*}
\overline U_\beta(y) = -\frac97 \left(\frac{1}{1296\beta}\right)^{1/9} \operatorname{sgn}(y) |y|^{7/9} + o(|y|^{7/9}), \quad |y|^{2/9} \overline U_\beta'(y) \to - \left(\frac{1}{1296\beta}\right)^{1/9}.
\end{equation*}
In particular, the far-field growth of $\overline U_\beta$ gives the sharp $C^{7/9}$ cusp in Theorem~\ref{mainthm}.

\subsection{Assumptions on the initial data}\label{Initial_subs}

We state the initial assumptions in terms of the transformed variable $u$ introduced in \eqref{u}. Through the relation \eqref{initu}, they determine the corresponding class of initial data $v_0$ for the original problem \eqref{Nov}. The resulting admissible class is discussed further in Appendix~\ref{app:admissible-data}, where we verify that it contains a relatively open set in an appropriate topology. This emphasizes that the class is not a collection of isolated or finely tuned examples, but a stable family of perturbations compatible with the modulation constraints.

For sufficiently small $\veps>0$, we assume that the initial data $u_0=u(\cdot,-\varepsilon)$ satisfy
\begin{equation}\label{in-H-C}
   u_0 + \frac{1}{2} \in H^5(\mathbb{R})\subseteq C^4(\mathbb{R}).
\end{equation}
We then place the maximal negative initial slope at $x=0$ and align the local expansion with that of the self-similar profile at the origin. We require
\begin{equation}\label{init_w_3-p}
u_0(0)=0, \quad
    \partial_xu_0(0)=-4\veps^{-1}, \quad
    \partial_x^2u_0(0)=0, \quad
    \partial_x^3u_0(0)=2^{25}\veps^{-10},
\end{equation}
along with the scale-consistent bounds
\begin{equation}\label{init_24-p}
\|\partial_x u_0\|_{L^{\infty}}\leq 4\veps^{-1}, \quad \|\partial_x^2u_0\|_{L^{\infty}}\leq C_1 \veps^{-11/2}, \quad \|\partial_x^3u_0\|_{L^{\infty}}\leq C_2 \veps^{-10}, \quad \|\partial_x^4u_0\|_{L^{\infty}}\leq C_3 \veps^{-29/2},
\end{equation}
where \(C_j>0\), $j=1,2,3$, are fixed constants independent of \(\varepsilon\), chosen sufficiently large; for definiteness, we take \(C_1=2^{12}\), \(C_2=2^{26}\), and \(C_3=2^{39}\). Let $m_0>0$ be a large constant specified later in \eqref{num_1}, and choose $N(m_0)>0$ sufficiently small so that
\begin{equation}\label{Nm0}
	N(m_0)(m_0^{2/9}+1) <\frac{1}{2}.
\end{equation}
We assume that the rescaled initial slope is close to $\overline{U}'$:
\begin{equation}\label{4.3a-p}
\left|\varepsilon(\partial_x u_0)(x)-\overline{U}' \left(\frac{x}{\veps^{9/2}}\right)\right|\leq \min\left\{\frac{ N(m_0)(\frac{x}{\veps^{9/2}})^2}{4(1+(\frac{x}{\veps^{9/2}})^2)}, \frac{5}{8\left(1+\left(\frac{x}{\veps^{9/2}}\right)^{2/9}\right)}\right\}
\end{equation}
for all $x\in\mathbb{R}$. We further impose the decay condition 
\begin{equation}\label{4.3a-p2}
\limsup_{|x|\rightarrow \infty}|x^{2/9}\partial_xu_0(x)|\leq \frac{1}{2}\left(\frac{5}{8}-\Theta\right),
\end{equation}
where $\textstyle \Theta:=\left(\frac{1}{1296}\right)^{1/9}<\frac{5}{8}$ is defined by \eqref{asymp-y-infty} with $\beta=1$. Finally, we impose
\begin{equation}\label{H0_bound}
    H(-\veps)
    =\left\|u_0+\frac{1}{2}\right\|^2_{L^2}
    +\left\|\partial_x u_0\right\|^2_{L^2}
    \leq \frac{81}{128},
\end{equation}
where $H(\cdot)$ is defined in Remark~\ref{Ham} below.

\begin{remark}[Hamiltonian conservation] \label{Ham}
The Novikov equation \eqref{Nov_1} possesses an invariant Hamiltonian integral. Multiplying \eqref{Nov_1} by $v$, we obtain
\begin{equation*}
\frac12 \frac{d}{dt} \int_{\mathbb{R}}  v^2\,dx-\int_{\mathbb{R}}vv_{xxt}\,dx =-  \int_{\mathbb{R}} 4 v^3v_x\,dx + \int_{\mathbb{R}}3v^2v_xv_{xx}\,dx+\int_{\mathbb{R}} v^3v_{xxx}\,dx.
\end{equation*}
For smooth solutions with sufficient decay at infinity, integration by parts gives
\begin{subequations}
\begin{align*}
&- \int_{\mathbb{R}} vv_{xxt}\,dx= \int_{\mathbb{R}} v_xv_{xt}\,dx =   \frac{1}{2} \frac{d}{dt} \int_{\mathbb{R}} v_x^2 \,dx, \\
&- \int_{\mathbb{R}} 4v^3v_x\,dx= - \int_{\mathbb{R}}(v^4)_x\,dx=0, \\
&\int_{\mathbb{R}}  3v^2v_xv_{xx}+v^3v_{xxx} \,dx =\int_{\mathbb{R}} (v^3v_{xx})_x\,dx=0. 
\end{align*}
\end{subequations} 
Thus, as long as the smooth solution persists, we have
\begin{equation}\label{H}
H(t):=\int_{\mathbb{R}} (v^2+v_x^2)(x,t) \, dx = H(-\veps) \quad \text{for all } t \geq -\varepsilon,
\end{equation}
This, together with a standard Sobolev inequality, yields the uniform bound on $v$:  
\begin{equation*}
|v(x,t)|\leq \frac{1}{\sqrt{2}} \| v (\cdot, t) \|_{H^1} = \frac{1}{\sqrt{2}} \sqrt{ H(-\veps ) }  =: C_v.
\end{equation*}
In particular, under the assumption \eqref{H0_bound}, we have
\begin{equation} \label{v_bound}
|v(x,t)| \leq \frac{9}{16} \quad \text{for all } x\in \mathbb{R}, \ t \geq -\varepsilon.
\end{equation}
\end{remark}

\subsection{Main result}
For a bounded open set $\Omega\subset\mathbb R$ and $\alpha\in(0,1]$, we write
\begin{equation*}
	[w]_{C^{\alpha}(\Omega) } := \sup_{x,y\in \Omega, x\ne y} \frac{|w(x) - w(y) | }{|x-y|^{\alpha}}
\end{equation*}
for the H\"older seminorm of $w$ on $\Omega$. We now state the main theorem.

\begin{theorem}\label{mainthm}
There exist constants $\ve_0>0$ and $C>0$ such that the following statement holds.

For each $\veps \in (0,\veps_0)$, let $u_0$ satisfy \eqref{in-H-C}--\eqref{H0_bound}, and define $\textstyle v_0(x) := u_0(x) + \frac{1}{2}$. Then the initial value problem \eqref{Nov} admits a unique solution
\begin{equation*}
v\in C([-\ve, T_*); H^5(\mathbb{R}) ) \cap  C^1([-\ve, T_*); H^4(\mathbb{R}) ),
\end{equation*}
which blows up in the $C^1$-norm at a finite time $T_*$ satisfying $|T_*| \leq C \veps^{11/6}$. Moreover, for any bounded open subset $\Omega\subset\mathbb{R}$, it holds that 
\begin{equation*}
\sup_{t<T_*} \left[ v(\cdot, t) \right]_{C^{7/9}(\Omega)}<\infty;
\end{equation*}
and for any $\textstyle \alpha \in (\frac79,1]$, 
\begin{equation*}
\begin{cases}
\displaystyle \lim_{t\nearrow T_*} \left[ v(\cdot,t) \right]_{C^\alpha(\Omega)}=\infty 
& \text{if } x_\ast\in\Omega,\\
\displaystyle \limsup_{t\nearrow T_*} \left[ v(\cdot,t) \right]_{C^\alpha(\Omega)}<\infty 
& \text{if } x_\ast\notin\overline{\Omega},
\end{cases}
\end{equation*}	
where $x_\ast$ is the blow-up location.\footnote{$x_\ast$ is defined in the proof of Theorem~\ref{mainthm}.} Furthermore, for any bounded open set $\Omega$ containing $x_\ast$ and for $\textstyle \alpha \in (\frac79,1]$, the temporal blow-up rate is given by 
\begin{equation}\label{blow-up-rate}
\left[ v(\cdot, t) \right]_{C^\alpha(\Omega)}\sim (T_*-t)^{-\frac{9\alpha-7}{2}}
\end{equation} 
for all $t$ sufficiently close to $T_*$.\footnote{Here, $A(t) \sim B(t)$ means that $C^{-1} B(t) \le A(t) \le C B(t)$ for some $C>0$ independent of $t$.}
\end{theorem}

The proof of Theorem~\ref{mainthm} is provided in Section~\ref{C13_subsec-p}. Theorem~\ref{mainthm} shows that solutions evolving from the smooth admissible initial data described in Section~\ref{Initial_subs} form $C^{7/9}$-type cusps at their first finite-time singularity, in the sense that the local $C^{7/9}$ H\"older seminorm remains uniformly bounded up to the blow-up time $T_*$, whereas, for every $\alpha>\frac{7}{9}$, the local $C^\alpha$ H\"older seminorm blows up as $t\nearrow T_*$ in any neighborhood of the singular point.

\subsection{Outline of the proof}

We briefly describe the main ideas of the proof. The argument begins with the construction of the self-similar profiles associated with the leading-order equation \eqref{leading}. Substituting the ansatz 
\begin{equation*}
w(x,t)=(T_\ast-t)^{7/2}\overline U((x-x_\ast)/(T_*-t)^{9/2})
\end{equation*}
into \eqref{leading} leads to the profile equation \eqref{ovU}. In Section~\ref{sec2}, we construct a one-parameter family of smooth odd monotone solutions to this equation. The key feature of these profiles is their far-field behavior \(\overline U(y)\sim -\operatorname{sgn}(y)|y|^{7/9}\), or equivalently \(\overline U'(y)\sim - |y|^{-2/9}\). This asymptotic behavior underlies the sharp \(C^{7/9}\) regularity.

With these profiles in hand, we then reformulate the Novikov equation \eqref{Nov-u} in modulated similarity variables \((y,s)\) adapted to the profile \(\overline U\), with rescaled unknown \(U=U(y,s)\). The modulation parameters \(\tau,\kappa,\xi\) track the blow-up time, the height, and the location of the singularity. They are chosen so that the normalized conditions 
\begin{equation*}
U(0,s)=0, \quad U_y(0,s)=-4, \quad U_{yy}(0,s)=0
\end{equation*}
are preserved for all \(s\). In these variables, the leading transport structure is governed by \(\frac92 y+\overline U\), while the remaining terms coming from the full Novikov equation appear as perturbative contributions involving the nonlocal pressure \(P\), the modulations, and the cubic remainder.

The main technical part of the proof is the stability analysis around \(\overline U\). We formulate a bootstrap scheme adapted both to the behavior near the self-similar origin \(y=0\) and to the far-field decay of the profile. The weighted estimate for \(U_y-\overline U'\) is the crucial bound: it propagates the \(|y|^{-2/9}\) decay of the slope and therefore determines the \(C^{7/9}\) cusp. A main difficulty is the treatment of the nonlocal term \(P\). First, the bounds for \(P_y,P_{yy},P_{yyy}\) require a careful use of the cancellations between the singular parts of the differentiated Green function and the cubic terms in the source; see Lemma~\ref{P_high_lem}. Second, \(P\) itself does not have the far-field decay needed for the weighted estimate on \(U_y-\overline U'\). We therefore exploit the cancellation in the combination
\begin{equation*}
P-\left(e^{-7s/2}U+\kappa+\frac12\right)^3,
\end{equation*}
proved in Lemma~\ref{Pdec_lem}. These two estimates are essential for closing the bootstrap estimates involving \(P\).

Finally, we use the stability estimates to prove Theorem~\ref{mainthm}. The lower bound for the slope near \(y=0\) yields divergence of the local \(C^\alpha\) seminorm for all \(\alpha \in (\frac79,1]\), whereas the far-field decay of \(U_y\) gives uniform \(C^{7/9}\) control and bounded \(C^\alpha\) seminorms away from the singular point. The same scaling argument gives the temporal blow-up rate \((T_*-t)^{-(9\alpha-7)/2}\).

\medskip

\emph{Plan of the paper}. Section~\ref{sec2} constructs the self-similar blow-up profiles and records the profile inequalities used in the stability analysis. Section~3 introduces the modulated similarity variables, formulates the bootstrap argument, and proves Theorem~\ref{mainthm} using the bootstrap bounds. Section~\ref{sec4} closes the bootstrap estimates. The appendix contains a local non-vanishing lemma, the construction of admissible initial data, auxiliary inequalities for \(\overline U\), and a maximum principle lemma used in the stability analysis.

\section{Construction of self-similar blow-up profiles} \label{sec2}
To construct self-similar solutions to the leading-order equation \eqref{leading}, we consider solutions of the form
\begin{equation}
\label{ss-general}
w(x,t)=(T_*-t)^\alpha \overline{U}(y), \quad y:=\frac{x-x_*}{(T_*-t)^{\alpha+1}},
\end{equation}
where \((x_*,T_*)\) denotes a prescribed spacetime blow-up point. Substituting \eqref{ss-general} into \eqref{leading}, we obtain
\begin{equation} \label{a-ODE}
\left(1+\frac14 \overline{U}'(y)\right)\overline{U}'(y) +\left((\alpha+1)y+\overline{U}(y)\right)\overline{U}''(y)=0.
\end{equation}

For the class of profiles considered here, the exponent \(\alpha\) is fixed by requiring the profile \(\overline U\) to be smooth and odd, with \(\overline U'(0)<0\) and \(\overline U'''(0)\neq0\). By oddness,
$\overline U(0)=\overline U''(0)=0$.
Evaluating \eqref{a-ODE} at \(y=0\), we obtain \(\overline U'(0)(1+\frac14\overline U'(0))=0\), and the condition \(\overline U'(0)<0\) therefore forces
$\overline U'(0)=-4$.
Differentiating \eqref{a-ODE} twice and evaluating at \(y=0\) yields \((2\alpha-7)\overline U'''(0)=0\). This, together with the nondegeneracy condition \(\overline U'''(0)\neq0\), gives
\begin{equation*}
\alpha=\frac72.
\end{equation*}

With this choice of exponent, we arrive at the self-similar ansatz
\begin{equation*}
w(x,t)= (T_*-t)^{7/2} \overline U\left( \frac{x-x_*}{(T_*-t)^{9/2}} \right).
\end{equation*}
The corresponding profile equation is \eqref{ovU}, which we recall as
\begin{equation*}
\left(1+\frac14\overline U'(y)\right)\overline U'(y) +\left(\frac92 y+\overline U(y)\right)\overline U''(y)=0, \quad y\in\mathbb{R}.
\end{equation*}
In the following proposition, we establish the existence of a one-parameter family \(\{\overline U_\beta\}_{\beta>0}\) of smooth odd monotone solutions to \eqref{ovU}. These profiles provide the leading-order description of singularity formation for the Novikov equation.

	\begin{proposition}\label{Profile-construct}
	The ODE \eqref{ovU} admits a one-parameter family of smooth solutions $\{ \overline{U}_\beta : \beta>0 \}$ such that 
	for each $\beta>0$, 
	$\overline{U}_\beta(y)$ is monotonically decreasing on $\mathbb{R}$ and is odd, i.e., $\UUb(-y) = - \UUb(y)$ for all $y\in\mathbb{R}$. 
	Furthermore, it holds that  
	\begin{equation}\label{U_taylor_small}
		\overline{U}_\beta (y) = -4y+\frac{2^{25} \beta}{6} y^3+\mathcal{O}(y^5) \quad\text{for }|y|\ll 1 ,
	\end{equation}
	and
	\begin{equation} \label{U_taylor_far}
		\overline{U}_\beta (y) = -\frac{9}{7}\Theta y^{7/9}+ o(y^{7/9}) \quad \text{for } |y|\gg 1,
	\end{equation}
	where $\textstyle \Theta:=\left(\frac{1}{1296\beta}\right)^{1/9}$.
	In particular,
	\begin{equation}\label{UU-0}
		\UU_\beta(0) = 0, \quad \UU_\beta'(0) =-4, \quad \UU_\beta''(0) = 0, \quad\UU_\beta^{(3)} (0) = 2^{25} \beta, 
	\end{equation} 
	and, as $|y|\rightarrow \infty$,
	\begin{equation}\label{asymp-y-infty}
		y^{-7/9}  \UU_\beta(y)  \to -\frac{9}{7}\Theta, \quad y^{2/9}  \UU_\beta'(y)  \to -\Theta, \quad y^{11/9}  \UU_\beta''(y)  \to \frac{2}{9}\Theta. 
	\end{equation} 
\end{proposition}

\begin{proof}
We first construct the desired one-parameter family of solutions. Since the equation \eqref{ovU}  is invariant under the odd reflection, it suffices to work on the half-line $\{y \in \mathbb{R} : y \geq 0 \}$. We set $\textstyle W(y) := \overline{U}(y) + \frac{9}{2}y$. Then \eqref{ovU} can be reformulated as the \emph{autonomous} ODE:
	\begin{equation}\label{Weq}
		\left(W'-\frac{1}{2}\right)\left(W'-\frac{9}{2}\right)+4WW''=0.
	\end{equation}
We look for solutions \(W\) to \eqref{Weq} satisfying \(W'=V(W)\) for some function \(V\). Using the relation $\textstyle W''=V\frac{dV}{dW}$, we obtain
	\begin{equation*} 
		\left(V-\frac{1}{2}\right) \left(V-\frac{9}{2}\right)+  4WV\frac{dV}{dW}=0.
	\end{equation*}
By separation of variables, this gives
	\begin{equation}\label{CW2}
		\beta W^2=\frac{\left(1-2V\right) }{\left(2V-9\right)^9}, 
	\end{equation}
where $\beta>0$ is a constant of integration. For each $\beta>0$, the relation \eqref{CW2} defines a smooth function $V(W)$ on $\mathbb{R}$ such that
	\begin{equation}\label{Vbound}
		\frac{1}{2}\leq V(W)<\frac{9}{2}, \quad \lim_{|W|\rightarrow \infty}V(W)=\frac{9}{2}, \quad V(0) = \frac12,
	\end{equation}
	and
	\begin{equation} \label{V'_rel}
		\frac{dV}{dW}=\frac{\beta W(2V-9)^{10}}{16V}.
	\end{equation}
By the construction of $V$, together with \eqref{Vbound} and \eqref{V'_rel}, $V$ is smooth and uniformly bounded on $\mathbb R$, and $\textstyle \frac{dV}{dW}$ is uniformly bounded on $\mathbb R$. Hence, by standard ODE theory, for each $\beta>0$, the initial value problem	
\begin{equation} \label{WVode}
W'=V(W),\quad W(0)=0
\end{equation}
admits a unique smooth solution on $[0,\infty)$. The odd extension of $W$ to $\mathbb R$ is then a smooth solution of \eqref{Weq} on the whole real line, and consequently $\textstyle \overline{U}_\beta(y):=W(y)-\frac{9}{2} y$ is an odd solution of \eqref{ovU} on $\mathbb{R}$. That is, 
	\begin{equation} \label{Ub-eq}
		\left(1+\frac{1}{4}\UUb'\right)\UUb'+\left(\UUb+\frac{9}{2}y\right)\UUb''=0, \quad y \in \mathbb{R}.
	\end{equation}
Moreover, since $\textstyle \overline{U}_\beta'=V(W)-\frac{9}{2}$, \eqref{Vbound} gives
\begin{equation} \label{Wbar_neg}
-4\le \overline{U}_\beta'(y)<0, \quad y\in\mathbb{R}.
\end{equation}
Thus $\overline{U}_\beta$ is decreasing.

	Next, we prove \eqref{U_taylor_small} for the function $\overline{U}_\beta$ constructed above. We first verify \eqref{UU-0}:
	\begin{subequations} 
		\begin{align}
			&\UUb(0)=W (0)=0, \label{UU-0-pf0}
			\\
			&\UUb'(0) = W'(0) -\frac92 = -4,   \label{UU-0-pf1}
			\\
			&\UUb''(0) = W''(0) = 0, \label{UU-0-pf2}
			\\
			&\UUb^{(3)} (0) = W^{(3)}(0) = 2^{25}\beta. \label{UU-0-pf3}
		\end{align}
	\end{subequations}
Indeed, the oddness of $W$ gives $W(0)=0$ and $W''(0)=0$, so \eqref{UU-0-pf0} and \eqref{UU-0-pf2} hold. The identity \eqref{UU-0-pf1} follows from \eqref{WVode} and $\textstyle V(0) = \frac{1}{2}$. Differentiating \eqref{V'_rel} with respect to $W$ gives $V''(0)=2^{27}\beta$. Hence, differentiating $W'=V(W)$ twice in $y$ and using $W''(0)=0$, we obtain
\begin{equation*}
W^{(3)}(0)=W'(0)^2V''(0)=2^{25}\beta,
\end{equation*}
which proves \eqref{UU-0-pf3}. Then, using Taylor's theorem together with \eqref{UU-0-pf0}--\eqref{UU-0-pf3}, we obtain \eqref{U_taylor_small}.

It remains to prove \eqref{U_taylor_far} and \eqref{asymp-y-infty}. From \eqref{Wbar_neg} and \eqref{UU-0-pf0}, we obtain
	\begin{equation}\label{U-rough}  
		y\UUb(y) \le 0, \quad | \UUb(y) | \le 4 |y|, 
		\quad y\in\mathbb{R}.
	\end{equation}
By the oddness of \(\UUb\), it suffices to consider \(y\ge0\). Then \eqref{U-rough} gives
$W(y) = \UUb(y) + \frac{9}{2} y\ge \frac{1}{2} y\ge0$.
Hence, using \eqref{CW2}, we obtain
	\begin{equation}\label{3.4}
		\UUb+\frac{9}{2}y= \left(\frac{4+\overline{U}_\beta'}{2^8\beta (-\overline{U}_\beta')^9}\right)^{1/2}, \quad y \geq 0.
	\end{equation}
By substituting \eqref{3.4} into \eqref{Ub-eq}, we have 
	\begin{equation}\label{U''_TEMP}
		\overline{U}_\beta''=4\sqrt{\beta}(4+\overline{U}_\beta')^{1/2}(-\overline{U}_\beta')^{11/2}, \quad y \ge 0.
	\end{equation}
Set
\begin{equation*}
\eta:=\left(\frac{4+\overline U_\beta'}{-\overline U_\beta'}\right)^{1/2}.
\end{equation*}
Then, by separation of variables and \eqref{UU-0-pf1}, we obtain
\begin{equation*}
\begin{split}
4\sqrt{\beta}y &=\frac{1}{512}\int_0^\eta (1+\xi^2)^4 \, d\xi \\ 
&=\frac{\eta}{512} \left( \frac{\eta^8}{9} + \frac{4}{7}\eta^6 + \frac{6}{5}\eta^4 +\frac{4}{3}\eta^2 + 1 \right) \\
&=\frac{(4+\overline{U}'_\beta)^{1/2}}{1260(-\overline{U}'_\beta)^{9/2}} \left( (\overline{U}_\beta')^4 -2(\overline{U}_\beta')^3 +6(\overline{U}_\beta')^2 -20\overline{U}_\beta' +70 \right),
	\end{split}
\end{equation*}
or, equivalently,
	\begin{equation}\label{U0FAR}
		-y^{2/9}\overline{U}_\beta'=(5040\sqrt{\beta})^{-2/9}(4+\overline{U}_\beta')^{1/9}\left((\overline{U}_\beta')^4-2(\overline{U}_\beta')^3+6(\overline{U}_\beta')^2-20(\overline{U}_\beta')+70\right)^{2/9}. 
	\end{equation}
We now use this identity to express the quantities appearing in \eqref{asymp-y-infty} in terms of  $\overline{U}_\beta'$. Combining \eqref{3.4} with \eqref{U0FAR}, we have
	\begin{equation*}
		\begin{split}
			\overline{U}_\beta&=-\frac{9}{2}y+\frac{(4+\overline{U}_\beta')^{1/2}}{16\sqrt{\beta}(-\overline{U}_\beta')^{9/2}}=\frac{-72\sqrt{\beta}(-y^{2/9}\overline{U}_\beta')^{9/2}+(4+\overline{U}_\beta')^{1/2}}{16\sqrt{\beta}(-\overline{U}_\beta')^{9/2}}
			\\
			&=\frac{1}{1120\sqrt{\beta}}\frac{(4+\overline{U}_\beta')^{1/2}\left((\overline{U}'_\beta)^3-2(\overline{U}'_\beta)^2+6(\overline{U}'_\beta)-20\right)}{(-\overline{U}_\beta')^{7/2}}.
		\end{split}
	\end{equation*}
	Together with \eqref{U0FAR}, this yields
	\begin{equation}\label{Ubar35}
		\begin{split}
			\frac{\overline{U}_\beta}{y^{7/9}}&=\frac{1}{1120\sqrt{\beta}}\frac{(4+\overline{U}_\beta')^{1/2}\left((\overline{U}'_\beta)^3-2(\overline{U}'_\beta)^2+6(\overline{U}'_\beta)-20\right)}{(-y^{2/9}\overline{U}_\beta')^{7/2}}
			\\
			&=\frac{(5040\sqrt{\beta})^{7/9}}{1120\sqrt{\beta}}\frac{(4+\overline{U}_\beta')^{1/9}\left((\overline{U}'_\beta)^3-2(\overline{U}'_\beta)^2+6(\overline{U}'_\beta)-20\right)}{\left((\overline{U}_\beta')^4-2(\overline{U}_\beta')^3+6(\overline{U}_\beta')^2-20(\overline{U}_\beta')+70\right)^{7/9}}.
		\end{split}
	\end{equation}
	Multiplying \eqref{U''_TEMP} by $y^{11/9}$ and using \eqref{U0FAR}, we obtain
	\begin{equation}\label{Ubar75}
		\begin{split}
			y^{11/9}\overline{U}_\beta''&=\frac{1}{1260(5040\sqrt{\beta})^{2/9}}(4+\overline{U}_\beta')^{10/9}\left((\overline{U}_\beta')^4-2(\overline{U}_\beta')^3+6(\overline{U}_\beta')^2-20(\overline{U}_\beta')+70\right)^{11/9}.
		\end{split}
	\end{equation}
Since $W(y)\to \infty$ as $y\to \infty$, \eqref{Vbound} and $ \textstyle \UUb'(y)=V(W(y))-\frac{9}{2}$ imply $\textstyle\lim_{y \to \infty} \overline{U}_\beta'(y)=0$. Taking $y\to\infty$ in \eqref{U0FAR}, \eqref{Ubar35}, and \eqref{Ubar75}, we obtain \eqref{asymp-y-infty} as $y\to\infty$. The case $y\to-\infty$ follows from the oddness of $\UUb$. Hence \eqref{asymp-y-infty} holds, and \eqref{U_taylor_far} follows.
\end{proof}

Throughout our analysis, we fix the normalized reference profile \(\overline U:=\overline U_1\). See Remark~\ref{rem:relax} for the role of the parameter \(\beta\) in the initial assumptions and for the corresponding relaxation of the normalized admissible class.

For later use, we collect several inequalities for \(\overline U\). Their proofs are deferred to Appendix~\ref{sec_U}. It holds that
\begin{subequations} \label{num_4236}
	\begin{align}
		&4+\overline{U}'(y) - \frac{3 y^2}{1+y^2}\geq 0, \quad y\in\mathbb{R},    \label{num_4}
		\\
		&				1+\frac{\overline{U}'}{2}+\frac{2}{y(1+y^2)}\left(\frac{9y}{2}+\overline{U}\right)\geq \frac{y^2}{2(1+y^2)}, \quad y\in\mathbb{R},  \label{num_2}
		\\
		&				\frac{11}{2}+\frac{3}{2}\overline{U}'+\frac{1}{y(1+y^2)}\left(\frac{9y}{2}+\overline{U}\right) \geq \frac{4y^2}{1+y^2}, \quad y\in\mathbb{R}.  \label{num_3}
	\end{align}
\end{subequations}
There exists $\textstyle \delta\in(\frac{1}{2},1)$ such that
\begin{equation}
	|\overline{U}''(y)|\frac{y^2+1}{y^2}\int^{|y|}_0{\frac{y'^2}{1+y'^2 }\,dy'}
	\leq  
	\delta \left(1+\frac{\overline{U}'}{2}+\frac{2}{y(1+y^2)}\left(\frac{9y}{2}+\overline{U}\right)\right)-\frac{\delta y^2}{48(1+y^2)}, \quad y \in \mathbb{R}. \label{num_6}
\end{equation}
Moreover, there exists a sufficiently large $m_0>0$ satisfying
\begin{multline}\label{num_1}
	\frac{63}{62}(y^{2/9}+1)|\overline{U}''(y)|\int^{|y|}_0\frac{dy'}{1+y'^{2/9}}
\\
\leq 1-\frac{1}{4(y^{2/9}+1)}+\frac{\overline{U}'}{2}-\frac{2y^{2/9}}{9(y^{2/9}+1)}\left(\frac{9}{2}+\frac{\overline U}{y}+\frac{1}{y}\int_{0}^{y}{\frac{dy'}{y'^{2/9}+1}}\right), \quad |y|\geq m_0.
\end{multline}

\section{Stability estimates and proof of main theorem} 
In this section, we establish stability estimates around the reference self-similar profile $\overline{U}=\overline{U}_1$ and use them to prove Theorem~\ref{mainthm}. We first introduce dynamically modulated self-similar variables adapted to $\overline{U}$ and rewrite the system \eqref{Nov-u} in these variables. We then obtain the stability estimates by a bootstrap argument. Finally, the proof of Theorem~\ref{mainthm} is given in Subsection~\ref{C13_subsec-p}.

\subsection{Self-similar variables and modulations}\label{modul_section}
We define dynamic modulation functions $\tau, \kappa, \xi: [-\varepsilon,T) \to \mathbb{R}$ for some $T>-\varepsilon$, satisfying a system of ODEs:
\begin{subequations}
\begin{align}
\dot{\tau}(t) & = -2\kappa(t) - \frac{(\tau(t)-t)^2p(\xi(t),t)}{4} + \frac{(\tau(t)-t)^2}{4} \left( \kappa(t) + \frac{1}{2} \right)^3, \label{tau'}\\
\dot{\kappa}(t) & = - \frac{4(\tau(t)-t)^{-11/2}}{\partial_x^3u(\xi(t),t)} \left( (\tau(t)-t)^{9/2} p_x(\xi(t),t) + 12 (\tau(t)-t)^{7/2} \left( \kappa(t) + \frac{1}{2} \right)^2 - 32(\tau(t)-t)^{3/2} \right) \label{kappa'}\\
& \quad - p_x(\xi(t),t) + 32 (\tau(t)-t)^{-3}, \nonumber\\
\dot{\xi}(t) & = \frac{(\tau(t)-t)^{-9/2}}{\partial_x^3u(\xi(t),t)} \left( (\tau(t)-t)^{9/2} p_x(\xi(t),t) + 12 (\tau(t)-t)^{7/2} \left( \kappa(t) + \frac{1}{2} \right)^2 - 32(\tau(t)-t)^{3/2} \right) \label{xi'}\\
& \quad + (\kappa(t)+1)\kappa(t)\nonumber
\end{align}
\end{subequations}
with the initial data
\begin{equation}\label{init_modul}
\tau(-\varepsilon) = 0, \quad \kappa(-\varepsilon) =0, \quad \xi(-\varepsilon)=0.
\end{equation}
Here we remark that the initial values $\tau(-\veps)$ and $\xi(-\veps)$ are chosen freely, while $\kappa(-\veps)$ is determined once $\xi(-\veps)=0$ is fixed, consistently with the self-similar change of variables introduced below. Since $\xi(-\varepsilon)=0$, the definition of $c$ gives
\begin{equation*}
	\kappa(-\veps)=u_0(0)=v_0(0)-c=0.
\end{equation*}
We also define $T_*$ as the first time at which $\tau(t)=t$, namely,
\begin{equation}\label{T-star}
T_* := \inf \{ t \in [-\varepsilon,\infty): \tau(t)=t \}.
\end{equation}
In Section~\ref{C13_subsec-p}, we will prove that \(T_*\) is well-defined and finite, and that it coincides with the first blow-up time of the solution.

We introduce the self-similar variables
\begin{equation} \label{ys}
y(x,t) = \frac{x - \xi(t)}{(\tau(t)-t)^{9/2}}, \quad s(t) = -\log (\tau(t)-t),
\end{equation}
and define new functions $U$ and $P$ by
\begin{equation} \label{UP}
u(x,t) - \kappa(t) = e^{-7s/2} U(y,s), \quad p(x,t) = P(y,s).
\end{equation}
By a straightforward calculation, one has
\begin{equation*}
\partial_t = \left( -\dot{\xi} e^{9s/2} + \frac{9}{2}y(1-\dot{\tau}) e^s \right) \partial_y + (1-\dot{\tau})e^s \partial_s, \quad \partial_x = e^{9s/2}\partial_y.
\end{equation*}
Using these relations, we rewrite \eqref{Nov-u} under the transformations \eqref{ys}--\eqref{UP} as
\begin{subequations} \label{Nov-UP}
	\begin{align}
		& U_s - \frac{7}{2} U + \left(\mathcal{U}+\frac{e^{11s/2}U_y^2}{2(1-\dot{\tau})} \right)U_y = - \frac{e^{5s/2}\dot{\kappa}}{1-\dot{\tau}} - \frac{e^{7s} P_y}{1-\dot{\tau}},\label{Ueq}
		\\
		& P - e^{9s} P_{yy} = \left(e^{-7s/2}U + \kappa + \frac{1}{2}\right)^3+\frac{3}{2}e^{2s}U_y^2\left(\frac{1}{2}+e^{-7s/2}U + \kappa + e^{11s/2}U_{yy}\right),\label{Peq}
	\end{align}
\end{subequations}
where
\begin{equation}\label{U}
	\mathcal{U}(y,s) := \frac{9}{2}y + \frac{(e^{-7s/2}U + 2\kappa + 1)U}{1-\dot{\tau}} + e^{7s/2} \frac{(\kappa+1)\kappa - \dot{\xi}}{1-\dot{\tau}}.
\end{equation}
Applying $\partial_y^n$ to \eqref{Ueq} for $n=1,2,3,4$, we have
\begin{subequations}
\begin{align}
& \left( \partial_s + 1 +  \frac{(e^{-7s/2}U + \kappa)U_y}{2(1-\dot{\tau})} + \frac{U_y}{4(1-\dot{\tau})}  \right) U_y + \mathcal{U} U_{yy}\label{Uy_der1}\\
& \quad = - \frac{e^{-2s}P}{1-\dot{\tau}} + \frac{e^{-2s}}{1-\dot{\tau}} \left(e^{-7s/2}U + \kappa + \frac{1}{2}\right)^3, \nonumber\\
& \left( \partial_s + \frac{11}{2} + \frac{3(e^{-7s/2}U + \kappa)U_y}{1-\dot{\tau}} + \frac{3U_y}{2(1-\dot{\tau})} \right) U_{yy} + \mathcal{U} \partial_y^3 U \label{Uy_der2}\\
& \quad = - \frac{e^{-2s}P_y}{1-\dot{\tau}} + \frac{3e^{-11s/2}U_y}{1-\dot{\tau}} \left(e^{-7s/2}U + \kappa + \frac{1}{2}\right)^2  - \frac{e^{-7s/2}U_y^3}{2(1-\dot{\tau})}, \nonumber\\
& \left( \partial_s + 10 + \frac{5(e^{-7s/2}U + \kappa)U_y}{1-\dot{\tau}} + \frac{5U_y}{2(1-\dot{\tau})} \right) \partial_y^3 U + \mathcal{U} \partial_y^4 U \label{Uy_der3}\\
& \quad = - \frac{e^{-2s}P_{yy}}{1-\dot{\tau}} + \frac{3e^{-11s/2}U_{yy}}{1-\dot{\tau}} \left(e^{-7s/2}U + \kappa + \frac{1}{2}\right)^2 + \frac{6e^{-9s}U_y^2}{1-\dot{\tau}} \left(e^{-7s/2}U + \kappa + \frac{1}{2}\right) \nonumber\\
& \qquad - \left( \frac{9e^{-7s/2}U_y^2}{2(1-\dot{\tau})} + \frac{3(e^{-7s/2}U +\kappa)U_{yy}}{1-\dot{\tau}} + \frac{3U_{yy}}{2(1-\dot{\tau})} \right) U_{yy}, \nonumber
\end{align}
\end{subequations}
and
\begin{equation} \label{Uy_der4}
\begin{split}
& \left( \partial_s + \frac{29}{2} + \frac{7(e^{-7s/2}U + \kappa)U_y}{1-\dot{\tau}} + \frac{7U_y}{2(1-\dot{\tau})}  \right) \partial_y^4 U + \mathcal{U} \partial_y^5 U \\
& \quad = - \frac{e^{-2s}\partial_y^3 P}{1-\dot{\tau}} + \frac{3e^{-11s/2}\partial_y^3 U}{1-\dot{\tau}} \left(e^{-7s/2}U + \kappa + \frac{1}{2}\right)^2 + \frac{18e^{-9s}U_yU_{yy}}{1-\dot{\tau}} \left(e^{-7s/2}U + \kappa + \frac{1}{2}\right) \\
& \qquad + \frac{6e^{-25s/2}U_y^3}{1-\dot{\tau}} - \left( \frac{19e^{-7s/2}U_y^2}{2(1-\dot{\tau})} + \frac{11(e^{-7s/2}U +\kappa)U_{yy}}{1-\dot{\tau}} + \frac{11U_{yy}}{2(1-\dot{\tau})} \right) \partial_y^3 U - \frac{12e^{-7s/2}U_yU_{yy}^2}{1-\dot{\tau}}.
\end{split}
\end{equation}

We impose the constraints
\begin{equation}\label{constraint-p}
U(0,s) = 0, \quad U_y(0,s) = -4, \quad U_{yy}(0,s) = 0 \quad \text{for all } s,
\end{equation}
so that $\overline{U}$ and $U(\cdot,s)$ have the same value and first two derivatives at $y=0$; see \eqref{UU-0}. Evaluating \eqref{Ueq}, \eqref{Uy_der1}, and \eqref{Uy_der2} at $y=0$ and using \eqref{constraint-p}, we obtain
\begin{subequations}
\begin{align}
& \dot{\tau} = - 2 \kappa - \frac{e^{-2s}P(0,s)}{4} + \frac{e^{-2s}}{4} \left( \kappa + \frac{1}{2} \right)^3, \label{tau}\\
& \dot{\kappa} =  - \frac{4e^{-9s/2}}{\partial_y^3U(0,s)} \left( P_y(0,s) + 12 e^{-7s/2} \left( \kappa + \frac{1}{2} \right)^2 - 32e^{-3s/2} \right) - e^{9s/2}P_y(0,s) + 32 e^{3s}, \label{kappa}\\
& \dot{\xi} = \frac{e^{-11s/2}}{\partial_y^3 U(0,s)} \left( P_y(0,s) + 12 e^{-7s/2} \left( \kappa + \frac{1}{2} \right)^2 - 32e^{-3s/2} \right) + (\kappa + 1)\kappa. \label{xi}
\end{align}
\end{subequations}
We note that the system of ODEs \eqref{tau}--\eqref{xi} is equivalent to \eqref{tau'}--\eqref{xi'}. A direct verification shows that 
\begin{equation*}
	(U(0,s),U_y(0,s), U_{yy}(0,s))\equiv (0,-4,0)
\end{equation*}
solves the equations \eqref{Ueq}, \eqref{Uy_der1}, and \eqref{Uy_der2}, while satisfying the condition \eqref{constraint-p} at the initial time $s=s_0:=-\log \veps$. By uniqueness, as long as the solution $(\tau,\kappa,\xi)$ to the initial problem \eqref{tau'}--\eqref{xi'} with \eqref{init_modul} exists, the corresponding function $U$ necessarily preserves the constraints \eqref{constraint-p} for all $s\geq s_0$. In Section~\ref{C13_subsec-p}, we will show that the modulation functions \((\tau,\kappa,\xi)\) exist locally in time and can be extended, by a continuation argument, to the interval \([-\varepsilon,T_*)\).

\subsection{Bootstrap argument}
We obtain the a priori stability estimates in self-similar variables by a bootstrap argument around the reference profile \(\overline U\). These estimates are the key input in the proof of Theorem~\ref{mainthm}: the weighted bound on \(U_y-\overline U'\), together with the far-field asymptotics of \(\overline U\), gives the global decay of \(U_y\) used in the \(C^{7/9}\) estimate, while the remaining bounds control the solution near \(y=0\) and provide the higher-derivative estimates needed for continuation.

In the proposition below, we write $\overline{U}(y):=\overline{U}_1(y)$, where $\overline{U}_1$ denotes the smooth solution to \eqref{ovU} constructed in Proposition~\ref{Profile-construct}. We also fix a sufficiently large constant $m_0>0$ satisfying \eqref{num_1}, and a sufficiently small constant $N(m_0)>0$ satisfying \eqref{Nm0}.

\begin{proposition}[Bootstrap improvement]\label{Boot}
	There exist constants $\veps_0>0$ and $M_0>0$ such that, if $M\ge M_0$ and $\veps\in(0,\veps_0)$, then the following statement holds.
	
Let $\sigma_1>s_0:=-\log\veps$, and let $U$ be a solution of \eqref{Nov-UP} on $\mathbb R\times[s_0,\sigma_1]$, satisfying \eqref{constraint-p}, with $(\tau,\kappa,\xi)$ defined by \eqref{tau}--\eqref{xi} and \eqref{init_modul}. Assume that, for all $(y,s)\in\mathbb R\times[s_0,\sigma_1]$,
	\begin{subequations}\label{Boot_2-p}
		\begin{align}
			&|U_y(y,s)-\overline{U}'(y)| \leq  \frac{ N(m_0)y^2}{1+y^2} ,\label{EP2_1D1-p}
			\\
			&|U_y(y,s)-\overline{U}'(y)| \le  \frac{1}{1+y^{2/9}}, \label{Utildey_M-p}
			\\
			& |U_{yy}(y,s)| \leq  \frac{M^{1/8} |y|}{(1+y^2)^{1/2}}, \label{EP2_1D3-p}
			\\
			&  |\partial_y ^3 U(0,s)-2^{25}| \leq 1, \label{EP2_1D2-p} 
			\\
			& \|\partial_y ^3 U(\cdot,s)\|_{L^{\infty}} \leq  2M^{3/4} ,\label{EP2_1D5-p}
			\\
			& \|\partial_y ^4 U(\cdot,s)\|_{L^{\infty}} \leq M. \label{EP2_1D4-p}
		\end{align}
	\end{subequations}
Then, for all $(y,s)\in\mathbb R\times[s_0,\sigma_1]$,
	\begin{subequations}\label{Boot_2-p'}
		\begin{align}
			&|U_y(y,s)-\overline{U}'(y)| \leq  \frac{ N(m_0)y^2}{2(1+y^2)} ,\label{EP2_1D1-p'}
			\\
			&|U_y(y,s)-\overline{U}'(y)| \le  \frac{3}{4(1+y^{2/9})}, \label{Utildey_M-p'}
			\\
			& |U_{yy}(y,s)| \leq  \frac{M^{1/8} |y|}{2(1+y^2)^{1/2}}, \label{EP2_1D3-p'}
			\\
			&  |\partial_y ^3 U(0,s)-2^{25}| \leq C\veps, \label{EP2_1D2-p'} 
			\\
			& \|\partial_y ^3 U(\cdot,s)\|_{L^{\infty}} \leq  \frac{5}{3}M^{3/4} ,\label{EP2_1D5-p'}
			\\
			& \|\partial_y ^4 U(\cdot,s)\|_{L^{\infty}} \leq \frac{M}{2}. \label{EP2_1D4-p'}
		\end{align}
	\end{subequations}
\end{proposition}

\begin{proof}
The proof is deferred to Section~\ref{sec4}, where it is decomposed into a sequence of lemmas. For the reader’s convenience, we indicate where each estimate is proved: \eqref{EP2_1D1-p'} is obtained in Lemma~\ref{Wy1_lem_p}, \eqref{Utildey_M-p'} in Lemma~\ref{mainprop_1-p}, \eqref{EP2_1D3-p'} in Lemma~\ref{Wy2_lem}, \eqref{EP2_1D2-p'} in Lemma~\ref{U30lem}, \eqref{EP2_1D5-p'} in Lemma~\ref{w-3-est},  and \eqref{EP2_1D4-p'} in Lemma~\ref{lem:Wy4}.
\end{proof}

We conclude this subsection by verifying that the initial data for $U$ satisfy the bootstrap assumptions.

\begin{remark}[Initial conditions for $U$]\label{init_rmk}
	The set of initial conditions \eqref{init_w_3-p}--\eqref{4.3a-p} implies that the self-similar unknown $U$ satisfies, at $s=s_0$,
	\begin{subequations}\label{EP-W-IC}
		\begin{align}
			& |U_y(y,s_0)-\overline U'(y)| \leq \frac{N(m_0)y^2}{4(1+y^2)},   \label{W-y2}
			\\
			&|U_y(y,s_0)-\overline{U}'(y)| \le \frac{5}{8(1+y^{2/9})},\label{Utildey_init}
			\\
			&\|U_{yy}(\cdot, s_0)\|_{L^{\infty}} \leq 2^{12}, \label{1D3}
			\\
			&|\partial_y ^3 U(0, s_0)-2^{25}| = 0,  \label{1D2} 
			\\
			&\|\partial_y ^3 U(\cdot, s_0) \|_{L^{\infty}} \leq 2^{26}, \label{1D5}
			\\
			&\|\partial_y ^4 U( \cdot, s_0) \|_{L^{\infty}} \leq 2^{39}. \label{1D4}
		\end{align}
	\end{subequations}	
	By Taylor's theorem, together with \eqref{constraint-p}, \eqref{1D2}, and \eqref{1D4},
	\begin{equation*}
		|U_{yy}(y,s_0)|\leq |y||\partial_y^3U(0,s_0)|+\frac{y^2}{2}\|\partial_y^4U(\cdot,s_0)\|_{L^{\infty}}\leq 2^{25}|y|+2^{38}y^2.
	\end{equation*}
	Combining this with \eqref{1D3}, we obtain
\begin{equation}\label{1D3'}
|U_{yy}(y,s_0)|\leq \min\left\{2^{25}|y|+2^{38}y^2,2^{12}\right\}\leq \frac{M^{1/8}|y|}{4(1+y^2)^{1/2}},\quad y\in\mathbb R,
\end{equation}
	provided $M$ is sufficiently large. Hence, from \eqref{EP-W-IC} and \eqref{1D3'}, the initial data for $U$ strictly satisfy the bootstrap assumptions \eqref{Boot_2-p}, provided $M$ is sufficiently large.
\end{remark}

\subsection{Modulation estimates}
We record the modulation estimates to be used in the proof of Theorem~\ref{mainthm}. These estimates are derived under the bootstrap assumptions, and their proof is given in Section~\ref{secPre}. We denote by $T_{\sigma_1}$ the time determined by $s(T_{\sigma_1})=\sigma_1$.

\begin{lemma}\label{lem:modulation-estimates}
Under the assumptions of Proposition~\ref{Boot}, there exists a constant $C>0$ such that
\begin{align}
&|\kappa(t)| \le C\veps^{5/6},\quad |\xi(t)|\le C\veps, \label{modul_dec} \\
&|\dot\tau(t)+2\kappa(t)| \le Ce^{-s}, \label{tau_kap} \\
&\left|(\kappa(t)+1)\kappa(t)-\dot\xi(t)\right| \le Ce^{-7s},  \label{kap_xi}
\end{align}
and
\begin{equation}\label{dottau}
|\dot\tau(t)|\le C\veps^{5/6},\quad 0\leq \frac{1}{1-\dot\tau(t)}\leq 1+\veps^{1/2}
\end{equation}
for all $t \in [-\varepsilon, T_{\sigma_1}]$.
\end{lemma}

\subsection{Proof of Theorem~\ref{mainthm}}\label{C13_subsec-p}
We split the proof into several steps.

\emph{Step 1.}
We prove the existence and regularity assertions using the self-similar formulation. By Lemma~\ref{local_exist}, the initial value problem \eqref{Nov} with initial data $v_0=u_0+\frac12$, where $u_0$ satisfies \eqref{in-H-C}--\eqref{H0_bound}, admits a unique local-in-time solution
\begin{equation*}
v\in C([-\veps,T];H^5(\mathbb R))\cap C^1([-\veps,T];H^4(\mathbb R))
\end{equation*}
for some $T>-\veps$. Hence the function $u$ defined by \eqref{u} solves \eqref{Nov-u}--\eqref{initu} on $[-\varepsilon,T]$. By the local well-posedness theory for ODEs, choosing $T>-\veps$ sufficiently close to $-\veps$ if necessary, the system \eqref{tau'}--\eqref{xi'} with initial data \eqref{init_modul} has a unique solution $(\tau,\kappa,\xi)$ on $[-\veps,T]$ satisfying $\tau(t)>t$. Therefore the change of variables \eqref{ys}--\eqref{UP} is well-defined. We denote by $U$ and $P$ the corresponding solution of \eqref{Nov-UP} on $[s_0,s(T)]$.

By Remark~\ref{init_rmk}, the initial data for $U$ strictly satisfy the bootstrap assumptions \eqref{Boot_2-p}. Hence, by continuity, these bounds hold on a nonempty interval. Proposition~\ref{Boot}, together with a standard continuation argument (see, e.g., Section~3.4 of \cite{KKY}), then implies that the self-similar solution $U$ is defined for all $s\ge s_0$ and satisfies \eqref{Boot_2-p}. The modulation estimates in Lemma~\ref{lem:modulation-estimates} hold along the same continuation.

We define $t=t(s)$ as the inverse of $s=s(t)$. Then it holds from \eqref{ys} that
\begin{equation*}
\frac{dt}{ds}=\frac{e^{-s}}{1-\dot\tau(t(s))}.
\end{equation*}
By \eqref{dottau}, we have $0\le (1-\dot\tau)^{-1}\le 1+\varepsilon^{1/2}$ for $\varepsilon>0$ sufficiently small. Hence
\begin{equation*}
t(s)=-\varepsilon+\int_{s_0}^s \frac{e^{-s'}}{1-\dot\tau(t(s'))}\,ds'
\end{equation*}
has a finite limit as $s\to\infty$. We denote this limit by $T_*$. Moreover, since
$\tau(t(s))-t(s)=e^{-s}\to0 \quad \text{as } s \to \infty$,
we have $\tau(t(s))\to T_*$ as $s\to\infty$. Thus, after continuous extension, $\tau(T_*)=T_*$. Therefore the above limit coincides with the time $T_*$ defined in \eqref{T-star}. In particular, $T_*<\infty$, and the solution $v$ exists on $[-\varepsilon,T_*)$.

We conclude this step by showing that $|T_*| \leq C \veps^{11/6}$. The bound in \eqref{dottau} gives
\begin{equation}\label{tau_small}
|\dot\tau(t)|\leq C\veps^{5/6},\quad t\in[-\veps,T_*).
\end{equation}
Using this together with $\frac{dt}{ds}=\frac{e^{-s}}{1-\dot\tau}$, we have
\begin{equation*}
|\tau(t(s))|\leq \int_{-\veps}^{t(s)}|\dot\tau(t')|\,dt'\leq C\veps^{5/6}\int_{s_0}^{s}\frac{e^{-s'}}{1-\dot\tau(t(s'))}\,ds'\leq C\veps^{5/6}\int_{s_0}^{s}e^{-s'}\,ds'\leq C\veps^{11/6}.
\end{equation*}
Letting $s\to\infty$, we have $|T_*| \leq C \veps^{11/6}$.

\emph{Step 2}: We show that $\textstyle \sup_{t<T_*} \left[ v(\cdot, t) \right]_{C^{7/9}(\Omega)}<\infty$ for any bounded open set $\Omega\subset \mathbb{R}$. We first prove the corresponding estimate for $u$ on an arbitrary bounded open set $\Omega_u\subset\mathbb R$. By \eqref{asymp-y-infty} and \eqref{Utildey_M-p}, we have
\begin{equation*}
|U_y(y,s)| \leq |\overline{U}'(y)| + |U_y(y,s)-\overline{U}'(y)| \leq \frac{C}{y^{2/9}+1}.
\end{equation*}
Thus,
\begin{equation}\label{Wydec_fin_2}
\sup_{y \in\mathbb{R}}  \left( (y^{2/9}+1)|U_y(y,s)| \right) \le C \quad  \text{ for all } s\ge s_0.
\end{equation} 
Using this, uniformly for $y,\wt{y}\in\mathbb{R}$ and $s\ge s_0$, we have
\begin{equation}\label{not-zero}
\frac{|U(y,s)-U(\wt{y},s)|}{{|y-\wt{y}|}^{7/9}}
= \frac{1}{{|y-\wt{y}|}^{7/9}}\left|\int_{\wt{y}}^y U_y(\hat{y},s)\,d\hat{y}\right|
\le \frac{C}{{|y-\wt{y}|}^{7/9}}\left|\int_{\wt{y}}^y (1+\hat{y}^2)^{-1/9}\,d\hat{y}\right|
\lesssim 1.
\end{equation}
Consider any two points $x\neq \wt{x}$ in $\Omega_u$. By the change of variables \eqref{ys}--\eqref{UP}, together with \eqref{not-zero}, we have 
\begin{equation*}
\frac{|u(x,t) - u(\tilde x, t ) | }{ | x - \tilde x|^{7/9}} = \frac{|U(y,s)-U(\wt{y},s)|}{{|y-\wt{y}|}^{7/9}} \lesssim 1
\end{equation*}
for all $t\in[-\ve, T_*)$. This yields
\begin{equation*}
\sup_{t<T_*} [u(\cdot, t)]_{C^{7/9}(\Omega_u)}<\infty
\end{equation*}
for any bounded open set $\Omega_u\subset\mathbb{R}$. We now transfer this estimate back to the original variable $v$. Let $\Omega\subset\mathbb R$ be a bounded open set. Since $T_*<\infty$ and $\Omega$ is bounded, there exists a bounded open interval $I_\Omega\subset\mathbb R$ such that $x-\frac14(t+\varepsilon)\in I_\Omega$ for all $x\in\Omega$ and all $t \in [-\varepsilon,T_*)$. By \eqref{u}, we have $v(x,t)=u(x-\frac14(t+\varepsilon),t)+\frac12$. Hence, for $x,\tilde x\in\Omega$ with $x\ne\tilde x$,
\begin{equation*}
\frac{|v(x,t)-v(\tilde x,t)|}{|x-\tilde x|^{7/9}}=\frac{|u(x-\frac14(t+\varepsilon),t)-u(\tilde x-\frac14(t+\varepsilon),t)|}{|(x-\frac14(t+\varepsilon))-(\tilde x-\frac14(t+\varepsilon))|^{7/9}}\le [u(\cdot,t)]_{C^{7/9}(I_\Omega)}.
\end{equation*}
Taking the supremum over $x,\tilde x\in\Omega$ and then over $t<T_*$, we obtain the desired bound
\begin{equation*}
\sup_{t<T_*}[v(\cdot,t)]_{C^{7/9}(\Omega)}\le \sup_{t<T_*}[u(\cdot,t)]_{C^{7/9}(I_\Omega)}<\infty.
\end{equation*}

\emph{Step 3}:
Note that, by \eqref{modul_dec} and \eqref{kap_xi}, $\dot\xi$ is uniformly bounded on $[-\varepsilon,T_*)$, which implies that the limit $\lim_{t\nearrow T_*}\xi(t)$ exists. We set
$x_\ast^u:=\lim_{t\nearrow T_*}\xi(t)$
and define the blow-up location of $v$ by
$x_\ast:=x_\ast^u+\frac14(T_*+\varepsilon)$.
In this step, we show that, for any $\alpha\in(\frac79,1]$ and any bounded open set $\Omega\subset\mathbb{R}$,
\begin{equation} \label{xcases}
\begin{cases}
\displaystyle \lim_{t\nearrow T_*}[v(\cdot,t)]_{C^\alpha(\Omega)}=\infty & \text{if } x_\ast\in\Omega, \\
\displaystyle \limsup_{t\nearrow T_*} \left[ v (\cdot, t) \right]_{C^{\alpha}(\Omega)} < \infty & \text{if } x_\ast\notin\overline{\Omega}.
\end{cases}
\end{equation}
To this end, it is enough to establish the corresponding estimates for \(u\), since \eqref{u} involves only a spatial translation and the addition of a constant, both of which preserve the H\"older seminorm.

For any $y,\widetilde y\in\mathbb{R}$ with $y\neq \widetilde y$, we set
\begin{equation*} 
x=\xi(t)+e^{-\frac92s}y,\quad \widetilde{x}=\xi(t)+e^{-\frac92s}\widetilde y.
\end{equation*}
Then, by \eqref{ys}--\eqref{UP},
\begin{equation} \label{scaleid}
\frac{|u(x,t)-u(\widetilde{x},t)|}{|x-\widetilde{x}|^\alpha} =e^{(\frac{9}{2}\alpha-\frac{7}{2})s} \frac{|U(y,s)-U(\widetilde y,s)|}{|y-\widetilde y|^{\alpha}}.
\end{equation}
To prove the blow-up of the H\"older seminorm in the \(u\)-variable, we apply this identity to a pair of points approaching the blow-up location \(x_\ast^u\). We take \(y=y_0\) and \(\wt{y}=0\), where \(y_0\in(-\eta_1,\eta_1)\setminus\{0\}\) and \(0<\eta_1\ll1\), and set
\begin{equation*}
x_0(t):=\xi(t)+e^{-\frac92s}y_0,\quad \wt{x}_0(t):=\xi(t).
\end{equation*}
Then \(x_0(t)\) and \(\wt{x}_0(t)\) converge to \(x_\ast^u\) as \(t\nearrow T_*\). By the mean value theorem,
\begin{equation*}
\frac{|U(y_0,s)-U(0,s)|}{|y_0|^{\alpha}}=|U_y(\overline{y},s)||y_0|^{1-\alpha}
\end{equation*}
for some \(\overline{y}\) lying between \(0\) and \(y_0\). From \eqref{constraint-p} and \eqref{EP2_1D1-p}, we see that $|U_y(\overline{y},s)| \geq \frac{1}{2}$ for sufficiently small $\eta_1 >0$, which implies that
\begin{equation} \label{Uybarlow}
|U_y(\overline{y},s)||y_0|^{1-\alpha}\ge c_0
\end{equation}
for some $c_0>0$.

Now suppose that $x_\ast\in\Omega$. We define the corresponding points in the original variables by
\begin{equation*}
X_0(t):=\frac14(t+\varepsilon)+x_0(t),\quad \wt{X}_0(t):=\frac14(t+\varepsilon)+\wt{x}_0(t).
\end{equation*}
Since \(X_0(t)\) and \(\wt{X}_0(t)\) converge to \(x_\ast\) as \(t\nearrow T_*\) and \(\Omega\) is open, both points belong to $\Omega$ for all $t$ sufficiently close to $T_*$. Therefore,
\begin{equation*}
[v(\cdot,t)]_{C^\alpha(\Omega)}
\ge \frac{|v(X_0(t),t)-v(\wt{X}_0(t),t)|}{|X_0(t)-\wt{X}_0(t)|^\alpha}
=\frac{|u(x_0(t),t)-u(\wt{x}_0(t),t)|}{|x_0(t)-\wt{x}_0(t)|^\alpha}.
\end{equation*}
Using \eqref{scaleid} and \eqref{Uybarlow}, we obtain
\begin{equation} \label{vHlow}
[v(\cdot,t)]_{C^\alpha(\Omega)}
\ge e^{(\frac{9}{2}\alpha-\frac{7}{2})s}
\frac{|U(y_0,s)-U(0,s)|}{|y_0|^\alpha}
\ge c_0e^{(\frac{9}{2}\alpha-\frac{7}{2})s}.
\end{equation}
Since $\alpha>\frac79$, we see that $[v(\cdot,t)]_{C^\alpha(\Omega)}\to\infty$ as $t\nearrow T_*$, i.e., $s\to\infty$.

We next prove the estimate for the second case in \eqref{xcases}. Set
\begin{equation*}
d_\ast:=d(x_\ast,\Omega)=\inf_{X\in\Omega}|X-x_\ast|>0.
\end{equation*}
Then we have, for all $t$ sufficiently close to $T_*$,
\begin{equation*}
\left|X-\frac14(t+\varepsilon)-\xi(t)\right|\ge \frac12d_\ast
\end{equation*}
for any $X\in\Omega$. Hence the corresponding self-similar variable
\begin{equation*}
Y_t(X):=e^{\frac92s}\left(X-\frac14(t+\varepsilon)-\xi(t)\right)
\end{equation*}
satisfies
\begin{equation*}
|Y_t(X)|\ge \frac12d_\ast e^{\frac92s}=:c_1e^{\frac92s}
\end{equation*}
for any $X\in\Omega$ and all $t$ sufficiently close to $T_*$. Now we fix $X,\widetilde X\in\Omega$ with $X\ne\widetilde X$, and define
\begin{equation*}
y:=Y_t(X),\quad \widetilde y:=Y_t(\widetilde X).
\end{equation*}
By \eqref{u} and the change of variables \eqref{ys}--\eqref{UP}, we have
\begin{equation*}
\frac{|v(X,t)-v(\widetilde X,t)|}{|X-\widetilde X|^\alpha}=e^{(\frac92\alpha-\frac72)s}\frac{|U(y,s)-U(\widetilde y,s)|}{|y-\widetilde y|^\alpha}.
\end{equation*}

We first consider the case where $y$ and $\widetilde y$ have the same sign. Then every $\widehat y$ between $y$ and $\widetilde y$ satisfies
$|\widehat y|\ge c_1e^{\frac92s}$.
This together with \eqref{Wydec_fin_2} gives
\begin{equation*}
|U_y(\widehat y,s)|\le C|\widehat y|^{-2/9}=C|\widehat y|^{\alpha-1}|\widehat y|^{\frac79-\alpha}\le C|\widehat y|^{\alpha-1}e^{(\frac72-\frac92\alpha)s}
\end{equation*}
for any $\alpha \in (\frac79,1]$. Therefore,
\begin{equation*}
\begin{split}
e^{(\frac92\alpha-\frac72)s}\frac{|U(y,s)-U(\widetilde y,s)|}{|y-\widetilde y|^\alpha} \le C\frac{\left|\int_{\widetilde y}^y |\widehat y|^{\alpha-1}\,d\widehat y\right|}{|y-\widetilde y|^\alpha} \le C\frac{\left||y|^\alpha-|\widetilde y|^\alpha\right|}{|y-\widetilde y|^\alpha} \le C.
\end{split}
\end{equation*}
We next consider $y$ and $\widetilde y$ with opposite signs. In this case, since
\begin{equation*}
\left| X-\frac14(t+\varepsilon)-\xi(t) \right| \ge \frac12d_\ast,\quad \left| \widetilde X-\frac14(t+\varepsilon)-\xi(t) \right| \ge \frac12d_\ast,
\end{equation*}
we have
$|X-\widetilde X|\ge d_\ast$.
Equivalently,
\begin{equation*}
|y-\widetilde y|=e^{\frac92s}|X-\widetilde X|\ge d_\ast e^{\frac92s}.
\end{equation*}
Using \eqref{not-zero}, we obtain
\begin{equation*}
\begin{split}
e^{(\frac92\alpha-\frac72)s}\frac{|U(y,s)-U(\widetilde y,s)|}{|y-\widetilde y|^\alpha} \le C e^{(\frac92\alpha-\frac72)s}|y-\widetilde y|^{\frac79-\alpha}  \le C e^{(\frac92\alpha-\frac72)s}\left(d_\ast e^{\frac92s}\right)^{\frac79-\alpha} \le C.
\end{split}
\end{equation*}
Combining the two cases yields
\begin{equation*}
\sup_{\substack{X,\widetilde X\in\Omega\\ X\ne\widetilde X}}\frac{|v(X,t)-v(\widetilde X,t)|}{|X-\widetilde X|^\alpha}\le C
\end{equation*}
for all $t$ sufficiently close to $T_*$. Hence
\begin{equation*}
\limsup_{t\nearrow T_*}\left[v(\cdot,t)\right]_{C^\alpha(\Omega)}<\infty
\end{equation*}
for any bounded open set $\Omega$ such that $x_\ast\notin\overline{\Omega}$ and any $\alpha \in (\frac79,1]$.

\emph{Step 4}: 
We finally prove the temporal blow-up rate \eqref{blow-up-rate}, i.e., for $x_\ast\in\Omega$ and $\alpha\in(\frac79,1]$,
\begin{equation*}
\left[v(\cdot,t)\right]_{C^\alpha(\Omega)}\sim (T_*-t)^{-\frac{9\alpha-7}{2}}
\end{equation*}
for all $t$ sufficiently close to $T_*$. We first prove the corresponding upper bound. By \eqref{Wydec_fin_2}, we have
\begin{equation*}
|U_y(y,s)|\le C(1+y^2)^{-1/9}\le C(1+y^2)^{(\alpha-1)/2}
\end{equation*}
for $\alpha\in(\frac79,1]$. Hence, for any $y,\widetilde y\in\mathbb R$,
\begin{equation*}
|U(y,s)-U(\widetilde y,s)|
\le \left|\int_{\widetilde y}^y |U_y(\widehat y,s)|\,d\widehat y\right|
\le C|y-\widetilde y|^\alpha.
\end{equation*}
Then from \eqref{ys} and \eqref{UP}, we have for any $x,\widetilde x\in\mathbb R$,
\begin{equation}\label{u_Calp2}
\frac{|u(x,t)-u(\widetilde x,t)|}{|x-\widetilde x|^\alpha}=e^{(\frac92\alpha-\frac72)s}\frac{|U(y,s)-U(\widetilde y,s)|}{|y-\widetilde y|^\alpha}\le Ce^{(\frac92\alpha-\frac72)s}.
\end{equation}
Thus, by \eqref{u} and \eqref{u_Calp2}, we obtain
\begin{equation*}
\left[v(\cdot,t)\right]_{C^\alpha(\Omega)}\le Ce^{(\frac92\alpha-\frac72)s}.
\end{equation*}
On the other hand, since $x_\ast\in\Omega$, the lower bound \eqref{vHlow} gives
\begin{equation*}
\left[v(\cdot,t)\right]_{C^\alpha(\Omega)}\ge ce^{(\frac92\alpha-\frac72)s}
\end{equation*}
for all $t$ sufficiently close to $T_*$. Consequently,
\begin{equation} \label{vasim}
\left[v(\cdot,t)\right]_{C^\alpha(\Omega)}\sim e^{(\frac92\alpha-\frac72)s}.
\end{equation}
Now it remains to rewrite this in terms of $T_*-t$. Note from \eqref{tau_small} that
\begin{equation*}
\frac12(T_*-t)\le \tau(t)-t=\int_t^{T_*}(1-\dot\tau(t'))\,dt'\le \frac32(T_*-t).
\end{equation*}
Since $e^{-s}=\tau(t)-t$, we obtain
\begin{equation*}
e^{(\frac92\alpha-\frac72)s}=(\tau(t)-t)^{-(\frac92\alpha-\frac72)}\sim (T_*-t)^{-(\frac92\alpha-\frac72)}=(T_*-t)^{-\frac{9\alpha-7}{2}}.
\end{equation*}
Combining this with \eqref{vasim} gives \eqref{blow-up-rate}. This completes the proof of Theorem~\ref{mainthm}.

\qed

\section{Closure of bootstrap assumptions}\label{sec4}

This section is devoted to proving the improved estimates \eqref{EP2_1D1-p'}--\eqref{EP2_1D4-p'} stated in Proposition~\ref{Boot}.

\subsection{Preliminary estimates} \label{secPre}

We provide some preliminary estimates and the proof of Lemma~\ref{lem:modulation-estimates}. We first present some inequalities that are used throughout our analysis. By \eqref{UP}, \eqref{u}, and \eqref{v_bound}, we have
\begin{equation}\label{u_bound}
|e^{-7s/2}U(y,s)+\kappa(t)|=\left|v\left(x+\frac{1}{4}(t+\varepsilon),t\right)-\frac{1}{2}\right|\leq \frac{9}{16} + \frac{1}{2} < C.
\end{equation} 
By \eqref{EP2_1D1-p}, \eqref{num_4}, and \eqref{Nm0}, we have a uniform bound on $U_y$ as  
\begin{equation} \label{Uy1-p}
	| U_y (y,s)| \le | \overline{U}'(y) | + | U_y(y,s) -\overline{U}'(y) | \le 4-\frac{3y^2}{1+y^2} +\frac{N(m_0) y^2}{1+y^2} \le 4.
\end{equation}
Also, by \eqref{Utildey_M-p} and \eqref{asymp-y-infty}, we get
\begin{equation}\label{Uy_M-p}
	|U_y(y,s)|\leq | \overline{U}'(y) | + | U_y(y,s) -\overline{U}'(y) |\leq \frac{C}{y^{2/9}+1},
\end{equation}
for some $C>0$.

Note that $\textstyle K(x):=\frac12 e^{-|x|}$ is the resolvent kernel for $(1-\partial_x^2)$ on $\mathbb{R}$; i.e., 
$( 1- \partial_x^2 ) K(x) = \delta_0(x)$,
where $\delta_0$ is the Dirac-delta measure concentrated at $x=0$. 
In view of \eqref{Nov2_2}, we can rewrite $p$ as 
\begin{equation*} 
	p(x,t) = \left(  K* \left(v^3+\frac{3}{2}vv_x^2+\frac{3}{2}v_x^2v_{xx}\right) \right) (x,t), 
\end{equation*}
where $*$ denotes convolution. 
Using \eqref{H}, \eqref{v_bound}, and \eqref{H0_bound}, we get 
\begin{equation*}
	\begin{split}
		|  p(x,t)  | &= \left|  K* \left(v^3 + \frac{3}{2}vv_x^2+\frac{3}{2}v_x^2v_{xx}\right) (x,t)  \right| 
		\\
		&\leq \frac{3}{4}\|v(\cdot,t)\|_{L^{\infty}}\int_{\mathbb{R}}\left(v^2+v_x^2\right)\,dx + \frac{1}{2}|K_x*v_x^3|
		\\
		&\leq C + \frac{1}{2}|K_x*v_x^3|.
	\end{split}
\end{equation*} 
Therefore, we have by \eqref{Uy1-p} and \eqref{H} that
\begin{equation}\label{p_bound}
	\begin{split}
		e^{-s}|P(y,s)|&\leq Ce^{-s}+C\|U_y(\cdot,s)\|_{L^{\infty}}\int_{\mathbb R}|K_x(x-z)|v_x^2(z,t)\,dz
		\\
		&\leq Ce^{-s}+C\int_{\mathbb{R}}v_x^2(z)\,dz
		\\
		&\leq Ce^{-s}+CH(-\veps) \leq C.
	\end{split}
\end{equation}

The following lemma gives exponential decay bounds for derivatives of $P$.

\begin{lemma} \label{P_high_lem}
Under the assumptions in Proposition~\ref{Boot}, there exists a constant $C>0$ such that
	\begin{equation} \label{P_high}
		\lVert \partial_y^j P(\cdot,s) \rVert_{L^\infty} \leq C e^{-3s/2}, \quad j=1,2,3
	\end{equation}
	for all $s \in [s_0,\sigma_1]$.
\end{lemma}

\begin{proof}
	Let $\lambda := e^{9s/2}$, and let $G_\lambda$ denote the Green function for the operator $(1-\lambda^2 \partial_{yy})$, i.e., $\textstyle G_\lambda(y) = \frac{1}{2\lambda} e^{-|y|/\lambda} $. It then follows that
	\begin{equation} \label{G'}
		G_\lambda'(y) = - \frac{1}{2\lambda^2} \frac{y}{|y|} e^{-|y|/\lambda}, \quad y \neq 0
	\end{equation}
	and
	\begin{equation} \label{G''}
		G_\lambda''(y) = \frac{1}{\lambda^2} G_\lambda(y) - \frac{1}{\lambda^2} \delta_0(y),
	\end{equation}
	where $'$ denotes $\textstyle \frac{d}{dy}$. Using $G_\lambda$, we can write the solution $P$ of \eqref{Peq} as
	\begin{equation*}
		P(y,s) = G_\lambda * \bigg[ \left(e^{-7s/2}U + \kappa + \frac{1}{2}\right)^3 + \frac{3}{2}e^{2s}U_y^2\left( e^{-7s/2}U + \kappa + \frac{1}{2} + e^{11s/2}U_{yy} \right) \bigg].
	\end{equation*}
Using \(U_y^2U_{yy}=\frac13\partial_y(U_y^3)\) and integrating by parts in the convolution, differentiating with respect to \(y\) gives
	\begin{equation*} 
		P_y(y,s) = G_\lambda' * \bigg[ \left(e^{-7s/2}U + \kappa + \frac{1}{2}\right)^3+\frac{3}{2}e^{2s}U_y^2\left( e^{-7s/2}U + \kappa + \frac{1}{2} \right) \bigg] + \frac{1}{2} e^{15s/2} G_\lambda'' * U_y^3.
	\end{equation*}
	We take the $L^\infty$-norm and use \eqref{G'}--\eqref{G''}, together with the bounds \eqref{u_bound} and \eqref{Uy1-p}, to obtain
	\begin{equation*}
		\begin{split}
			\left| P_y (y,s) \right| & \leq \bigg| \int_\mathbb{R} \frac{1}{2\lambda^2} e^{-|y-z|/\lambda} \left(e^{-7s/2}U + \kappa + \frac{1}{2}\right)^3 (z,s) \, dz \bigg| \\
			& \quad + \bigg| \int_\mathbb{R} \frac{3}{4\lambda^2} e^{-|y-z|/\lambda} e^{2s}U_y^2\left( e^{-7s/2}U + \kappa + \frac{1}{2} \right) (z,s) \, dz \bigg| \\
			& \quad + \bigg| \int_\mathbb{R} \frac{1}{4\lambda^3} e^{-|y-z|/\lambda} e^{15s/2} U_y^3(z,s) \, dz \bigg| + \bigg| \frac{e^{15s/2}}{2\lambda^2} U_y^3(y,s) \bigg| \\
			& \leq C \left( \frac{1}{\lambda} + \frac{e^{2s}}{\lambda} + \frac{e^{15s/2}}{\lambda^2} \right) \leq C e^{-3s/2},
		\end{split} 
	\end{equation*}
	which proves \eqref{P_high} for $j=1$.
	For the case $j=2$, we estimate \eqref{Peq} as
	\begin{equation*}
	\begin{split}
		|P_{yy}(y,s)| & = e^{-9s} \left| - P + \left(e^{-7s/2}U + \kappa + \frac{1}{2}\right)^3 +\frac{3}{2}e^{2s}U_y^2\left(\frac{1}{2}+e^{-7s/2}U + \kappa + e^{11s/2}U_{yy}\right) \right| \\
		& \leq C e^{-9s} \left( e^s + 1 + e^{2s} + e^{15s/2} \right) \leq C e^{-3s/2},
	\end{split}
	\end{equation*}
	where we used \eqref{p_bound}, \eqref{Uy1-p}, \eqref{u_bound}, and \eqref{EP2_1D3-p}. 
	Similarly, we obtain the bound on $|\partial_y^3P|$ by differentiating \eqref{Peq} and using in addition \eqref{EP2_1D5-p}:
	\begin{equation*}
	\begin{split}
		|\partial_y^3 P (y,s)| & = e^{-9s} \bigg| - P_y + 3 e^{-7s/2}U_y \left(e^{-7s/2}U + \kappa + \frac{1}{2}\right)^2 + 3 e^{2s} U_y U_{yy} \left(\frac{1}{2}+e^{-7s/2}U + \kappa + e^{11s/2}U_{yy}\right) \\
		& \qquad \quad + \frac{3}{2} e^{2s} U_y^2 \left( e^{-7s/2}U_y + e^{11s/2} \partial_y^3 U \right) \bigg| \\
		& \leq C e^{-9s} \left( e^{-s/4} + e^{-7s/2} + e^{2s} + e^{15s/2} \right) \leq C e^{-3s/2}.
	\end{split}
	\end{equation*}
	This completes the proof.
\end{proof}

Now we prove Lemma~\ref{lem:modulation-estimates}. The proof is first carried out on any interval on which the change of variables \eqref{ys} is well-defined and \(|\dot{\tau}|\le\frac12\). Such intervals are nonempty by local well-posedness, since \eqref{tau}, \eqref{init_modul}, and \eqref{p_bound} imply \(|\dot{\tau}(-\veps)|\le C\veps\). On such an interval \(s(t)\) is strictly increasing, so \(T_\sigma\) is well-defined whenever \(s(T_\sigma)=\sigma\). The estimate \eqref{dottau}, proved below, improves the auxiliary bound \(|\dot\tau|\le\frac12\). Therefore, after reducing \(\veps_0\) if necessary, a standard continuation argument extends the estimates up to \(\sigma=\sigma_1\).

\begin{proof}[Proof of Lemma~\ref{lem:modulation-estimates}]
From \eqref{kappa}, we have
\begin{equation*}
\begin{split}
|\dot{\kappa}| & \leq \frac{4e^{-9s/2}}{|\partial_y^3U(0,s)|} \left( |P_y(0,s)| + 12 e^{-7s/2} \left( \kappa + \frac{1}{2} \right)^2 + 32e^{-3s/2} \right)  + \left| e^{9s/2} P_y(0,s) - 32 e^{3s} \right| \\
&=:I_1+I_2.
\end{split}
\end{equation*}
First, by \eqref{P_high}, \eqref{u_bound}, and \eqref{EP2_1D2-p}, the term $I_1$ is bounded as
$I_1 \leq C e^{-6s}$.
Next we claim that
\begin{equation}\label{Py0}
I_2=e^{9s/2}\left|P_y(0,s) - 32 e^{-3s/2}\right|\leq Ce^{s/6}.
\end{equation}
As in the proof of Lemma~\ref{P_high_lem}, we write $P_y(y,s)$ using the Green function $\textstyle G_\lambda(y) := \frac{1}{2\lambda} e^{-|y|/\lambda}$ for the operator $(1-\lambda^2 \partial_{yy})$, with $\lambda=e^{9s/2}$. Evaluating the representation of $P_y$ at $y=0$, we obtain
\begin{equation}\label{Py0_1}
\begin{split}
P_y(0,s) - 32e^{-3s/2}  &=  \frac{1}{2\lambda^2} \int_{\mathbb{R}}\frac{z}{|z|} e^{-|z|/\lambda} \left(e^{-7s/2}U(z,s) + \kappa + \frac{1}{2}\right)^3 \,dz \\
&\quad + \frac{3e^{2s}}{4\lambda^2} \int_{\mathbb{R}} \frac{z}{|z|} e^{-|z|/\lambda} U_y^2(z,s)\left( e^{-7s/2}U(z,s) + \kappa + \frac{1}{2}  \right)\,dz \\
&\quad+   \frac{e^{15s/2}}{4\lambda^3} \int_{\mathbb{R}} e^{-|z|/\lambda}  U_y^3(z,s)\,dz=: J_1+J_2+J_3.
\end{split}
\end{equation}
Here the first term $J_1$ is estimated as
\begin{equation}\label{Py0_2}
\left|J_1\right|\leq C\lambda^{-1}=Ce^{-9s/2}
\end{equation}
using the bound \eqref{u_bound}. To estimate the next two terms, we observe that
\begin{equation}\label{**'}
e^{-5s/2}\int_{\mathbb{R}}|U_y(y,s)|^2\,dy=\int_{\mathbb{R}}|u_x(x,t)|^2\,dx \leq C
\end{equation}
for some constant $C>0$, where in the last inequality we used \eqref{H}. Using this together with \eqref{u_bound}, we have
\begin{equation}\label{Py0_3}
\left|J_2\right|\leq  Ce^{-7s} \int_{\mathbb{R}}   U_y^2(z,s)\,dz = Ce^{-9s/2} \int_{\mathbb{R}}   v_x^2\,dx \leq Ce^{-9s/2}.
\end{equation}
Lastly, we note that by \eqref{Uy_M-p}, \eqref{**'}, and H\"older's inequality, we have 
\begin{equation*}
\begin{split}
\int_{\mathbb{R}}|U_y(z,s)|^3\,dz &\leq C\int_{\mathbb{R}}\left(\frac{1}{z^{2/9}+1}\right)^{5/3}|U_y(z,s)|^{4/3}\,dz \\
&\leq C\left(\int_{\mathbb{R}}\left(\frac{1}{(z^{2/9}+1)^{5/3}}\right)^3\,dz\right)^{1/3}\left(\int_{\mathbb{R}}\left(|U_y(z,s)|^{4/3}\right)^{3/2}\,dz\right)^{2/3} \\
&\leq C\left(\int_{\mathbb{R}}|U_y(z,s)|^2\,dz\right)^{2/3}\leq Ce^{5s/3}.
\end{split}
\end{equation*}
Therefore, we obtain
\begin{equation}\label{Py0_4}
\left|J_3\right|\leq  Ce^{-6s} \int_{\mathbb{R}}   |U_y(z,s)|^3 \,dz  \leq Ce^{-13s/3}.
\end{equation}
Combining \eqref{Py0_1} with \eqref{Py0_2}, \eqref{Py0_3}, and \eqref{Py0_4}, we close the claim \eqref{Py0}. Hence,
\begin{equation*}
|\dot{\kappa}| \leq C e^{s/6}.
\end{equation*}
From this and \eqref{ys}, we have
\begin{equation}\label{kap_bound}
|\kappa(t)| \leq  \int_{-\varepsilon}^t |\dot{\kappa}| \, dt' \leq  \int_{s_0}^s |\dot{\kappa}| \frac{e^{-s'}}{1-\dot{\tau}} \, ds' \leq  \int_{s_0}^s C e^{-5s'/6} \, ds' \leq  C \varepsilon^{5/6}.
\end{equation}

Next, applying \eqref{P_high}, \eqref{u_bound}, \eqref{EP2_1D2-p}, and \eqref{kap_bound} to \eqref{xi}, we obtain
\begin{equation*}
\begin{split}
| \dot{\xi}| = \left| \frac{e^{-11s/2}}{\partial_y^3 U(0,s)} \left( P_y(0,s) + 12 e^{-7s/2} \left( \kappa + \frac{1}{2} \right)^2 - 32e^{-3s/2} \right) + (\kappa + 1)\kappa \right| \leq C.
\end{split}
\end{equation*}
Therefore, similarly to \eqref{kap_bound}, we obtain
$|\xi(t)| \leq C\veps$.
This concludes \eqref{modul_dec}.

		Applying \eqref{p_bound} and \eqref{modul_dec} to \eqref{tau}, we obtain \eqref{tau_kap}:
		\begin{equation*}
			|\dot{\tau}+2\kappa| \leq   \left|\frac{e^{-2s}P(0,s)}{4}\right| + \left|\frac{e^{-2s}}{4} \left( \kappa + \frac{1}{2} \right)^3\right| \leq Ce^{-s}.
		\end{equation*}
		Next we apply the bounds \eqref{EP2_1D2-p}, \eqref{P_high}, and \eqref{modul_dec} to \eqref{xi} to obtain
		\begin{equation*}
			\begin{split}
				\left|(\kappa + 1)\kappa-\dot{\xi}\right| &= \frac{e^{-11s/2}}{|\partial_y^3 U(0,s)|} \left| P_y(0,s) + 12 e^{-7s/2} \left( \kappa + \frac{1}{2} \right)^2 - 32e^{-3s/2} \right| \leq Ce^{-7s}.
			\end{split}
		\end{equation*}
		Finally, the bounds in \eqref{dottau} follow from \eqref{tau_kap} and \eqref{modul_dec} for sufficiently small $\veps>0$.
\end{proof}

We also record a far-field property of the transport velocity \(\mathcal U\) appearing in the differentiated equations \eqref{Uy_der1}--\eqref{Uy_der4}. The lemma below verifies the far-field hypothesis used in Section~\ref{proof_boot} when applying the decay estimate Lemma~\ref{rmk2}; the resulting far-field control is then used together with the maximum principle Lemma~\ref{max_2}.

\begin{lemma}\label{UW_far}
Under the assumptions in Proposition~\ref{Boot}, it holds that
\begin{equation} \label{eq:UW_far}
\inf_{\{|y|>1, s\in[s_0, \sigma_1]\}} \mathcal{U}(y,s) \frac{y}{|y|} \geq \frac{1}{8},
\end{equation}
where $\mathcal{U}(y,s)$ is defined by \eqref{U}.
\end{lemma}

\begin{proof}
We first observe that
\begin{equation}\label{tempvbound}
\left|e^{-7s/2}U(y,s)+\kappa+\frac12\right|=\left|v\left(x+\frac14(t+\varepsilon),t\right)\right|\le\frac{9}{16}
\end{equation}
by \eqref{v_bound} and \eqref{H0_bound}. Set \(A:=e^{-7s/2}U+\kappa+\frac12\). By \eqref{constraint-p} and \eqref{Uy1-p}, we have \(|U(y,s)|\le4|y|\). We recall the definition of \(\mathcal U\) from \eqref{U}:
\begin{equation*}
\mathcal U=\frac92y+\frac12U+AU+\frac{A\dot\tau U+(\frac12\dot\tau+\kappa)U}{1-\dot\tau}+e^{7s/2}\frac{(\kappa+1)\kappa-\dot\xi}{1-\dot\tau}.
\end{equation*}
Using \eqref{tempvbound}, we obtain
\begin{equation*}
\frac92y+\left(\frac12+A\right)U\ge\frac92y-\frac{17}{16}|U|\ge\frac14y,\quad y\ge0.
\end{equation*}
Moreover, by \eqref{kap_xi}, \eqref{dottau}, and \eqref{modul_dec},
\begin{equation*}
\left|\frac{A\dot\tau U+(\frac12\dot\tau+\kappa)U}{1-\dot\tau}\right|\le C\veps^{5/6}|y|,\quad \left|e^{7s/2}\frac{(\kappa+1)\kappa-\dot\xi}{1-\dot\tau}\right|\le Ce^{-7s/2}.
\end{equation*}
Hence, for \(y\ge1\),
\begin{equation} \label{calU_lower}
\mathcal U(y,s)\ge\left(\frac14-C\veps^{5/6}\right)y-Ce^{-7s/2}\ge\frac18y
\end{equation}
provided \(\veps>0\) is sufficiently small. Similarly, for \(y\le-1\), we have
\begin{equation} \label{calU_upper}
\mathcal U(y,s)\le\left(\frac14-C\veps^{5/6}\right)y+Ce^{-7s/2}\le\frac18y.
\end{equation}
The estimates \eqref{calU_lower}--\eqref{calU_upper} imply \eqref{eq:UW_far}.
\end{proof}

\subsection{Proof of \texorpdfstring{\eqref{EP2_1D1-p'}--\eqref{EP2_1D4-p'}}{(3.13a)--(3.13f)}} \label{proof_boot}

We prove the improved estimates in the order needed for the argument. We first close the higher-derivative bounds, then establish the weighted estimates for \(U_y-\overline U'\) and \(U_{yy}\). The final global weighted estimate \eqref{Utildey_M-p'} is obtained after the auxiliary far-field decay estimates required for the maximum-principle argument. 

In the proofs below, we first choose \(M>0\) sufficiently large. All subsequent smallness assumptions on \(\veps\) are allowed to depend on this choice of \(M\).

\begin{lemma}\label{U30lem} 
Under the assumptions in Proposition~\ref{Boot}, it holds that 
	\begin{equation}\label{str_U3-p}
		|\partial_y^3U(0,s)-2^{25}|\leq C\veps
	\end{equation}
	for all $s\in [s_0,\sigma_1]$.
\end{lemma}

\begin{proof}
	Evaluating \eqref{Uy_der3} at $y=0$ and using \eqref{constraint-p} and $\textstyle \frac{1}{1-\dot{\tau}} = 1 + \frac{\dot{\tau}}{1-\dot{\tau}}$, we obtain
	\begin{equation*}
		\partial_s \partial_y^3U(0, s) + D^U_3(s) \partial_y^3U(0, s) = F^U_3(s),
	\end{equation*} 
	where
		\begin{subequations}
			\begin{align*}
				D^U_3(s) &:=  - \frac{10 (2\kappa +\dot{\tau})}{1-\dot{\tau}},
				\\
				F^U_3(s) &:=- \frac{e^{-2s}P_{yy}(0,s)}{1-\dot{\tau}}  + \frac{96e^{-9s}}{1-\dot{\tau}} \left(\kappa + \frac{1}{2}\right) - \left(  e^{7s/2} \frac{(\kappa+1)\kappa - \dot{\xi}}{1-\dot{\tau}}\right) \partial_y^4 U(0,s) .
			\end{align*}
		\end{subequations}
	By \eqref{tau_kap} and \eqref{dottau}, $D^U_3$ is bounded as
		\begin{equation} \label{DU3}
			| D^U_3 (s) |\leq 10(1+\veps^{1/2})|\dot{\tau}+2\kappa|\leq  C e^{-s}. 
		\end{equation}
		Applying \eqref{constraint-p}, \eqref{p_bound}, and  \eqref{modul_dec} to \eqref{Peq}, we find
		\begin{equation*}
			e^{9s} |P_{yy}(0,s)| = \left|P(0,s) - \left( \kappa + \frac{1}{2}\right)^3-24e^{2s}\left(\kappa+\frac{1}{2}  \right)\right| \leq Ce^{2s}.
		\end{equation*}
		This, together with \eqref{EP2_1D4-p} and \eqref{kap_xi}, gives
			\begin{equation}\label{FU3}
				\begin{split}
					|F^U_3(s)|&=\left|- \frac{e^{-2s}P_{yy}(0,s)}{1-\dot{\tau}}  + \frac{96e^{-9s}}{1-\dot{\tau}} \left(\kappa + \frac{1}{2}\right) - \left(  e^{7s/2} \frac{(\kappa+1)\kappa - \dot{\xi}}{1-\dot{\tau}}\right) \partial_y^4 U(0,s) \right| \leq Ce^{-7s/2}. 
				\end{split}
			\end{equation}
		By \eqref{EP2_1D2-p}, \eqref{DU3}, and \eqref{FU3}, we have 
		\begin{equation*}
			\left|\partial_s \partial_y^3U(0,s)\right|\leq \left|D^U_3(s) \partial_y^3 U(0,s) \right|+\left|F^U_3 (s) \right|\leq  Ce^{-s}.
		\end{equation*} 
		Hence,
		\begin{equation*} 
			\left|\partial_y^3U (0, s)-2^{25}\right|  = \left|  \partial_y^3U (0, s) - \partial_y^3U (0, s_0 )  \right|  \leq \int_{s_0}^{s} \left|\partial_{s'} \partial_y^3 U (0, s' ) \right|\,ds' \leq C\int_{s_0}^s {e^{-s'}\,ds'}\leq C\veps,
		\end{equation*}
		which proves \eqref{str_U3-p}.		
\end{proof}

\begin{lemma}\label{lem:Wy4}
Under the assumptions in Proposition~\ref{Boot}, it holds that 
\begin{equation*}
\lVert \partial_y^4 U(\cdot,s) \rVert_{L^\infty} \leq \frac{M}{2}
\end{equation*}
for all $s \in [s_0,\sigma_1]$.

\end{lemma}

\begin{proof}
	From \eqref{Uy_der4}, we have
\begin{equation}\label{Uyyyy-eq}
	\partial_s \partial_y^4U + D^U_4 \partial_y^4U + \mathcal{U} \partial_y^5 U  = F^U_4, 
\end{equation}
where $\mathcal{U}$ is defined by \eqref{U}, and
	\begin{subequations}
		\begin{align*}
			D^U_4 (y,s) &:= \frac{29}{2}  + \frac{7(1+2\kappa)U_y}{2(1-\dot{\tau})} ,
			\\
			F^U_4 (y,s) &:= - \frac{e^{-2s}\partial_y^3 P}{1-\dot{\tau}} + \frac{3e^{-11s/2}\partial_y^3 U}{1-\dot{\tau}} \left(e^{-7s/2}U + \kappa + \frac{1}{2}\right)^2 + \frac{18e^{-9s}U_yU_{yy}}{1-\dot{\tau}} \left(e^{-7s/2}U + \kappa + \frac{1}{2}\right) \\
			& \qquad + \frac{6e^{-25s/2}U_y^3}{1-\dot{\tau}} - \left( \frac{19e^{-7s/2}U_y^2}{2(1-\dot{\tau})} + \frac{11(e^{-7s/2}U+\kappa) U_{yy}}{1-\dot{\tau}} + \frac{11U_{yy}}{2(1-\dot{\tau})} \right) \partial_y^3 U \\
			& \qquad - \frac{12e^{-7s/2}U_yU_{yy}^2}{1-\dot{\tau}}-\frac{7e^{-7s/2}U U_y}{1-\dot{\tau}}\partial_y^4U.
		\end{align*}
	\end{subequations}
	It follows from \eqref{dottau}, \eqref{modul_dec}, and \eqref{Uy1-p} that
\begin{equation}\label{DU4}
	D^U_4 (y,s) = \frac{29}{2}  + \frac{7(1+2\kappa)U_y}{2(1-\dot{\tau})} \geq \frac{29}{2}-\frac{7\cdot 4}{2}(1+\veps^{1/2})(1+C\veps^{5/6})\geq \frac{1}{4}
\end{equation}
for sufficiently small $\veps>0$.
	Also, by \eqref{dottau}, \eqref{P_high},  \eqref{EP2_1D5-p}, \eqref{Uy1-p}, \eqref{EP2_1D3-p}, and \eqref{EP2_1D4-p}, we have 
	\begin{equation*}
		\begin{split}
			|F^U_4(y,s)|&\leq  Ce^{-7s/2}+\left|e^{-7s/2}U+\kappa  + \frac{1}{2} \right| \frac{11|U_{yy}||\partial_y^3 U|}{1-\dot{\tau}} + \frac{7e^{-7s/2}|U| |U_y|}{1-\dot{\tau}}|\partial_y^4U|
			\\
			&\leq Ce^{-7s/2}+ \frac{99}{8}(1+\veps^{1/2})M^{7/8} + 7(1+\veps^{1/2})Me^{-7s/2}|U| |U_y|.
		\end{split}
	\end{equation*}
	Here, the last inequality used the fact that
	\begin{equation*}
		\left|e^{-7s/2}U+\kappa  + \frac{1}{2} \right|= \left| v\left(x+\frac{1}{4}(t+\varepsilon),t \right) \right|\leq  \frac{9}{16},
	\end{equation*}
	coming from \eqref{v_bound} and \eqref{H0_bound}.
	Note that \eqref{Uy_M-p}, \eqref{u_bound}, and \eqref{Uy1-p} imply
	\begin{equation*}
		e^{-7s/2}|U(y,s)||U_y(y,s)|\leq Ce^{-s/2}|e^{-7s/2}U|^{6/7} \frac{|U|^{1/7}}{y^{2/9}+1} \leq Ce^{-s/2}\frac{|y|^{1/7}}{y^{2/9}+1}\leq \veps^{1/4}.
	\end{equation*}
	Thus, we conclude that
		\begin{equation}\label{FU4}
		\begin{split}
			|F^U_4(y,s)|&\leq 13 M^{7/8} + 8M \veps^{1/4}. 
		\end{split}
	\end{equation}

Let  $\psi$ be the characteristic curve satisfying $\partial_s \psi=\mathcal{U} (\psi, s)$ and $\psi(y,s_0)=y$. 
By integrating \eqref{Uyyyy-eq} along $\psi$, and using \eqref{DU4}, \eqref{FU4}, and \eqref{1D4}, we have 
\begin{equation*}
\|\partial_y^4U(\cdot,s)\|_{L^{\infty}} \le 2^{39}e^{-(s-s_0)/4}+4\left(13M^{7/8}+8M\varepsilon^{1/4}\right)\leq \frac{M}{2}
\end{equation*}
where, in the last inequality, we have chosen $M>0$ sufficiently large and then $\varepsilon>0$ sufficiently small. This completes the proof.
\end{proof}

\begin{lemma}\label{Wy1_lem_p}
Under the assumptions in Proposition~\ref{Boot}, it holds that 
	\begin{equation}
		|U_y(y,s)-\overline{U}'(y)|\leq \frac{Ny^2}{2(1+y^2)} \label{Burgers_B.1-p}
	\end{equation}
	for all $y\in\mathbb{R}$ and $s\in[s_0,\sigma_1]$, where $N$ is the constant chosen in \eqref{Nm0}.
\end{lemma}

\begin{proof}
	Let $\wt{U} := U- \overline{U}$.  Then, from \eqref{Uy_der1} and \eqref{ovU}, we have the equation for $\wt{U}_y$: 
	\begin{equation}\label{Eq_diff-p}
		\begin{split}
			&\wt{U}_{ys}+\left(1+\frac{\wt{U}_y+2\overline{U}'}{4(1-\dot{\tau})}\right)\wt{U}_y+\mathcal{U}\wt{U}_{yy}=
			-\frac{\wt{U}\overline{U}''}{1-\dot{\tau}}-\frac{e^{-2s}P}{1-\dot{\tau}}+\frac{e^{-2s}}{1-\dot{\tau}}\left(e^{-7s/2}U+\kappa+\frac{1}{2}\right)^3
			\\
			&\quad-\left(\frac{(e^{-7s/2}U+2\kappa)U}{1-\dot{\tau}}+\frac{\dot{\tau}\overline{U}}{1-\dot{\tau}}+e^{7s/2}\frac{(\kappa+1)\kappa-\dot{\xi}}{1-\dot{\tau}}\right)\overline{U}''-\frac{\dot{\tau}(\overline{U}')^2}{4(1-\dot{\tau})}-\frac{(e^{-7s/2}U+\kappa)}{2(1-\dot{\tau})}U_y^2,
		\end{split}
	\end{equation}
	where $\mathcal{U}$ is defined by \eqref{U}.	Defining a new weighted function
	\[Z (y,s) :=\frac{y^2+1}{y^2}\wt{U}_y(y,s),\]
	we have from \eqref{Eq_diff-p} that 
	\begin{equation*} 
		\partial_s Z +D^Z (y,s) Z +\mathcal{U}(y,s) Z_y = F^Z(y,s)+\int_{\mathbb{R}} Z(y',s) K^Z(y,s;y') \,dy',
	\end{equation*}
	where
\begin{subequations}
\begin{align*}
D^Z(y,s) &:= 1+\frac{\wt{U}_y+2\overline{U}'}{4(1-\dot{\tau})}+\frac{2}{y(1+y^2)}\left(\frac{9}{2}y + \frac{(2\kappa+1)U}{1-\dot{\tau}} \right),\\
F^Z_1(y,s) &:=  -\frac{y^2+1}{y^2}\frac{e^{-2s}P}{1-\dot{\tau}}+\frac{e^{-2s}}{1-\dot{\tau}}\frac{y^2+1}{y^2}\left(e^{-7s/2}U+\kappa+\frac{1}{2}\right)^3
			\\
& \quad-\frac{y^2+1}{y^2}\left(\frac{(e^{-7s/2}U+2\kappa)U}{1-\dot{\tau}}+\frac{\dot{\tau}\overline{U}}{1-\dot{\tau}}+e^{7s/2}\frac{(\kappa+1)\kappa-\dot{\xi}}{1-\dot{\tau}}\right)\overline{U}'',
			\\
F^Z_2(y,s) &:= -\frac{y^2+1}{y^2}\frac{\dot{\tau}(\overline{U}')^2}{4(1-\dot{\tau})}-\frac{y^2+1}{y^2}\frac{(e^{-7s/2}U+\kappa)}{2(1-\dot{\tau})}U_y^2 \\
& \quad -\frac{2}{y(1+y^2)}\left(\frac{e^{-7s/2}U^2}{1-\dot{\tau}}+e^{7s/2}\frac{(\kappa+1)\kappa-\dot{\xi}}{1-\dot{\tau}}\right)Z,
\end{align*}
\end{subequations}
and
\begin{equation*}
K^Z(y,s;y') := -\frac{1}{1-\dot{\tau}}\frac{y^2+1}{y^2}\overline{U}''(y)\mathbb{I}_{[0,y]}(y')\frac{y'^2}{1+y'^2},
\end{equation*}
where \(\mathbb I_{[0,y]}\) is the oriented indicator of the interval from \(0\) to \(y\), namely
\begin{equation*}
\mathbb I_{[0,y]}(y')=
\begin{cases}
1, & 0\le y'\le y,\quad y>0,\\
-1, & y\le y'\le 0,\quad y<0,\\
0, & \text{otherwise}.
\end{cases}
\end{equation*}
First, we estimate $|Z|$ near $y=0$. By the Taylor expansion, 
	\begin{equation*}
		U_y(y,s)=-4+\frac{y^2}{2}\partial_y^3U(0,s)+\frac{y^3}{6}\partial_y^4U(y',s) \quad \text{ for some } |y'|<|y|<l,
	\end{equation*} 
	where $\textstyle l:= \frac{N}{M}>0$. 
	Using \eqref{EP2_1D4-p} and \eqref{str_U3-p}, we have
	\begin{equation*}
		\begin{split}
			| U_y(y, s) + 4 - 2^{24}y^2 | 
			& = \left| \frac{y^2}{2}\left(\partial_y^3U( 0, s) - 2^{25} \right)+\frac{y^3}{6}\partial_y^4U(y',s)\right|
			\leq C\veps  y^2 + \frac{M}{6}  |y|^3  \leq y^2\left( C \veps + \frac{ M | l | }{ 6  } \right)
		\end{split}
	\end{equation*}
	for all $| y | \le l$. 
	By this together with \eqref{U_taylor_small},  we have for sufficiently small $l=l(M)>0$, 
	\begin{equation*}
		\begin{split}
			|U_y(y, s) -\overline{U}'(y) |&\leq y^2\left(C\veps +\frac{M|l|}{6}\right)+C y^4
			\\
			&\leq y^2\left(C\veps +\frac{M|l|}{6}+C\cdot l^2\right)\leq \frac{ Ny^2}{4(1+y^2)} \quad \text{for}\quad |y|\leq l, 
		\end{split}
	\end{equation*} 
	which yields 
	\begin{equation}\label{V-l-p}
		|Z(y,s)|\leq \frac{N}{4}  \quad \text{for} \quad |y|\leq l.
\end{equation}
	
	We obtain a lower bound on $D^Z$. For this, we apply $U=\widetilde{U}+\overline{U}$, $\textstyle \frac{1}{1-\dot{\tau}} = 1 + \frac{\dot{\tau}}{1-\dot{\tau}}$, and \eqref{num_2} to the definition of $D^Z$:
	\begin{equation}\label{DV}
		\begin{split}
			D^Z(y,s) &= 1+\frac{\overline{U}'}{2}+\frac{2}{y(1+y^2)}\left(\frac{9y}{2}+\overline{U}\right) + \frac{\wt{U}_y+2\dot{\tau}\overline{U}'}{4(1-\dot{\tau})}+\frac{2(2\kappa+1)}{y(1+y^2)}\frac{\wt{U}+\dot{\tau}\overline{U}}{1-\dot{\tau}}
			\\
			&\geq \frac{y^2}{2(1+y^2)} -\left|\frac{\wt{U}_y+2\dot{\tau}\overline{U}'}{4(1-\dot{\tau})}+\frac{2(2\kappa+1)}{y(1+y^2)}\frac{\wt{U}+\dot{\tau}\overline{U}}{1-\dot{\tau}}\right|.
		\end{split}
	\end{equation}
	Using \eqref{EP2_1D1-p}, \eqref{Wbar_neg}, \eqref{U-rough}, \eqref{modul_dec}, and \eqref{dottau}, we see that 
	\begin{equation}\label{DV1}
		\begin{split}
			\left|\frac{\wt{U}_y+2\dot{\tau}\overline{U}'}{4(1-\dot{\tau})}+\frac{2(2\kappa+1)}{y(1+y^2)}\frac{\wt{U}+\dot{\tau}\overline{U}}{1-\dot{\tau}}\right|&\leq\left|\frac{\wt{U}_y}{4(1-\dot{\tau})}\right|+\left|\frac{2(2\kappa+1)}{y(1+y^2)}\frac{\wt{U}}{1-\dot{\tau}}\right|+\left|\frac{\dot{\tau}\overline{U}'}{2(1-\dot{\tau})}\right|+\left|\frac{2(2\kappa+1)}{1-\dot{\tau}}\frac{\dot{\tau}\overline{U}}{y}\right|
			\\
			&\leq \frac{ N(1+\veps^{1/2})y^2}{4(1+y^2)} + \left|\frac{2(1+\veps^{1/2})(1+C\veps^{5/6})}{y(1+y^2)}\right|\left| \int_0^{|y|} \frac{Ny'^2}{1+y'^2}\,dy'\right| +C\veps^{5/6}
			\\
			&\leq \frac{11}{12}\frac{Ny^2}{1+y^2}+C\veps^{1/2}\leq \frac{23}{24}\frac{Ny^2}{1+y^2}, \qquad |y|\geq l
		\end{split}
	\end{equation}
	for sufficiently small $\veps>0$.
	Here in the second inequality we have used the fact that 
	\begin{equation*}
		| \wt{U}(y,s) | \le \int_0^{|y|} |  \wt{U}_y (y', s) | dy'  \le \int_0^{|y|} \frac{Ny'^2}{1+y'^2 } dy',  
	\end{equation*}
	which follows from $\wt{U}(0,s)=0$ and \eqref{EP2_1D1-p}. 
	Hence, 
	\begin{equation}\label{damp-p}
		\begin{split}
			D^Z (y, s) &\geq \left(\frac{1}{2}-\frac{23N}{24}\right)\frac{ y^2}{1+y^2}\geq \frac{ y^2}{48(1+y^2)}, \qquad |y|\geq l,
		\end{split}
	\end{equation}
	as $\textstyle N(m_0)<\frac{1}{2}$ from \eqref{Nm0}.
	Next,  we compare $K^Z$ with $D^Z$. 
	From \eqref{num_6} and \eqref{DV}--\eqref{DV1}, we have 
	\begin{equation*}
		\begin{split}
			\int_\mathbb{R}{| K^Z (y,s;y')|\,dy'}&\leq |\overline{U}''(y)|\frac{y^2+1}{y^2}\int^{|y|}_0{\frac{y'^2}{1+y'^2 }\,dy'}
			\\
			&\leq  
			\delta \left(1+\frac{\overline{U}'}{2}+\frac{2}{y(1+y^2)}\left(\frac{9y}{2}+\overline{U}\right)\right)-\frac{\delta y^2}{48(1+y^2)}
			\\
			&\leq  \delta D^Z(y,s), \quad |y|\geq l 
		\end{split}
	\end{equation*}
	for some constant $\textstyle \delta\in(\frac{1}{2},1)$.

	We claim that 
	\begin{equation}\label{F_claim}
		\|F_1^Z(\cdot,s)\|_{L^{\infty}(|y|\geq l)}+\|F_2^Z(\cdot,s)\|_{L^{\infty}(|y|\geq l)}\leq C\veps^{1/2}.
	\end{equation}
	First we estimate $|F_1^Z|$ using \eqref{p_bound}, \eqref{u_bound}, \eqref{modul_dec}, \eqref{dottau}, \eqref{kap_xi}, \eqref{asymp-y-infty}, and $\textstyle \left|\frac{U+\overline{U}}{y}\right|\leq 8$ from \eqref{Uy1-p} and \eqref{U-rough}:
	\begin{equation*}
		\begin{split}
			\left|F^Z_1(y,s)\right|&\leq Ce^{-s}+e^{-7s/2}U^2(y,s)|\overline{U}''(y)|+C\veps^{5/6}\left|\frac{U+\overline{U}}{y}\right||y\overline{U}''(y)|+e^{-7s/2}|\overline{U}''(y)|
			\\
			&\leq C\veps^{5/6}+e^{-7s/2}U^2(y,s)|\overline{U}''(y)|, \qquad |y|\geq l.
		\end{split}
	\end{equation*}
	The second term on the right-hand side is bounded as
	\begin{equation*}
		\begin{split}
			e^{-7s/2}U^2(y,s)|\overline{U}''(y)|&\leq Ce^{-7s/2}(y^2+1)^{-11/18}U^2(y,s)
			\\
			&=Ce^{-s/2}|e^{-7s/2}U(y,s)|^{6/7}(y^2+1)^{-11/18}|U(y,s)|^{8/7}
			\\
			&\leq Ce^{-s/2}\frac{|y|^{8/7}}{(y^2+1)^{11/18}} \leq Ce^{-s/2}
		\end{split}
	\end{equation*}
	by \eqref{asymp-y-infty}, \eqref{Uy1-p}, and $|e^{-7s/2}U(y,s)|\leq C$ from \eqref{u_bound} and \eqref{modul_dec}.
	Therefore, we get $\|F_1^Z(\cdot,s)\|_{L^{\infty}(|y|\geq l)}\leq C\veps^{1/2}$.
	Similarly, by \eqref{dottau}, \eqref{modul_dec}, \eqref{Uy1-p}, \eqref{Wbar_neg}, \eqref{EP2_1D1-p}, and \eqref{kap_xi}, 
	\begin{equation*}
		\begin{split}
			|F^Z_2(y,s)|&\leq C|\dot{\tau}||\overline{U}'|^2+Ce^{-7s/2}|U|U_y^2+C|\kappa|U_y^2+Ce^{-7s/2}\frac{U^2}{|y|^3}\wt{U}_y+Ce^{7s/2}|(\kappa+1)\kappa-\dot{\xi}|
			\\
			&\leq C\veps^{5/6} + Ce^{-7s/2}|U|U_y^2, \qquad |y|\geq l.
		\end{split}
	\end{equation*}
	 By \eqref{Uy_M-p} and $|e^{-7s/2}U(y,s)|\leq C$, it holds that
	\begin{equation*}
		\begin{split}
			Ce^{-7s/2}|U(y,s)|U_y^2(y,s)&\leq Ce^{-7s/2}(y^2+1)^{-2/9}|U(y,s)|
			\\
			&=Ce^{-s/2}|e^{-7s/2}U|^{6/7}(y^2+1)^{-2/9}|U|^{1/7}
			\\
			&\leq Ce^{-s/2}\frac{|y|^{1/7}}{(y^2+1)^{2/9}}\leq Ce^{-s/2},
		\end{split}
	\end{equation*}
	where the second inequality holds by \eqref{Uy1-p}.
	Thus, $\|F_2^Z(\cdot,s)\|_{L^{\infty}(|y|\geq l)}\leq C\veps^{1/2}$. Combining this with the estimate for $F_1^Z$, we obtain \eqref{F_claim}.
	
		On the other hand, by \eqref{4.3a-p} and \eqref{Utildey_M-p}, we have 
		\begin{equation}\label{Wy_claim}
			|Z(y,s_0)|\leq \frac{N}{4}, \quad \limsup_{|y|\rightarrow \infty}|Z (y,s)|=0,
		\end{equation} 
		respectively. 
		Appealing to Lemma~\ref{max_2} in Appendix along with \eqref{V-l-p} and \eqref{damp-p}--\eqref{Wy_claim}, we conclude 
		$\|Z (\cdot,s)\|_{L^{\infty}(\mathbb{R})}\leq \frac{N}{2}$.
		We remark that the condition \eqref{D-cond} in Lemma~\ref{max_2} holds for sufficiently small $\ve>0$. This gives the desired estimate \eqref{Burgers_B.1-p}.	
\end{proof}

\begin{lemma}\label{Wy2_lem}
Under the assumptions in Proposition~\ref{Boot}, it holds that
		\begin{equation}\label{Uyy-bd}
		|U_{yy} (y, s) |\leq  \frac{M^{1/8}|y|}{2(1+y^2)^{1/2}}
	\end{equation}
		for all $y\in\mathbb{R}$ and $s \in [s_0,\sigma_1]$.
\end{lemma}

\begin{proof}
	Setting $\wt{Z}(y,s) := (y^2+1)^{1/2} y^{-1}   U_{yy}(y,s)$, 
	we obtain from \eqref{Uy_der2} that 
	\begin{equation*}
		\partial_s\wt{Z} + D^{\wt{Z}} \wt{Z} + \mathcal{U} \wt{Z}_y=F^{\wt{Z}},
	\end{equation*}
	where $\mathcal{U}$ is given in \eqref{U}, and
		\begin{subequations}
		\begin{align*}
			D^{\wt{Z}} (y, s) &:= \frac{11}{2} + \frac{3U_y}{2(1-\dot{\tau})}+\frac{3\kappa U_y}{1-\dot{\tau}}+\frac{1}{y(1+y^2)}\left(\frac{9}{2}y + \frac{(2\kappa + 1)U}{1-\dot{\tau}}  + e^{7s/2} \frac{(\kappa+1)\kappa - \dot{\xi}}{1-\dot{\tau}}\right),
			\\
			F^{\wt{Z}} (y, s) &:= -\frac{(y^2+1)^{1/2}}{y(1-\dot{\tau})}\left(  e^{-2s}P_y - 3e^{-11s/2}U_y \left(e^{-7s/2}U  + \kappa+\frac{1}{2}\right)^2  + \frac{e^{-7s/2}U_y^3}{2} + 3e^{-7s/2}UU_yU_{yy} \right)
			\\
			&\qquad-\frac{U_{yy}}{y^2(1+y^2)^{1/2}}\frac{e^{-7s/2}U^2}{1-\dot{\tau}}.
		\end{align*}
	\end{subequations}
		First, in the region $\textstyle \{ y : |y|\leq \frac{1}{M} =: l \}$ for sufficiently large $M>0$,  by the Taylor expansion together with \eqref{EP2_1D2-p} and \eqref{EP2_1D4-p}, we have 
	\begin{equation*}
		|U_{yy}(y,s)|\leq |y||\partial_y^3U(0,s)|+\frac{y^2}{2}|\partial_y^4U(y',s)|\leq (2^{25}+1)|y|+\frac{M}{2}y^2\leq \frac{M^{1/8} |y|}{4(1+y^2)^{1/2}}  \quad \text{for}\quad |y|\leq l,
	\end{equation*}
	where $y'\in (-| y |, | y |)$.
	This gives 
	\begin{equation}\label{Uyy_local}
		\|\wt{Z}(\cdot,s)\|_{L^{\infty}(|y|\leq l)}\leq  \frac{M^{1/8}}{4}.
	\end{equation}

	Next, we obtain a lower bound on $D^{\wt{Z}}$. We note that, by \eqref{num_3}, 
	\begin{equation*}
		\frac{11}{2}+\frac{3}{2}\overline{U}'+\frac{1}{y(1+y^2)}\left(\frac{9y}{2}+\overline{U}\right) \geq \frac{4y^2}{1+y^2}.
	\end{equation*}
	In light of this, we split  $D^{\wt{Z}}$ into two parts as 
	\begin{equation*}
		\begin{split}
			D^{\wt{Z}}&=\left(\frac{11}{2}+\frac{3}{2}\overline{U}'+\frac{1}{y(1+y^2)}\left(\frac{9y}{2}+\overline{U}\right)\right)
			\\
			&\quad +\left(\frac{3(\wt{U}_y+\dot{\tau}\overline{U}')}{2(1-\dot{\tau})}+\frac{3\kappa U_y}{1-\dot{\tau}}+\frac{1}{y(1+y^2)}\left(\wt{U}+\frac{(2\kappa+\dot{\tau})U}{1-\dot{\tau}}+e^{7s/2}\frac{(\kappa+1)\kappa-\dot{\xi}}{1-\dot{\tau}}\right)\right). 
		\end{split}
	\end{equation*}
	Using \eqref{EP2_1D1-p}, \eqref{Wbar_neg}, \eqref{U-rough}, \eqref{modul_dec}, \eqref{Uy1-p}, and \eqref{dottau}, we obtain
	\begin{equation*}
		\begin{split}
			&\left|\frac{3(\wt{U}_y+\dot{\tau}\overline{U}')}{2(1-\dot{\tau})}+\frac{3\kappa U_y}{1-\dot{\tau}}+\frac{1}{y(1+y^2)}\left(\wt{U}+\frac{(2\kappa+\dot{\tau})U}{1-\dot{\tau}}+e^{7s/2}\frac{(\kappa+1)\kappa-\dot{\xi}}{1-\dot{\tau}}\right)\right| 
			\\
			&\quad \leq \left|\frac{3(\wt{U}_y+\dot{\tau}\overline{U}')}{2(1-\dot{\tau})}\right|+\left|\frac{3\kappa U_y}{1-\dot{\tau}}\right|+\frac{|\wt{U}|}{|y|(1+y^2)}+\frac{1}{|y|(1+y^2)}\left|\frac{(2\kappa+\dot{\tau})U}{1-\dot{\tau}}+e^{7s/2}\frac{(\kappa+1)\kappa-\dot{\xi}}{1-\dot{\tau}}\right| 
			\\
			&\quad \leq \left|\frac{3N(1+\veps^{1/2})y^2}{2(1+y^2)}+C\veps^{5/6}\right|+C\veps^{5/6} + \frac{N}{|y|(1+y^2)}\left|\int^y_0\frac{ y'^2}{1+y'^2}\,dy'\right| + C\veps
			\\
			&\quad \leq \frac{11Ny^2}{6(1+y^2)} + C\veps^{5/6} \leq \frac{y^2}{1+y^2}, \qquad |y|\geq l
		\end{split}
	\end{equation*}
	for sufficiently small $\veps>0$.
	Here, the last inequality holds since \(N(m_0)<\frac12\) by \eqref{Nm0}.
	Thus, we have
	\begin{equation}\label{Uyy_D}
		\begin{split}
			D^{\wt{Z}}(y, s)	&\geq \frac{3y^2}{1+y^2}
			\quad \text{for } |y|\geq l.
		\end{split}
	\end{equation}

	We estimate $F^{\wt{Z}}$. By \eqref{dottau}, \eqref{u_bound}, \eqref{Uy1-p}, \eqref{P_high}, and \eqref{EP2_1D3-p}, it holds for $|y|\geq l$ that 
	\begin{equation*}
		\begin{split}
			|F^{\wt{Z}}(y,s)|&\leq Ce^{-2s}|P_y|+Ce^{-11s/2}|U_y|+Ce^{-7s/2}|U_y|^3+Ce^{-7s/2}|U||U_y||U_{yy}|+Ce^{-7s/2}\frac{|U_{yy}|U^2}{y^2}
			\\
			&\leq Ce^{-7s/2}+Ce^{-7s/2}|U||U_y|,
		\end{split}
	\end{equation*}
	where we also used $\textstyle |U(y,s)|\leq \left|\int^y_0 U_y(\wt{y},s)\,d\wt{y}\right|\leq 4|y|$ from \eqref{Uy1-p}. Furthermore,
	\begin{equation*}
		\begin{split}
			Ce^{-7s/2}|U(y,s)||U_y(y,s)|&\leq Ce^{-7s/2}(y^2+1)^{-1/9}|U(y,s)|
			\\
			&=Ce^{-s/2}|e^{-7s/2}U|^{6/7}(y^2+1)^{-1/9}|U|^{1/7}
			\\
			&\leq Ce^{-s/2}\frac{|y|^{1/7}}{(y^2+1)^{1/9}}\leq Ce^{-s/2}
		\end{split}
	\end{equation*}
	by \eqref{Uy_M-p} and \eqref{u_bound}. Thus we have
	\begin{equation}\label{Uyy_F}
		\|F^{\wt{Z}}(\cdot,s)\|_{L^{\infty}(|y|\geq l)}\leq Ce^{-s/2}.
	\end{equation}

		In addition, 
		\begin{equation*}
			|U_{yy}(y,s_0)|\leq \min\left\{2^{25}|y|+2^{38}y^2,2^{12}\right\}\leq \frac{M^{1/8}|y|}{4(1+y^2)^{1/2}},\quad y \in \mathbb{R}
		\end{equation*}
		by \eqref{init_w_3-p} and \eqref{init_24-p}. Therefore, we have  
	\begin{equation*}
		\|\wt{Z}(\cdot,s_0)\|_{L^{\infty}(\mathbb{R})}\leq \frac{M^{1/8}}{4}.
	\end{equation*}
Finally, we claim that $\textstyle\limsup_{|y|\rightarrow \infty}|\wt{Z} (y,s ) |\leq \frac14 M^{1/8}$ for sufficiently large $M>0$. 
Thanks to \eqref{Uyy_D} and \eqref{Uyy_F}, we have  
\begin{equation*}
\inf_{|y|\geq 10}D^{\wt{Z}}(y,s)\geq \frac{ 3\cdot 10^2}{1+ 10^2}=:\lambda_D, 
\end{equation*}
and
\begin{equation*}
\|F^{\wt{Z}}(\cdot,s)\|_{L^{\infty}(|y|\geq 10)}\leq Ce^{-s/2}. 
\end{equation*}
Applying Lemma~\ref{rmk2} with Lemma~\ref{UW_far}, we obtain the claimed inequality: 
\begin{equation}\label{Uyy_far}
\limsup_{|y|\rightarrow \infty}|\wt{Z}(y,s)|\leq \limsup_{|y|\rightarrow \infty}|\wt{Z}(y,s_0)|e^{-\lambda_D(s-s_0)}+  Ce^{-\lambda_Fs}  \le  \frac{M^{1/8}}{4},
\end{equation}
where we have used the fact that $\textstyle \|\widetilde Z(\cdot,s_0)\|_{L^\infty(|y|\ge10)}\leq 2^{12}\frac{\sqrt{101}}{10}<\frac14 M^{1/8}$, which follows from $\|U_{yy}(\cdot,s_0)\|_{L^\infty}\leq 2^{12}$.

We then apply Lemma~\ref{max_2} with \eqref{Uyy_local}-\eqref{Uyy_far} to conclude that $\textstyle \|\wt{Z} (\cdot, s) \|_{L^{\infty}}\leq \frac12 M^{1/8}$ for all $s\in[s_0,\sigma_1]$, which yields the desired estimate \eqref{Uyy-bd}. Here the condition \eqref{D-cond} in Lemma~\ref{max_2} follows by choosing $\varepsilon=\varepsilon(M)>0$ sufficiently small.
This completes the proof. 
\end{proof}

\begin{lemma}\label{w-3-est}
Under the assumptions in Proposition~\ref{Boot}, it holds that
\begin{equation*}
\|\partial_y^3U(\cdot,s)\|_{L^{\infty}}\leq \frac{5}{3}M^{3/4}
\end{equation*}
for all $s \in [s_0,\sigma_1]$.
\end{lemma}

\begin{proof}
	From \eqref{Uy_der3}, we have
	\begin{equation*}
		\partial_s \partial_y^3U + D^U_3 \partial_y^3U + 	\mathcal{U} \partial_y^4 U  =  F^U_3,
	\end{equation*}
	where $\mathcal{U}$ is defined by \eqref{U}, and 
		\begin{subequations}
			\begin{align*}
				D^U_3 &:= 10 +  \frac{5(1+2\kappa)U_y}{2(1-\dot{\tau})} ,
				\\
				F^U_3 &:=  - \frac{e^{-2s}P_{yy}}{1-\dot{\tau}} + \frac{3e^{-11s/2}U_{yy}}{1-\dot{\tau}} \left(e^{-7s/2}U + \kappa + \frac{1}{2}\right)^2 + \frac{6e^{-9s}U_y^2}{1-\dot{\tau}} \left(e^{-7s/2}U + \kappa + \frac{1}{2}\right) \nonumber\\
				& \qquad - \left( \frac{9e^{-7s/2}U_y^2}{2(1-\dot{\tau})} + \frac{3(e^{-7s/2}U +\kappa)U_{yy}}{1-\dot{\tau}} + \frac{3U_{yy}}{2(1-\dot{\tau})} \right) U_{yy}-\frac{5e^{-7s/2}U U_y \partial_y^3U}{1-\dot{\tau}}.
			\end{align*}
		\end{subequations}
	Due to \eqref{EP2_1D2-p} and \eqref{EP2_1D4-p}, it holds  for all  $\textstyle |y|\leq \frac{2}{3} M^{-1/4}$ that 
	\begin{equation} \label{Uy3loc}
		|\partial_y^3U(y,s)|\leq |\partial_y^3U(0,s)|+|y|\|\partial_y^4U(\cdot,s)\|_{L^{\infty}} \leq (2^{25}+1)+|y|M\leq \frac{5M^{3/4}}{6} 
	\end{equation}
	for sufficiently large $M>0$. 
	
	Now we consider the region $\textstyle |y|\geq \frac{2}{3} M^{-1/4}$. 
	Thanks to \eqref{num_4} and \eqref{EP2_1D1-p}, we have 
	\begin{equation*}
		10+\frac{5}{2}\overline{U}'+\frac{5}{2}\wt{U}_y \geq 10+\frac{5}{2}\overline{U}'-\frac{5N y^2}{2(1+y^2)} \geq
		\frac{25y^2}{4(1+y^2)},
	\end{equation*}
	as $\textstyle N(m_0)<\frac{1}{2}$ by \eqref{Nm0}. 
	Using this, together with \eqref{tau_kap}, \eqref{dottau}, and \eqref{Uy1-p}, we obtain 
	\begin{equation}\label{Uy3D'}
		\begin{split}
			D^U_3(y,s)&= \left(10+\frac{5}{2}\overline{U}'+\frac{5}{2}\wt{U}_y\right)+\frac{5(\dot{\tau}+2\kappa)U_y}{2(1-\dot{\tau})}
			\geq \frac{25y^2}{4(1+y^2)}-Ce^{-s}\geq \frac{8}{ 3M^{1/2}} 
		\end{split}
	\end{equation} 
	for all $\textstyle |y| \ge \frac{2}{3} M^{-1/4}$ and sufficiently small $\varepsilon=\varepsilon(M)>0$.

	We then show that \(F^U_3\) satisfies the following uniform bound:
	\begin{equation}\label{FU3'}
		\|F^U_3(\cdot,s)\|_{L^{\infty}(|y|\geq \frac{2}{3} M^{-1/4})}\leq \frac{5}{2}M^{1/4} .
	\end{equation}
	By \eqref{P_high}, \eqref{dottau}, \eqref{u_bound}, \eqref{Uy1-p}, and \eqref{EP2_1D3-p}, it holds that
	\begin{equation}\label{FU3'-1}
	\begin{split}
		& \left|- \frac{e^{-2s}P_{yy}}{1-\dot{\tau}} + \frac{3e^{-11s/2}U_{yy}}{1-\dot{\tau}} \left(e^{-7s/2}U + \kappa + \frac{1}{2}\right)^2 + \frac{6e^{-9s}U_y^2}{1-\dot{\tau}} \left(e^{-7s/2}U + \kappa + \frac{1}{2}\right)-\frac{9e^{-7s/2}U_y^2U_{yy}}{2(1-\dot{\tau})}\right| \\
		& \quad \leq C M^{1/8} e^{-7s/2}.
	\end{split}
	\end{equation}
	Also, by \eqref{dottau}, \eqref{u_bound}, \eqref{EP2_1D3-p}, \eqref{v_bound}, and \eqref{H0_bound}, we have
	\begin{equation} 
		\left|  \frac{3(e^{-7s/2}U+\kappa)U_{yy}^2}{1-\dot{\tau}} + \frac{3U_{yy}^2}{2(1-\dot{\tau})} \right| = \left| \frac{3v\left(x+\frac14(t+\varepsilon),t\right)}{1-\dot{\tau}} U_{yy}^2 \right| \leq  \frac{27}{16}(1+\varepsilon^{1/2})M^{1/4}.
	\end{equation}
To estimate the last term of \(F_3^U\), we use \eqref{Uy_M-p}, the bound \(|e^{-7s/2}U(y,s)|\le C\) following from \eqref{u_bound} and \eqref{modul_dec}, and \(|U(y,s)|\le4|y|\), which follows from \eqref{constraint-p} and \eqref{Uy1-p}:
		\begin{equation*}
		\begin{split}
			e^{-7s/2}|U(y,s)||U_y(y,s)|&\leq e^{-7s/2}(y^2+1)^{-1/9}|U(y,s)|
			\\
			&=Ce^{-s/2}|e^{-7s/2}U|^{6/7}(y^2+1)^{-1/9}|U|^{1/7}
			\\
			&\leq Ce^{-s/2}\frac{|y|^{1/7}}{(y^2+1)^{1/9}}\leq Ce^{-s/2},
		\end{split}
	\end{equation*}
	which, together with \eqref{dottau} and \eqref{EP2_1D5-p}, yields
	\begin{equation}\label{FU3'-3}
		\left|\frac{5e^{-7s/2}U U_y \partial_y^3U}{1-\dot{\tau}}\right|\leq C M^{3/4} \veps^{1/2}.
	\end{equation}
	Combining \eqref{FU3'-1}--\eqref{FU3'-3}, we get \eqref{FU3'} for sufficiently small $\varepsilon>0$.

	Next, we claim that $\textstyle\limsup_{|y|\rightarrow \infty}|\partial_y^3U(y,s)|\leq \frac56 M^{3/4} $. 
	Due to \eqref{Uy3D'} and \eqref{FU3'}, 
	it holds that for $|y|\geq 10$, 
	\begin{subequations}
		\begin{align*}
			&\inf_{|y|\geq 10}D^U_3(y,s)\geq \frac{ 25\cdot 10^2}{4(1+10^2)}-C\veps =:\lambda_{D_3^U} >0, 
			\\
			&\|F^U_3(\cdot,s)\|_{L^{\infty}(|y|\geq 10)}\leq CM^{1/4}.
		\end{align*}
	\end{subequations}
	Having these and Lemma~\ref{UW_far}, we apply Lemma~\ref{rmk2} to get
	\begin{equation}\label{Uy3dec} 
		\limsup_{| y |\rightarrow \infty}|\partial_y^3U( y,s)|  \leq
		\limsup_{|y|\rightarrow \infty}|\partial_y^3U(y,s_0)|e^{-  \lambda_{ D_3^U} (s-s_0)}+ CM^{1/4} ( \lambda_{ D_3^U})^{-1} \leq \frac{5M^{3/4}}{6}
	\end{equation}
	for some constant $M>0$ sufficiently large. 
	Here we have used the initial condition \eqref{init_24-p}, which yields $\textstyle \|\partial_y^3U(\cdot,s_0)\|_{L^{\infty}}\leq 2^{26}<\frac{1}{2}M^{3/4}$.
	
	We are ready to complete the proof by applying Lemma~\ref{max_2} with \eqref{init_24-p} and \eqref{Uy3loc}--\eqref{Uy3dec}. Since $K \equiv 0$, the condition \eqref{max_2_2} holds for any $\delta \in (0,1)$; we fix $\delta = 10^{-3}$. With $\textstyle d_0= \frac{5}{6}M^{3/4}$, $\textstyle \lambda_D = \frac{8}{3} M^{-1/2}$, and $\textstyle F_0 = \frac{5}{2} M^{1/4}$, we have
	\begin{equation*}
	d_0 \lambda_D = \frac{20}{9}M^{1/4} > \frac{5 }{4(1-10^{-3}) }M^{1/4} =  \frac{F_0}{2(1-\delta)},
	\end{equation*}
hence \eqref{D-cond} holds. Therefore, by Lemma~\ref{max_2}, we conclude that
	$\|\partial_y^3U(\cdot,s)\|_{L^{\infty}}\leq \frac{5M^{3/4}}{3}$.
\end{proof}

The next three lemmas are devoted to the proof of \eqref{Utildey_M-p'}. First, Lemma~\ref{Pdec_lem} gives a weighted estimate for the nonlocal forcing term appearing in the equation for the weighted derivative. Lemma~\ref{dec-lem} then uses this estimate, together with the far-field behavior of the transport velocity, to obtain the boundary control at spatial infinity. Finally, Lemma~\ref{mainprop_1-p} combines this far-field control with a maximum-principle argument to close \eqref{Utildey_M-p'}.

\begin{lemma} \label{Pdec_lem}
Under the assumptions in Proposition~\ref{Boot}, it holds that 
\begin{equation}\label{P-dec}
\bigg\|e^{-2s}(y^{2/9}+1)\bigg(P-\left(e^{-7s/2}U + \kappa + \frac{1}{2}\right)^3\bigg)\bigg\|_{L^{\infty}_y } \leq Ce^{-5s/6}
\end{equation} 
for all $s\in[s_0,\sigma_1]$.
\end{lemma}
	
\begin{proof}
	From \eqref{p_bound} and \eqref{u_bound}, we first observe that
	\begin{equation*}
	\bigg\lVert e^{-2s} \bigg(P-\left(e^{-7s/2}U + \kappa + \frac{1}{2}\right)^3\bigg)  \bigg\rVert_{L^\infty} \leq C e^{-s}.
	\end{equation*}
	Thus, for the proof of this lemma, it is enough to show that
	\begin{equation}\label{P-dec'}
		y^{2/9}\bigg|P-\left(e^{-7s/2}U+\kappa+\frac{1}{2}\right)^3\bigg| \leq Ce^{7s/6}.
	\end{equation}
	
	Note that the resolvent kernel $K_a(x)$ for the operator $(a^2 - \partial_x^2)$ with $a>0$ is given by $\textstyle K_a (x)=\frac{1}{2a}e^{-a|x|}$. Let $a=e^{-9s/2}$. We then have from \eqref{Peq} that
		\begin{equation*}
			\begin{split}
				P(y,s)&=\frac{e^{-9s/2}}{2}\int_{\mathbb{R}}e^{-e^{-9s/2}|y-z|} \\
				& \qquad \qquad \quad \times \bigg[ \left(e^{-7s/2}U + \kappa + \frac{1}{2}\right)^3 + \frac{3}{2}e^{2s}U_y^2\left(\frac{1}{2}+e^{-7s/2}U + \kappa + e^{11s/2}U_{yy}\right)\bigg](z,s)\,dz
				\\
				&=\frac{e^{-9s/2}}{2}\int_{\mathbb{R}}e^{-e^{-9s/2}|y-z|}\bigg[ \left(e^{-7s/2}U + \kappa +  \frac{1}{2}\right)^3+\frac{3}{2}e^{2s}U_y^2\left(\frac{1}{2}+e^{-7s/2}U + \kappa \right)\bigg](z,s)\,dz 
				\\
				&\quad -\frac{e^{-3s/2}}{4}\int_{\mathbb{R}} \operatorname{sgn}(y-z)e^{-e^{-9s/2}|y-z|}U_y^3(z,s)\,dz,
			\end{split}
		\end{equation*}
		where $\operatorname{sgn}(y)$ is defined by $1$ if $y>0$,  $-1$ if $y<0$.
		Using this representation of $P$, together with the fact $\textstyle \frac{a}{2}\int_{\mathbb{R}} e^{-a|x|} dx =1$, one can bound the quantity on the left-hand side of \eqref{P-dec'} as
		\begin{equation} \label{P-dec''}
			\begin{split}
				y^{2/9}&\bigg|P-\left(e^{-7s/2}U+\kappa+\frac{1}{2}\right)^3\bigg|
				\\
				&\leq \frac{e^{-9s/2}}{2}y^{2/9}\int_{\mathbb{R}}e^{-e^{-9s/2}|y-z|}\bigg|\left(e^{-7s/2}U+\kappa+\frac{1}{2}\right)^3(z,s)-\left(e^{-7s/2}U+\kappa+\frac{1}{2}\right)^3(y,s)\bigg|\,dz 
				\\
				&\quad+ \frac{3e^{-5s/2}}{4}y^{2/9}\int_{\mathbb{R}}e^{-e^{-9s/2}|y-z|}|U_y(z,s)|^2\left|\frac{1}{2}+e^{-7s/2}U(z,s)+\kappa\right|\,dz 
				\\
				&\quad+ \frac{e^{-3s/2}}{4}y^{2/9}\int_{\mathbb{R}}e^{-e^{-9s/2}|y-z|}|U_y(z,s)|^3\,dz
				\\
				&\leq
				Ce^{-8s}y^{2/9}\int_{A}e^{-e^{-9s/2}|y-z|}\left|U(z,s)-U(y,s)\right|\,dz 
				+Ce^{-9s/2}y^{2/9}\int_{A^c}e^{-e^{-9s/2}|y-z|}\,dz 
				\\
				&\quad+ Ce^{-5s/2}y^{2/9}\int_{\mathbb{R}}e^{-e^{-9s/2}|y-z|}|U_y(z,s)|^2\,dz 
				\\
				&\quad+ Ce^{-3s/2}y^{2/9}\int_{\mathbb{R}}e^{-e^{-9s/2}|y-z|}|U_y(z,s)|^3\,dz
				\\
				&=:I_{A} + I_{A^c} + J_1 + J_2,
			\end{split}
		\end{equation}
		where $\textstyle A:=\{z\in\mathbb{R} \;|\; |y-z|<\frac{|y|}{2}\}$. In the second inequality, we have used the following inequalities, which hold true due to \eqref{u_bound}: 
		\begin{equation*}
			\left|\left(e^{-7s/2}U+\kappa+\frac{1}{2}\right)^3(z,s)-\left(e^{-7s/2}U+\kappa + \frac{1}{2}\right)^3(y,s)\right|\leq Ce^{-7s/2}\left|U(z,s)-U(y,s)\right|
		\end{equation*} for $z\in A$, and 
		\begin{equation*}
			\left|\left(e^{-7s/2}U+\kappa+\frac{1}{2}\right)^3(z,s)-\left(e^{-7s/2}U+\kappa+\frac{1}{2}\right)^3(y,s)\right|\leq C
		\end{equation*}
		for $z\in A^c$.

		Now we estimate the terms $I_A$, $I_{A^c}$, $J_1$, and $J_2$. First, using the bound \eqref{Uy_M-p}, we obtain
		\begin{equation*}
			\begin{split}
				I_{A}
				&\leq Ce^{-8s}\int_{A}e^{-e^{-9s/2}|y-z|}y^{2/9}|y-z|\left|U_y(\wt{z},s)\right|\,dz 
				\\
				&\leq Ce^{-8s}\int_{A}e^{-e^{-9s/2}|y-z|}\frac{y^{2/9}}{\wt{z}^{2/9}+1}|y-z|\,dz 
			\end{split}
		\end{equation*}
		for some $\wt{z}$ between $y$ and $z$.
		Note that if $z\in A$ , then $\textstyle |\wt{z}-y|< \frac{|y|}{2}$, which yields $\textstyle \frac{|y|}{2} < |\wt{z}|<3\frac{|y|}{2}$. Hence,		
		\begin{equation}\label{I01}
			\begin{split}
				I_{A}&\leq Ce^{-8s}\int_{A}e^{-e^{-9s/2}|y-z|}\frac{y^{2/9}}{(\frac{|y|}{2})^{2/9}+1}|y-z|\,dz 
				\\
				&\leq Ce^{-8s}\int_{A}e^{-e^{-9s/2}|y-z|}|y-z|\,dz
				\\
				&\leq Ce^{s},
			\end{split}
		\end{equation} 
		where we used the identity $\textstyle \int_{\mathbb{R}}e^{-a|y-z|}|y-z|\,dz=2a^{-2}$.
		Next we estimate $I_{A^c}$ as follows:
		\begin{equation}\label{I02}
			\begin{split}
				I_{A^c}
				&\leq Ce^{-9s/2}\int_{A^c}e^{-e^{-9s/2}|y-z|/2}\left(e^{-e^{-9s/2}|y|/4}y^{2/9}\right)\,dz
				\\
				&\leq Ce^{-7s/2}\int_{A^c}e^{-e^{-9s/2}|y-z|/2}\,dz 
				\\
				&\leq Ce^s.
			\end{split}
		\end{equation}
		Here we used $\textstyle e^{-a |y | /4}y^{2/9}\leq Ce^{s}\left(e^{-a | y | /4}(a | y | )^{2/9}\right) \leq Ce^s$ and $\textstyle \int_{\mathbb{R}}e^{-a|y-z|/2}\,dz=\frac{4}{a}$.
		Similarly,
		\begin{equation} \label{J_1}
			\begin{split}
				J_1&=Ce^{-5s/2}y^{2/9}\int_{A}e^{-e^{-9s/2}|y-z|}|U_y(z,s)|^2\,dz+Ce^{-5s/2}y^{2/9}\int_{A^c}e^{-e^{-9s/2}|y-z|}|U_y(z,s)|^2\,dz 
				\\
				&\leq Ce^{-5s/2}\int_{A}e^{-e^{-9s/2}|y-z|}\frac{y^{2/9}}{(\frac{y}{2})^{2/9}+1}|U_y(z,s)|\,dz
				\\
				&\quad+Ce^{-5s/2}\int_{A^c}e^{-e^{-9s/2}|y-z|/2}\left(e^{-e^{-9s/2}|y|/4}y^{2/9}\right)|U_y(z,s)|^2\,dz  
				\\
				&\leq Ce^{-5s/2}\int_{A}e^{-e^{-9s/2}|y-z|}|U_y(z,s)|\,dz+Ce^{-3s/2}\int_{A^c}e^{-e^{-9s/2}|y-z|/2}|U_y(z,s)|^2\,dz 
				\\
				&\leq Ce^{-5s/2}\left(\int_{\mathbb{R}}e^{-2e^{-9s/2}|y-z|}\,dz\right)^{1/2}\left(\int_{\mathbb{R}}|U_y(z,s)|^2\,dz\right)^{1/2}+Ce^{-3s/2}\int_{A^c}|U_y(z,s)|^2\,dz 
				\\
				&= Ce^{-s/4}\left(\int_{\mathbb{R}}|U_y(z,s)|^2\,dz\right)^{1/2}+Ce^{-3s/2}\int_{A^c}|U_y(z,s)|^2\,dz \leq Ce^s.
			\end{split}
		\end{equation}
		In the last inequality, we used
		\begin{equation}\label{**}
			e^{-5s/2}\int_{\mathbb{R}}|U_y(y,s)|^2\,dy=\int_{\mathbb{R}}|u_x(x,t)|^2\,dx \leq C,
		\end{equation}
		which follows from \eqref{H}. Finally, we estimate $J_2$. As above, using $\textstyle e^{-e^{-9s/2} |y | /4}y^{2/9}\leq Ce^s$, we obtain
		\begin{equation}\label{J2}
			\begin{split}
				J_2&=Ce^{-3s/2}y^{2/9}\int_{A}e^{-e^{-9s/2}|y-z|}|U_y(z,s)|^3\,dz + Ce^{-3s/2}y^{2/9}\int_{A^c}e^{-e^{-9s/2}|y-z|}|U_y(z,s)|^3\,dz
				\\
				&\leq e^{-3s/2}\int_A e^{-e^{-9s/2}|y-z|}\frac{y^{2/9}}{(y/2)^{2/9}+1}|U_y(z,s)|^2\,dz 
				\\
				&\quad+ e^{-3s/2}\int_{A^c}e^{-e^{-9s/2}|y-z|}|y-z|^{2/9}|U_y(z,s)|^3\,dz
				\\
				&\leq Ce^{-3s/2}\int_{\mathbb{R}}|U_y(z,s)|^2\,dz + e^{-s/2}\int_{\mathbb{R}}|U_y(z,s)|^3\,dz.
			\end{split}
		\end{equation}
		By \eqref{Uy_M-p} and \eqref{**}, we have 
		\begin{equation*}
			\begin{split}
				e^{-s/2}\int_{\mathbb{R}}|U_y(z,s)|^3\,dz &\leq Ce^{-s/2}\int_{\mathbb{R}}\left(\frac{1}{z^{2/9}+1}\right)^{5/3}|U_y(z,s)|^{4/3}\,dz 
				\\
				&\leq Ce^{-s/2}\left(\int_{\mathbb{R}}\left(\frac{1}{(z^{2/9}+1)^{5/3}}\right)^3\,dz\right)^{1/3}\left(\int_{\mathbb{R}}\left(|U_y(z,s)|^{4/3}\right)^{3/2}\,dz\right)^{2/3}
				\\
				&\leq Ce^{-s/2}\left(\int_{\mathbb{R}}|U_y(z,s)|^2\,dz\right)^{2/3}\leq Ce^{7s/6}.
			\end{split}
		\end{equation*}
		Therefore, applying this and \eqref{**} to \eqref{J2}, we get
		\begin{equation} \label{J_2}
			J_2\leq Ce^{7s/6}.
		\end{equation}
		Combining the estimates \eqref{I01}, \eqref{I02}, \eqref{J_1}, and \eqref{J_2} with \eqref{P-dec''}, we deduce \eqref{P-dec'}, and hence \eqref{P-dec}. 
\end{proof}

\begin{lemma}\label{dec-lem}
Under the assumptions in Proposition~\ref{Boot}, it holds that 
\begin{equation}\label{Ut_dec_p}
\limsup_{|y|\rightarrow \infty} |(y^{2/9}+1)(U_y (y, s) -\overline{U}' (y) ) | <  \frac{5}{8}
\end{equation}
for all $s\in[s_0,\sigma_1]$.
\end{lemma}

\begin{proof}
	We note that 
	\begin{equation*}
		\begin{split}
			\limsup_{|y|\rightarrow \infty}|(y^{{2/9}}+1)(U_y-\overline{U}')|&\leq \limsup_{|y|\rightarrow \infty}|(y^{2/9}+1)U_y|+\lim_{|y|\rightarrow \infty}|(y^{2/9}+1)\overline{U}'|
			\\
			&= \limsup_{|y|\rightarrow \infty}|(y^{2/9}+1)U_y|+\Theta
		\end{split}
	\end{equation*}
	for $\textstyle\Theta=\left(\frac{1}{1296}\right)^{1/9}<\frac{1}{2}$, 
	where we used \eqref{asymp-y-infty} with $\beta=1$. 
	Therefore, it suffices to show that
	\begin{equation*} 
		\limsup_{|y|\rightarrow \infty} |(y^{2/9}+1)U_y(y,s)|< \frac{5}{8}-\Theta. 
	\end{equation*}
	
	Let $\mu(y,s):=(y^{2/9}+1)U_y(y,s)$. Then we have from \eqref{Uy_der1} that
	\begin{equation*}
		\partial_s \mu+D^{\mu}\mu+\mathcal{U} \partial_y \mu =F^{\mu},
	\end{equation*}
	where $\mathcal{U}$ is defined by \eqref{U}, and
	\begin{subequations}
		\begin{align*}
			D^{\mu}(y,s) &:= 1  + \frac{U_y}{4} -\frac{2y^{2/9}}{9(y^{2/9}+1)}\left(\frac{9}{2}+\frac{U}{y}\right),
			\\
			F^{\mu}(y,s) &:=  -(y^{2/9}+1)\frac{e^{-2s}}{1-\dot{\tau}} \left(P-\left(e^{-7s/2}U + \kappa + \frac{1}{2}\right)^3\right) - (y^{2/9}+1)\frac{e^{-7s/2}U U_y^2}{2(1-\dot{\tau})}
			\\
			&\qquad+\frac{2y^{2/9}}{9(y^{2/9}+1)}\frac{1}{y}\left(\frac{(e^{-7s/2}U + 2\kappa + \dot{\tau})U}{1-\dot{\tau}} + e^{7s/2} \frac{(\kappa+1)\kappa - \dot{\xi}}{1-\dot{\tau}}\right)\mu-\frac{(\dot{\tau}+2\kappa)U_y}{4(1-\dot{\tau})}\mu.
		\end{align*}
	\end{subequations}
	Using \eqref{Utildey_M-p} and $\wt{U}(0,s)=0$, where $\wt{U}(y,s):=U(y,s)-\overline{U}(y)$, we have
	\begin{equation*}
		\begin{split}
			D^{\mu}(y,s)&=1+\frac{\wt{U}_y}{4}+\frac{\overline{U}'}{4}-\frac{2y^{2/9}}{9(1+y^{2/9})}\left(\frac{9}{2}+\frac{\wt{U}+\overline{U}}{y}\right)
			\\
			&\geq 1-\frac{1}{4(1+y^{2/9})}+\frac{\overline{U}'}{4}-\frac{2y^{2/9}}{9(1+y^{2/9})}\left(\frac{9}{2}+\frac{\overline U}{y}+\frac{1}{y}\int_{0}^{y}{\frac{dy'}{1+{y'}^{2/9}}}\right), \quad |y| \geq m_0,
		\end{split}
	\end{equation*}
	where $m_0>0$ is the sufficiently large constant from \eqref{Utildey_M-p}.
	Thanks to \eqref{num_1} and $\overline{U}'\leq 0$ from \eqref{Wbar_neg}, we see that $D^\mu$ is non-negative for $|y|\geq m_0$, i.e., \begin{equation}\label{Dmu}
		\inf_{|y|\geq m_0}D^{\mu} (y,s) \geq 0.
	\end{equation} 
	
	Next we show that 
	\begin{equation}\label{F-mu}
		\|F^\mu(\cdot,s)\|_{L^{\infty}( |y|\geq m_0 )}  \leq C(m_0)e^{-s/2}.
	\end{equation}
	By \eqref{P-dec}, \eqref{Utildey_M-p}, \eqref{dottau}, \eqref{Uy1-p}, \eqref{tau_kap}, and \eqref{kap_xi}, we have
	\begin{equation}
		\begin{split}
			\left|F^\mu(y,s)\right|&\leq  Ce^{-5s/6} +Ce^{-7s/2}(y^{2/9}+1)|U| U_y^2 +C\left(e^{-7s/2}\frac{U^2}{|y|} +  e^{-s} \right), \quad |y|\geq m_0.
		\end{split}
	\end{equation}
	We first claim
	\begin{equation}\label{Fnu_1}
		\left\|e^{-7s/2}(y^{2/9}+1)UU_y^2\right\|_{L^{\infty}_y}\leq Ce^{-s/2}.
	\end{equation}
	By \eqref{u_bound} and \eqref{modul_dec}, we have $|e^{-7s/2}U(y,s)|\leq C$. Moreover, \eqref{Uy_M-p} gives
	\begin{equation} 
		|U(y,s)|=\left|\int^y_0 U_y(\wt{y},s)\,d\wt{y}\right|\leq C\int^{|y|}_0\frac{1}{\wt{y}^{2/9}+1}\,d\wt{y}\leq C|y|^{7/9}.
	\end{equation}
	Using these bounds, we estimate
	\begin{equation*}
		\begin{split}
			|e^{-7s/2}(y^{2/9}+1)UU_y^2|&\leq C|e^{-7s/2}(y^{2/9}+1)^{-1}U(y,s)|
			\\
			&=Ce^{-s/2}|e^{-7s/2}U|^{6/7}(y^{2/9}+1)^{-1}|U|^{1/7} 
			\\
			&\leq C e^{-s/2}\frac{|y|^{1/9}}{(y^{2/9}+1)}\leq Ce^{-s/2}.
		\end{split}
	\end{equation*}
	This proves \eqref{Fnu_1}. 
	Similarly, we obtain
	\begin{equation}\label{Fnu_2}
		\left|e^{-7s/2}\frac{U^2}{y}\right|=e^{-s/2}|e^{-7s/2}U|^{6/7}\left|\frac{U^{8/7}}{y} \right|\leq Ce^{-s/2} \frac{|y|^{8/9}}{|y|} \leq Ce^{-s/2} \quad \text{for }|y|\geq m_0.
	\end{equation}
	Therefore, by \eqref{P-dec}, \eqref{Fnu_1}, and \eqref{Fnu_2}, we get \eqref{F-mu}.
	
	Now we apply Lemma~\ref{rmk2} along with Lemma~\ref{UW_far}, \eqref{Dmu}, and \eqref{F-mu} to obtain
	\begin{equation*}
		\limsup_{|y|\rightarrow \infty}|\mu(y,s)|\leq \limsup_{|y|\rightarrow \infty}|\mu(y,s_0)|+C\veps^{1/2}  < \frac{5}{8}-\Theta,
	\end{equation*}
	for sufficiently small $\varepsilon>0$.
	Here, we used \eqref{4.3a-p2} and \eqref{Uy_M-p}, that is equivalent to 
	\begin{equation*}
		\limsup_{|y|\rightarrow \infty}|\mu(y,s_0)|\leq \frac{1}{2}\left(\frac{5}{8}-\Theta\right).
	\end{equation*}
	This completes the proof. 
\end{proof}

\begin{lemma}\label{mainprop_1-p} 
Under the assumptions in Proposition~\ref{Boot}, it holds that 
\begin{equation*}
(y^{2/9}+1) | {U}_y (y, s) - \UU'(y) | \leq \frac{3}{4}
\end{equation*}
for all $y\in\mathbb{R}$ and $s\in[s_0,\sigma_1]$.
\end{lemma}

\begin{proof}
	Let $\nu(y,s):=(y^{2/9}+1)\wt{U}_y(y,s)$, where $\wt{U}(y,s):=U(y,s)-\overline{U}(y)$. Then, by \eqref{Uy_der1} and \eqref{ovU}, $\nu$ satisfies
	\begin{equation}\label{nu_eq1-p}
		\begin{split}
			&\nu_s+D^\nu(y,s)\nu + \mathcal{U}\nu_y
			=\int_{\mathbb{R}}\nu(y',s)K^\nu(y,s;y')\,dy'+F_1^\nu(y,s)+F_2^\nu(y,s),
		\end{split}
	\end{equation}
	where $\mathcal{U}$ is defined by \eqref{U}, and
	\begin{subequations}
		\begin{align*}
			D^\nu(y,s)&:=1+\frac{\wt{U}_y+2\overline{U}'}{4}-\frac{2y^{2/9}}{9(y^{2/9}+1)}\left(\frac{9}{2} + \frac{U}{y} \right),
			\\
			F^\nu_1(y,s)&:=-(y^{2/9}+1)\frac{e^{-2s}}{1-\dot{\tau}}\left(P-\left(e^{-7s/2}U+\kappa+\frac{1}{2}\right)^3\right)-(y^{2/9}+1)\frac{e^{-7s/2}UU_y^2}{2(1-\dot{\tau})}
			\\
			&\quad+\frac{2y^{2/9}}{9(y^{2/9}+1)}\frac{1}{y}\left( \frac{(e^{-7s/2}U + 2\kappa + \dot{\tau})U}{1-\dot{\tau}} + e^{7s/2} \frac{(\kappa+1)\kappa - \dot{\xi}}{1-\dot{\tau}}\right)\nu-\frac{(\dot{\tau}+2\kappa)(y^{2/9}+1)U_y^2}{4(1-\dot{\tau})},
			\\
			F^\nu_2(y,s)&:=
			-\left(\frac{e^{-7s/2}U^2}{1-\dot{\tau}}+\frac{(2\kappa+\dot{\tau})\overline{U}}{1-\dot{\tau}}+e^{7s/2}\frac{(\kappa+1)\kappa-\dot{\xi}}{1-\dot{\tau}}\right)(y^{2/9}+1)\overline{U}'',
		\end{align*}
	\end{subequations}
	and
\begin{equation*}
			K^\nu(y,s;y'):=-\frac{1+2\kappa}{1-\dot{\tau}}\overline{U}''(y)\mathbb{I}_{[0,y]}(y')\frac{y^{2/9}+1}{y'^{2/9}+1},
\end{equation*}
	where \(\mathbb I_{[0,y]}\) is the oriented indicator of the interval from \(0\) to \(y\).
	We first claim that
	\begin{equation}\label{D-K-p}
		D^\nu(y,s)\geq \int_{\mathbb{R}}|K^\nu(y,s;y')|\,dy', \qquad  |y|\geq m_0
	\end{equation}
	for some constant $m_0>0$.
	By the bootstrap assumption \eqref{Utildey_M-p}, we have
	\begin{equation}\label{nu_D-p}
		\begin{split}
			D^\nu(y,s)&=1+\frac{1}{4}\wt{U}_y+\frac{1}{2}\overline{U}'-\frac{2y^{2/9}}{9(y^{2/9}+1)}\left(\frac{9}{2}+\frac{U}{y}\right)
			\\
			&\geq 1-\frac{1}{4(y^{2/9}+1)}+\frac{\overline{U}'}{2}-\frac{2y^{2/9}}{9(y^{2/9}+1)}\left(\frac{9}{2}+\frac{\overline U}{y}+\frac{1}{y}\int_{0}^{y}{\frac{dy'}{y'^{2/9}+1}}\right) 
			=: D^\nu_{-}(y),
		\end{split}
	\end{equation}  
	Moreover, by \eqref{dottau}, \eqref{modul_dec},
	\begin{equation}\label{nu_K-p}
		\begin{split}
			\int_{\mathbb{R}}{|K^\nu(y, s; y')|\,dy'}
			&\leq (1+C\veps^{5/6})(y^{2/9}+1)|\overline{U}''(y)|\int^{|y|}_{0}{\frac{dy'}{(y'^{2/9}+1)}} =: K^\nu_{+}(y).
		\end{split}
	\end{equation} 
	Then it follows from \eqref{num_1} that
	\begin{equation}\label{D-K+}
		D^\nu_-(y) \ge  K^\nu_+(y), \quad  |y|\geq m_0
	\end{equation} 
	for sufficiently small $\veps>0$ and sufficiently large $m_0>0$, which implies \eqref{D-K-p}.

	Next we show that 
	\begin{equation}\label{F-nu-p}
		\|F^\nu_1(\cdot,s)\|_{L^{\infty}( |y|\geq m_0 )} + \|F^\nu_2(\cdot,s)\|_{L^{\infty}( |y|\geq m_0 )}  \leq C(m_0)e^{-s/2}.
	\end{equation}
	Here $F^{\nu}_1(y,s)$ has the same algebraic structure as $F^{\mu}$ in the proof of Lemma~\ref{dec-lem}, and hence $F^\nu_1$ satisfies a similar decay estimate $\lVert F_1^{\nu}(\cdot,s) \rVert_{L^\infty(|y| \geq m_0)} \leq C(m_0)e^{-s/2}$ as shown in \eqref{F-mu}--\eqref{Fnu_2}. We omit the details to avoid repetition. Therefore, it is enough to show $\textstyle \|F^\nu_2(\cdot,s)\|_{L^{\infty}( |y|\geq m_0 )}  \leq C(m_0)e^{-s/2}$. By \eqref{dottau}, \eqref{tau_kap}, \eqref{kap_xi}, and \eqref{U-rough}, we have
	\begin{equation*}
		\begin{split}
			|F^\nu_2(y,s)|&\leq C\left(\frac{e^{-7s/2}U^2}{|y|}+\frac{|2\kappa+\dot{\tau}||\overline{U}|}{|y|}+e^{7s/2}\frac{|(\kappa+1)\kappa-\dot{\xi}|}{|y|}\right)|y|(y^{2/9}+1)|\overline{U}''|
			\\
			&\leq C\left(\frac{e^{-7s/2}U^2}{|y|}+Ce^{-s}\right)|y|(y^{2/9}+1)|\overline{U}''|, \qquad |y|\geq m_0.
		\end{split}
	\end{equation*}
	As in the proof of Lemma~\ref{dec-lem}, using \eqref{u_bound}, \eqref{modul_dec}, and \eqref{Uy_M-p}, we obtain $\textstyle \frac{e^{-7s/2}U^2}{|y|}\leq Ce^{-s/2}$. Moreover, by \eqref{asymp-y-infty} we have $|y|(y^{2/9}+1)|\overline{U}''(y)|\le C$ for some constant $C>0$. Plugging these bounds into the estimate above yields \eqref{F-nu-p}.
	
	To conclude the proof, we claim that
	\begin{equation}\label{claim_reg-p}
		\|\nu(\cdot,s)\|_{L^{\infty}}\leq \frac{3}{4}, \qquad s\in[s_0,\sigma_1].
	\end{equation}
	Suppose to the contrary that \eqref{claim_reg-p} fails. 
	Thanks to the fact that $\nu \in C([s_0, \sigma_1]\times \mathbb{R})$ and the initial condition \eqref{4.3a-p}, i.e., $\textstyle \|\nu(\cdot,s_0)\|_{L^{\infty}}\leq \frac{5}{8}$,  
	we have that 
	\begin{equation}\label{s2-def}
		s_2:= \min \left\{ s \in [s_0, \sigma_1]: \|\nu(\cdot,s)\|_{L^{\infty}} = \frac{3}{4} \right\}
	\end{equation} 
	is well-defined, and that there exists $s_1\in [s_0, s_2)$ such that
	\begin{equation*} 
		\frac{5}{8} = \|\nu(\cdot,s_1)\|_{L^{\infty}}  \le \|\nu(\cdot,s)\|_{L^{\infty}} < \frac{3}{4} \quad \text{ for all } s\in(s_1,s_2). 
	\end{equation*}
	Then, for each \(s\in[s_1,s_2]\), we have \(\|\nu(\cdot,s)\|_{L^\infty}\ge \frac{5}{8}\), and Lemma~\ref{dec-lem} implies that
\begin{equation*}
\limsup_{|y|\to\infty}|\nu(y,s)|<\frac58\le \|\nu(\cdot,s)\|_{L^\infty}.
\end{equation*}
Hence \(|\nu(\cdot,s)|\) attains its maximum at some finite point. Moreover, by \eqref{EP2_1D1-p} and \eqref{Nm0},
\begin{equation*}
\|\nu(\cdot,s)\|_{L^\infty(|y|\le m_0)}
\le (m_0^{2/9}+1)\|\widetilde U_y(\cdot,s)\|_{L^\infty(|y|\le m_0)}
\le N(m_0)(m_0^{2/9}+1)<\frac12.
\end{equation*}
Since \(\|\nu(\cdot,s)\|_{L^\infty}\ge\frac58\), any maximum point lies in \(\{|y|>m_0\}\). In particular, \(\nu\) is smooth at such a point and its \(y\)-derivative vanishes there. We may therefore choose
\begin{equation} \label{y-star}
y_\ast(s)\in\{ |y|>m_0:\ |\nu(y_\ast(s),s)|=\|\nu(\cdot,s)\|_{L^\infty}\},
\end{equation}
with \(\partial_y\nu(y_\ast(s),s)=0\).
	By \eqref{D-K-p}, \eqref{nu_D-p}, \eqref{nu_K-p}, \eqref{D-K+}, and \eqref{y-star},  
	if $\nu(y_*(s),s)> 0$, then
	\begin{equation}
		\begin{split}\label{nu_DK-p}
			D^\nu(y_*(s),s)\nu(y_*(s),s)&\geq D^\nu_{-}(y_*(s))\|\nu(\cdot,s)\|_{L^{\infty}}
			\\
			&\geq K^\nu_{+}(y_*(s))\|\nu(\cdot,s)\|_{L^{\infty}}\geq \left|\int_{\mathbb{R}}K^\nu(y_*(s),s;y')\nu(y',s)\,dy'\right|.
		\end{split}
	\end{equation}
	Similarly, if  $\nu(y_*(s),s) < 0$, we have 
	\begin{equation}
		\begin{split}\label{nu_DK2-p}
			D^\nu(y_*(s),s)\nu(y_*(s),s)&\le - D^\nu_{-}(y_*(s))\|\nu(\cdot,s)\|_{L^{\infty}}
			\\
			&\le - K^\nu_{+}(y_*(s))\|\nu(\cdot,s)\|_{L^{\infty}} \le - \left|\int_{\mathbb{R}}K^\nu(y_*(s),s;y')\nu(y',s)\,dy'\right|.
		\end{split}
	\end{equation}
	Using \eqref{F-nu-p}, \eqref{nu_DK-p}, \eqref{nu_DK2-p}, and noting that \(\partial_y\nu(y_\ast(s),s)=0\) by the choice of \(y_\ast(s)\) in \eqref{y-star}, we deduce from \eqref{nu_eq1-p} that, if \(\nu(y_\ast(s),s)>0\), then
	\begin{equation}\label{nu_temp-p}
		\begin{split}
			\partial_s \nu(y_*(s) ,s) & \leq  Ce^{-s/2}+\int_{\mathbb{R}} \nu(y',s)K^\nu(y_*(s),s;y')\,dy'-D^\nu(y_*(s),s)\nu(y_*(s),s)
			\leq Ce^{-s/2}.
		\end{split}
	\end{equation}
	On the other hand, if $\nu(y_*(s),s)< 0$, then
	\begin{equation}\label{nu_temp--p}
		\begin{split}
			\partial_s \nu(y_*(s) ,s) & \ge - Ce^{-s/2}.
		\end{split}
	\end{equation}
	Fix any $s > s_0$. By the definition of $y_*$, it holds that
	\begin{equation*}
		\|\nu(\cdot,s-h)\|_{L^\infty} \geq |\nu(y_*(s) - h\mathcal{U} (y_*(s),s), s - h)| 
	\end{equation*}
	for any sufficiently small $h>0$. 
	Then, it is straightforward to check that 
	\begin{equation}\label{AP_R1} 
		\begin{split}
			& \lim_{h \to 0^+} \frac{\|\nu(\cdot,s-h)\|_{L^\infty} - \|\nu(\cdot,s)\|_{L^\infty}}{-h} 
			\\
			& \quad \leq \lim_{h \to 0^+}  \frac{ |\nu(y_*(s) - h\mathcal{U} (y_*(s),s), s - h)| -  |\nu(y_*(s), s )| }{-h} 
			\\
			& \quad = (\partial_s + \mathcal{U} (y_*(s),s)\partial_y)|\nu|(y,s)|_{y=y_*(s)}, 
		\end{split}
	\end{equation} 
	provided that the limit in the first line exists. 
	Note that $ \| \nu(\cdot, s)\|_{L^{\infty}}$ is Lipschitz continuous in $s$.\footnote{This can readily be shown using the fact that $\nu(y,s)$ is at least $C^1$ except at $y=0$ and \eqref{Ut_dec_p}. For more details, we refer to the argument in Remark 3.12 of \cite{BKK2}.} Thus by Rademacher's theorem, it is differentiable at almost every $s\in[s_1, s_2]$, ensuring that the limit in the first line exists for almost all $s$. 
	Moreover, the limit in the second line of \eqref{AP_R1} exists due to \eqref{y-star} and the smoothness of $\nu(y, s)$ except at $y=0$.
	
	Since $\partial_y \nu(y_*(s), s) =0$, we deduce from \eqref{AP_R1} that 
	\begin{equation} \label{L-thm-p}
		\begin{split}
			\frac{d}{ds} \| \nu(\cdot, s)\|_{L^{\infty}} 
			& \leq  \begin{cases}
				\partial_s \nu(y, s)|_{y=y_*(s)} & \text{if } \nu(y_*(s),s)>0 \\
				-\partial_s \nu(y, s)|_{y=y_*(s)} & \text{if } \nu(y_*(s),s)<0
			\end{cases}
		\end{split}
	\end{equation} 
	for almost all $s\in[s_1, s_2]$. 
	Hence, integrating \eqref{L-thm-p} with respect to $s$ over $[s_1, s_2]$, and using  \eqref{nu_temp-p} and \eqref{nu_temp--p}, we have 
	\begin{equation*}
		\begin{split} 
			\| \nu(\cdot, s_2) \|_{L^{\infty}} 
			& =  \| \nu(\cdot, s_1) \|_{L^{\infty}} +  \int_{s_1}^{s_2} \frac{d}{ds} \| \nu(\cdot, s)\|_{L^{\infty}} \; ds  
			\\
			& \le \| \nu(\cdot, s_1) \|_{L^{\infty}} + \int_{s_1}^{s_2} C e^{-s/2} \; ds \\
			&  \leq  \| \nu(\cdot, s_1) \|_{L^{\infty}} + C \veps^{1/2} 
			\\
			&= \frac{5}{8} + C \ve^{1/2} < \frac{3}{4}, 
		\end{split} 
	\end{equation*}
provided that $\ve>0$ is sufficiently small. 
This leads to a contradiction to the choice of $s_2$ in \eqref{s2-def}. This proves \eqref{claim_reg-p}, and completes the proof of Lemma~\ref{mainprop_1-p}. 
\end{proof}

		\section{Appendix}

\subsection{A non-vanishing property} \label{sec1.2}

We discuss a local non-vanishing property associated with the Riccati-type mechanism. This supports the nonzero-background condition underlying the Galilean-type change of variables introduced in Section~\ref{sec1.1}. To this end, we first recall some results from \cite{HH, JN}, which will also be used in the continuation argument in Section~\ref{C13_subsec-p}.

\begin{lemma}[Local existence {\cite[Theorem~1]{HH}}]\label{local_exist}
		Let $v_0\in H^s(\mathbb{R})$ with $\textstyle s>\frac{3}{2}$. There exists a maximal existence time $T_{\max}>-\varepsilon$ and a unique solution $v$ to \eqref{Nov} such that
		\begin{equation*}
			v\in C([-\varepsilon,T_{\max});H^s(\mathbb{R}))\cap C^1([-\varepsilon,T_{\max});H^{s-1}(\mathbb{R})).
		\end{equation*} 
\end{lemma}

\begin{lemma}[Blow-up criterion {\cite[Theorem~2.1]{JN}}]
\label{JN_thm}
Let $v_0 \in H^s(\mathbb R)$ with $s \ge 2$, and let $T_{\max}$ be the maximal existence time of the solution $v(x,t)$ to \eqref{Nov}. Then the corresponding solution blows up in finite time (i.e., $T_{\max} < \infty$) if and only if
\begin{equation*}
\lim_{t \nearrow T_{\max}} \inf_{x\in\mathbb R} (v v_x)(x,t) = -\infty.
\end{equation*}
\end{lemma}

By Lemma~\ref{local_exist} and the Sobolev embedding, if the initial data are sufficiently smooth, then there exists a maximal existence time \(T_{\max}>-\varepsilon\) such that the corresponding solution \(v\) remains \(C^1\) in space for all \(t<T_{\max}\). Lemma~\ref{JN_thm} further shows that finite-time blow-up can occur only if \(vv_x\) becomes unbounded from below as \(t\nearrow T_{\max}\). Since this criterion is formulated through a spatial infimum, it does not by itself identify a pointwise blow-up location.

We therefore consider a time \(T_*\leq T_{\max}\) and a point \(x_*\in\mathbb{R}\) such that \(v\) extends continuously to \((x_*,T_*)\) and
\begin{equation*}
\lim_{\substack{(x,t)\to(x_*,T_*)\\ t<T_*}}(vv_x)(x,t)=-\infty.
\end{equation*}
Under this condition, one necessarily has \(T_*=T_{\max}\). Indeed, if \(T_*<T_{\max}\), then the local well-posedness and Sobolev embedding imply that \(v\) remains \(C^1\) near \(t=T_*\), contradicting the above divergence. The following lemma shows that such a point cannot belong to the zero set of \(v\).

\begin{lemma} \label{x=0}
Let $T_*>-\varepsilon$ and $x_* \in \mathbb{R}$. Assume that a solution $v$ to \eqref{Nov} satisfies the following:
\begin{enumerate}[(i)]
\item for each $t<T_*$, $v(\cdot,t)\in C^1(\mathbb R)$;
\item $v$ is continuous at $(x_\ast,T_*)$ and satisfies
\begin{equation*}
\lim_{\substack{(x,t)\to(x_\ast,T_*)\\ t<T_*}} (vv_x)(x,t) = -\infty.
\end{equation*}
\end{enumerate}
Then
\begin{equation*}
v(x_\ast,T_*)\neq 0.
\end{equation*}
\end{lemma}

\begin{proof}
Define $g(x,t):=v(x,t)^2$. Then $g(x,t)\ge 0$ for all $(x,t) \in \mathbb{R} \times [-\varepsilon,T_*)$. Moreover, by the assumptions, $g$ is continuous at $(x_\ast,T_*)$
and
\begin{equation*}
g_x(x,t)=2v(x,t)v_x(x,t) \quad \text{for each }t<T_*.
\end{equation*}
Hence
\begin{equation} \label{limgx}
\lim_{\substack{(x,t)\to(x_\ast,T_*)\\ t<T_*}} g_x(x,t)=-\infty.
\end{equation}

Suppose for contradiction that $g(x_\ast,T_*)=0$, and fix a constant $M>0$. By \eqref{limgx}, there exists $\delta>0$ such that
\begin{equation*}
0<|(x,t)-(x_\ast,T_*)|<\delta \quad \Longrightarrow \quad g_x(x,t)\le -M
\end{equation*}
for $t<T_*$. By continuity of $g$ at $(x_\ast,T_*)$ and the assumption $g(x_\ast,T_*)=0$, we may choose $t_0 <T_*$ sufficiently close to $T_*$ so that
\begin{equation*}
g(x_\ast,t_0)<\frac{M\delta}{4}
\end{equation*}
and $\textstyle |t_0-T_*| < \frac{\delta}{4}$. Set $\textstyle x_0 = x_\ast  + \frac{\delta}{4}$. Then, for every \(\xi\in[x_\ast,x_0]\), we have
\begin{equation*}
0<|(\xi,t_0)-(x_\ast,T_*)|<\delta,
\end{equation*}
and therefore \(g_x(\xi,t_0)\le -M\). Since \(g(\cdot,t_0)\in C^1(\mathbb R)\), the fundamental theorem of calculus gives
\begin{equation*}
g(x_0,t_0)
= g(x_\ast,t_0) +\int_{x_\ast}^{x_0} g_x(\xi, t_0) \, d\xi 
\le g(x_\ast,t_0)-\frac{M\delta}{4}
< 0,
\end{equation*}
which contradicts the nonnegativity of $g$.
Hence $g(x_\ast,T_*)\neq 0$, i.e., $v(x_\ast,T_*)\neq 0$.
\end{proof}

It follows from Lemma~\ref{x=0} that any point at which \(vv_x\to-\infty\) locally uniformly as \((x,t)\to(x_*,T_*)\) must lie away from the zero set of \(v\). This supports the nondegeneracy assumption \(v_0(0)\neq0\) and the Galilean-type change of variables introduced in Section~\ref{sec1.1}. For the solutions considered in this paper, the blow-up location is in fact uniformly separated from the zero set of \(v\), as shown in the following remark.

\begin{remark} \label{center}
Let $\kappa, \xi$ be the modulation functions as in Section~\ref{modul_section}, and let
\begin{equation*}
X(t):=\frac14(t+\varepsilon)+\xi(t).
\end{equation*}
For solutions whose initial data satisfy the assumptions in Section~\ref{Initial_subs} the constraint \(U(0,s)=0\) in \eqref{constraint-p}, together with \eqref{u} and \eqref{UP}, gives
\begin{equation*}
v(X(t),t)=u(\xi(t),t)+\frac12=\kappa(t)+\frac12, \quad -\varepsilon\le t<T_*.
\end{equation*}
By \eqref{modul_dec}, we have \(|\kappa(t)|\le C\varepsilon^{5/6}\) for all \(t<T_*\). Hence, choosing \(\varepsilon>0\) sufficiently small, we obtain
\begin{equation*}
v(X(t),t)\ge \frac14, \quad -\varepsilon\le t<T_*.
\end{equation*}
Thus the modulated blow-up location remains uniformly separated from the zero set of \(v\).
\end{remark}

\subsection{Admissible initial data}\label{app:admissible-data}

For \(\veps>0\) sufficiently small, we construct initial data \(u_0\) satisfying \eqref{in-H-C}--\eqref{H0_bound}, and show that the set of such data has non-empty relative interior in the class satisfying \eqref{in-H-C} and \eqref{init_w_3-p}, with respect to the topology induced by the \(H^5\) norm together with the weighted quantity in \eqref{4.3a-p} and the derivative bounds in \eqref{init_24-p}.

\subsubsection{Bounds for the reference profile \texorpdfstring{\(\overline U\)}{Ubar}}

We first record derivative bounds for \(\overline U\), which justify our choice of the constants \(C_1=2^{12}\), \(C_2=2^{26}\), and \(C_3=2^{39}\) in \eqref{init_24-p}.

\begin{lemma}\label{lem:U-profile-norms}
Let \(\overline U=\overline U_1\) be the profile constructed in Proposition~\ref{Profile-construct}. Then
\begin{equation*}
\|\overline U'\|_{L^\infty}=4,\quad \|\overline U''\|_{L^\infty}<2^{12},\quad \|\overline U^{(3)}\|_{L^\infty}=2^{25},\quad \|\overline U^{(4)}\|_{L^\infty}<2^{39}.
\end{equation*}
\end{lemma}

\begin{proof}
By oddness, it suffices to work on \(y>0\). Set
\begin{equation*}
q(y):=-\overline U'(y),\quad A(y):=\overline U(y)+\frac92y.
\end{equation*}
Then \(q(0)=4\), \(q(y)\to0\) as \(y\to\infty\), and \eqref{ovU} gives
\begin{equation*}
Aq'=-\frac14q(4-q),\quad A'=\frac92-q.
\end{equation*}
Moreover, \eqref{3.4} with \(\beta=1\) in the proof of Proposition~\ref{Profile-construct} gives
\begin{equation*}
A=\frac{(4-q)^{1/2}}{16q^{9/2}}.
\end{equation*}
Hence
\begin{equation*}
q'=-4(4-q)^{1/2}q^{11/2}.
\end{equation*}
Differentiating with respect to \(y\), and using \(q=-\overline U'\), we obtain
\begin{equation*}
\overline U''=4(4-q)^{1/2}q^{11/2},\quad \overline U^{(3)}=32q^{10}(3q-11),\quad \overline U^{(4)}=128q^{29/2}(4-q)^{1/2}(110-33q).
\end{equation*}
Since \(q\) decreases from \(4\) to \(0\), the desired bounds follow by taking the maximum of the above expressions, with absolute values, on \([0,4]\). In particular,
\begin{equation*}
\max_{q\in[0,4]}|\overline U^{(3)}|=\max_{q\in[0,4]}32q^{10}|3q-11|=32\cdot4^{10}=2^{25},
\end{equation*}
where the maximum is attained at \(q=4\).
\end{proof}

\subsubsection{Compatibility of the initial conditions} 

We construct \(\varphi\in C_c^\infty(\mathbb R)\) and an even function \(g\in C_c^\infty(\mathbb R)\), and define the initial datum by
\begin{equation}\label{admissible-data-def-Nov}
u_{0,\veps}(x):=-\frac12+\varphi(x)+\veps^{7/2}\Phi\left(\frac{x}{\veps^{9/2}}\right),
\end{equation}
where
$\Phi(y):=\int_{-\infty}^y g(y')\,dy'$. 
We then verify that \(u_{0,\veps}\) satisfies \eqref{in-H-C}--\eqref{H0_bound}: \(\varphi\) is used to enforce the background value and the Hamiltonian bound \eqref{H0_bound}, while \(g\) is chosen to enforce the matching conditions \eqref{init_w_3-p} and the self-similar slope condition \eqref{4.3a-p}.

Let \(\Theta=(\frac{1}{1296})^{1/9}\), and set
\begin{equation*}
b(y):=\min\left\{\frac{N(m_0)y^2}{4(1+y^2)},\frac{5}{8(1+|y|^{2/9})}\right\}.
\end{equation*}
Choose \(\rho\in(\frac{8\Theta}{5},1)\) sufficiently close to \(\frac{8\Theta}{5}\), and then choose \(\rho_0\) and \(\theta\) such that
\begin{equation*}
\frac{8 \Theta}{5}<\rho_0<\rho<1,\quad \theta<\frac58(1-\rho).
\end{equation*}
By taking \(\rho\) sufficiently close to \(\frac{8\Theta}{5}\), the number \(\theta\) may be chosen sufficiently close to \(\frac58-\Theta\). 

We first choose \(\varphi\in C_c^\infty(\mathbb R)\) such that
\begin{equation*}
\varphi(0)=\frac12,\quad \varphi\equiv\frac12\ \text{near }0,\quad \|\varphi\|_{H^1}^2<\frac{81}{128},\quad |x|^{2/9}|\varphi'(x)|\leq \theta\ \text{for all }x\in\mathbb R.
\end{equation*}
Such a function \(\varphi\) is obtained by smoothing an even cutoff which decreases from \(\frac12\) to \(0\) at the maximal rate allowed by the weight \(|x|^{-2/9}\); choosing \(\theta\) sufficiently close to \(\frac58-\Theta\) gives the strict \(H^1\) bound above.

Next, we construct \(g\) so that
\begin{equation}\label{g-admissible-Nov}
g(y)=\overline U'(y)\ \text{near }y=0,\quad \int_{\mathbb R}g(y)\,dy=0,\quad |g(y)-\overline U'(y)|\leq \rho b(y).
\end{equation}
By \eqref{asymp-y-infty}, there exists \(R\gg1\) such that
\begin{equation*}
|\overline U'(y)|\leq \rho_0 b(y)\quad\text{for }|y|\geq R.
\end{equation*}
Let \(\chi_R\in C_c^\infty(\mathbb R)\) be an even cutoff function such that \(0\leq\chi_R\leq1\), \(\chi_R=1\) for \(|y|\leq R\), and \(\chi_R=0\) for \(|y|\geq2R\). Set \(g_1:=\chi_R\overline U'\). Then \(g_1=\overline U'\) near the origin, \(g_1\in C_c^\infty(\mathbb R)\) is even, and
\begin{equation*}
|g_1(y)-\overline U'(y)|\leq \rho_0 b(y).
\end{equation*}
Indeed, the left-hand side vanishes for \(|y|\leq R\), while for \(|y|\geq R\) it is bounded by \(|\overline U'(y)|\). Moreover, since \(\overline U'\leq0\) on \(\mathbb R\) and \(\overline U'<0\) near the origin, we have
\begin{equation*}
\int_{\mathbb R}g_1(y)\,dy<0.
\end{equation*}
We now add a nonnegative correction to \(g_1\) in the far field so that the resulting function satisfies the zero-integral condition in \eqref{g-admissible-Nov}. Since \(b(y)\sim\frac58|y|^{-2/9}\) as \(|y|\to\infty\), we have
\begin{equation*}
\int_R^L b(y)\,dy\to\infty \quad\text{as }L\to\infty.
\end{equation*}
By taking \(R\) large enough and then choosing \(L>2R\) sufficiently large, we may find an even nonnegative function \(g_2\in C_c^\infty(\mathbb R)\) supported in \(\{2R\leq |y|\leq L\}\) such that
\begin{equation*}
0\leq g_2(y)\leq(\rho-\rho_0)b(y),\quad \int_{\mathbb R}g_2(y)\,dy=-\int_{\mathbb R}g_1(y)\,dy,
\end{equation*}
and such that \(\lVert g_1+g_2-\overline U'\rVert_{W^{3,\infty}}\) is as small as needed away from the prescribed far-field envelope. Setting \(g:=g_1+g_2\), we obtain an even function \(g\in C_c^\infty(\mathbb R)\) satisfying \eqref{g-admissible-Nov}. By Lemma~\ref{lem:U-profile-norms}, \(g\) may be chosen so that
\begin{equation}\label{g-derivative-bounds-Nov}
\|g'\|_{L^\infty}<2^{12},\quad \|g''\|_{L^\infty}<2^{26},\quad \|g'''\|_{L^\infty}<2^{39}.
\end{equation}
Since \(g\) is compactly supported and has zero integral, \(\Phi\in C_c^\infty(\mathbb R)\). Hence \(u_{0,\veps}+\frac12\in C_c^\infty(\mathbb R)\), and \eqref{in-H-C} holds. Moreover,
\begin{equation}\label{rescaled-slope-Nov}
\veps\partial_xu_{0,\veps}(\veps^{9/2}y)=g(y)+\veps\varphi'(\veps^{9/2}y).
\end{equation}

It remains to verify the admissibility conditions for the data defined by \eqref{admissible-data-def-Nov}. Since \(g\) is even and \(\int_{\mathbb R}g=0\), \(\Phi(0)=0\). Since \(g=\overline U'\) near the origin and \(\varphi\equiv\frac12\) near the origin, Proposition~\ref{Profile-construct} gives
\begin{equation*} 
u_{0,\veps}(0)=0,\quad \partial_xu_{0,\veps}(0)=-4\veps^{-1},\quad \partial_x^2u_{0,\veps}(0)=0,\quad \partial_x^3u_{0,\veps}(0)=2^{25}\veps^{-10}.
\end{equation*}
Thus \eqref{init_w_3-p} holds. We then prove \eqref{4.3a-p}. Since \(\varphi'\equiv0\) near the origin and \(|x|^{2/9}|\varphi'(x)|\leq\theta\), the choice of \(\theta\) and \(\rho\) implies, for all sufficiently small \(\veps>0\),
\begin{equation*}
\veps|\varphi'(\veps^{9/2}y)| \leq (1-\rho)b(y)\quad\text{for all }y\in\mathbb R.
\end{equation*}
Combining this with \eqref{g-admissible-Nov} and \eqref{rescaled-slope-Nov}, we obtain
\begin{equation*}
|\veps\partial_xu_{0,\veps}(\veps^{9/2}y)-\overline U'(y)| \leq b(y),
\end{equation*}
which gives \eqref{4.3a-p}. Since \(g\) and \(\varphi'\) are compactly supported, \eqref{rescaled-slope-Nov} also gives
\begin{equation*}
\limsup_{|x|\to\infty}|x|^{2/9}|\partial_xu_{0,\veps}(x)|=0,
\end{equation*}
and hence \eqref{4.3a-p2} holds.

We verify the derivative bounds in \eqref{init_24-p}. By \eqref{num_4} and the definition of \(b\),
\begin{equation*}
-4\leq \overline U'(y)-b(y)\leq \veps\partial_xu_{0,\veps}(\veps^{9/2}y)\leq \overline U'(y)+b(y)<4
\end{equation*}
for sufficiently small $N(m_0)$.
Thus \(\|\partial_xu_{0,\veps}\|_{L^\infty}\leq4\veps^{-1}\). For \(j=2,3,4\), we have
\begin{equation*}
\partial_x^ju_{0,\veps}(x)=\varphi^{(j)}(x)+\veps^{(7-9j)/2}g^{(j-1)}\left(\frac{x}{\veps^{9/2}}\right).
\end{equation*}
Using \eqref{g-derivative-bounds-Nov} and the fact that \(\varphi\) is fixed independently of \(\veps\), we obtain, for all sufficiently small \(\veps>0\),
\begin{equation*}
\|\partial_x^2u_{0,\veps}\|_{L^\infty}\leq2^{12}\veps^{-11/2},\quad \|\partial_x^3u_{0,\veps}\|_{L^\infty}\leq2^{26}\veps^{-10},\quad \|\partial_x^4u_{0,\veps}\|_{L^\infty}\leq2^{39}\veps^{-29/2}.
\end{equation*}
Thus \eqref{init_24-p} holds.

Finally, we verify the Hamiltonian constraint. From \eqref{admissible-data-def-Nov},
\begin{equation*}
u_{0,\veps}+\frac12=\varphi+\veps^{7/2}\Phi\left(\frac{x}{\veps^{9/2}}\right).
\end{equation*}
Since \(\Phi\) is compactly supported and \(\varphi\) is fixed, the self-similar correction contributes \(o(1)\) to the \(H^1\) norm as \(\veps\to0\). Hence
\begin{equation*}
H(-\veps)=\left\|u_{0,\veps}+\frac12\right\|_{L^2}^2+\|\partial_xu_{0,\veps}\|_{L^2}^2=\|\varphi\|_{H^1}^2+o(1)<\frac{81}{128}
\end{equation*}
for all sufficiently small \(\veps>0\). This proves \eqref{H0_bound}.

The construction above also gives non-empty relative interior. Indeed, by the strict choices in \eqref{g-admissible-Nov} and in the estimate for \(\varphi'\), there exists \(\eta>0\) such that
\begin{equation*}
\left|\veps\partial_xu_{0,\veps}(\veps^{9/2}y)-\overline U'(y)\right|\le(1-\eta)b(y)
\end{equation*}
for all \(y\), with the equality at \(y=0\) forced by \eqref{init_w_3-p}. The bounds in \eqref{init_24-p} for \(j=2,3,4\) and the Hamiltonian bound are strict. The equality in the case \(j=1\) of \eqref{init_24-p} occurs only at \(x=0\). At this point, the matching conditions \eqref{init_w_3-p} give
\begin{equation*}
\partial_x^2\bigl(\partial_xu_{0,\veps}+4\veps^{-1}\bigr)(0) =\partial_x^3u_{0,\veps}(0) =2^{25}\veps^{-10}>0.
\end{equation*}
Together with the strict inequalities above, this gives a relative neighborhood of \(u_{0,\veps}\) contained in the admissible class.

\begin{remark}[The parameter \(\beta\) and relaxation of the admissible class]\label{rem:relax}
The value \(2^{25}\) in \eqref{init_w_3-p} comes from the normalization \(\overline U=\overline U_1\). Since
\begin{equation*}
\overline U_\beta(y)=\beta^{-1/2}\overline U(\beta^{1/2}y),\quad \overline U_\beta^{(3)}(0)=2^{25}\beta,
\end{equation*}
replacing \(\overline U\) by \(\overline U_\beta\) changes the third matching condition to
\begin{equation*}
\partial_x^3u_0(0)=2^{25}\beta\varepsilon^{-10}.
\end{equation*}
For \(\beta\) in a sufficiently small open interval containing \(1\), the strict bounds used in the construction of initial data above remain valid, the far-field inequality
\begin{equation*}
\left(\frac{1}{1296\beta}\right)^{1/9}<\frac58
\end{equation*}
continues to hold, and the scaled form of the profile lower bound remains compatible with the fixed slope bound \(\|\partial_xu_0\|_{L^\infty}\le4\varepsilon^{-1}\). Thus the same construction applies after \(\overline U\) is replaced by \(\overline U_\beta\) in \eqref{4.3a-p}. The bootstrap argument would then use \(\overline U_\beta\) in place of \(\overline U\) as the comparison profile. Consequently, the normalized quantity \(\varepsilon^{10}\partial_x^3u_0(0)\) may vary in an open interval containing \(2^{25}\).

More generally, one could formulate a fully \(\beta\)-dependent version of the assumptions, with \(\overline U_\beta\), \(\beta\)-scaled weights, and \(\beta\)-scaled derivative bounds, and then repeat the same argument for any fixed \(\beta>0\) (cf. \cite{KKY}).
\end{remark}

\subsection{\texorpdfstring{Inequalities related to $\overline{U}$ and its derivatives}{Inequalities related to Ubar and its derivatives}}\label{sec_U}

In this subsection, we verify \eqref{num_4236}--\eqref{num_1} for $\overline{U}=\overline{U}_\beta$ with $\beta=1$, the solution to the ODE problem \eqref{ovU}. We begin with \eqref{num_4}.

\begin{lemma} 
Let $\overline{U} := \overline{U}_1$ be the self-similar profile constructed in
Proposition~\ref{Profile-construct}. Then it holds that
\begin{equation}\label{num-4-0}
4+\overline{U}'(y) - \frac{3y^2}{1+y^2}\geq 0
\end{equation}
for all $y \in \RR$. 
\end{lemma}
			
\begin{proof}
We observe that
\begin{equation*}
\frac{y^2}{1+y^2}=\int^y_0 \frac{2y'}{(1+y'^2)^2}\,dy'
\end{equation*}
and, by $\overline{U}'(0) = -4$,
\begin{equation*}
4+\overline{U}'(y)=\int^y_0 \overline{U}''(y')\,dy'.
\end{equation*}
Hence \eqref{num-4-0} can be rewritten as
\begin{equation*}
\int^y_0 \left( \overline{U}''(y')-\frac{6 y'}{(1+y'^2)^2} \right) \,dy' \geq 0.
\end{equation*}
Since $\overline{U}(y)$ is an odd function, without loss of generality, we restrict our consideration to the case $y>0$. Therefore, it is enough to show that 
\begin{equation*}
\overline{U}''(y)\geq \frac{6 y}{(1+y^2)^2}, \qquad y> 0.
\end{equation*}
First, we consider the region $\{ y > 0 : -4\leq \overline{U}'(y)\leq -1\}$.  
From \eqref{U0FAR}, we have 
\begin{equation*}
\left(\frac{4+\overline{U}'}{y^2}\right)^{1/2}=\frac{5040(-\overline{U}')^{9/2}}{(\overline{U}')^4-2(\overline{U}')^3+6(\overline{U}')^2-20\overline{U}'+70}.
\end{equation*}
By this and \eqref{U''_TEMP}, we have
\begin{equation*}
\begin{split}
\frac{(1+y^2)^2}{y}\overline{U}''(y)&=\frac{4(1+y^2)^2}{y}(4+\overline{U}')^{1/2}(-\overline{U}')^{11/2} \geq 4\left(\frac{4+\overline{U}'}{y^2}\right)^{1/2}(-\overline{U}')^{11/2} \\
&=\frac{20160(-\overline{U}')^{10}}{(\overline{U}')^4-2(\overline{U}')^3 + 6(\overline{U}')^2 - 20\overline{U}'+70} \geq  \frac{20160\cdot 1^{10}}{1^4+2\cdot 1^3+6\cdot 1^2+20\cdot 1+70}\geq 6
\end{split}
\end{equation*}
for $y \in \{ y > 0 : -4\leq \overline{U}'(y)\leq -1\}$. Next, for $y \in \{ y>0 : -1\leq \overline{U}'(y) \leq 0\}$, we obtain \eqref{num-4-0} directly as
\begin{equation*}
4+\overline{U}'-\frac{3y^2}{1+y^2}\geq 3-\frac{3 y^2}{1+y^2}\geq 0.
\end{equation*}
Therefore, \eqref{num-4-0} holds for all $y>0$ and it is trivial for $y=0$. This concludes the proof.
\end{proof}

			Next, we prove \eqref{num_2}--\eqref{num_3}.
			\begin{lemma} 
				Let $\overline{U} := \overline{U}_1$ be the self-similar profile constructed in
Proposition~\ref{Profile-construct}. Then it holds that
				\begin{equation*} 
					1+\frac{\overline{U}'}{2}+\frac{2}{1+y^2}\left(\frac{9}{2}+\frac{\overline 	U}{y}\right)\geq \frac{y^2}{2(1+y^2)} 
				\end{equation*}
				and 
				\begin{equation*} 
					\frac{11}{2}+\frac{3\overline{U}'}{2}+\frac{1}{1+y^2}\left(\frac{9}{2}+\frac{\overline{U}}{y}\right) \geq \frac{4y^2}{1+y^2}
				\end{equation*}
				for all $y \in \RR$. 
			\end{lemma}
			\begin{proof}
				Thanks to $\overline{U}$ being an odd function, we can restrict our consideration to the case $y\geq 0$. 
				By \eqref{num_4} and $\overline{U}(0)=0$, we obtain 
				\begin{equation*}
					\overline{U}(y)=\int^y_0 \overline{U}'(y')\,dy'\geq \int^y_0\frac{3y'^2}{1+y'^2}\,dy'-4y= - y-3 \tan^{-1}(y).
				\end{equation*} 
				Hence by \eqref{num_4}, 
				\begin{equation*}
					\begin{split}
						1+\frac{\overline{U}'}{2}+\frac{2}{1+y^2}\left(\frac{9}{2}+\frac{\overline{U}}{y}\right)&\geq \frac{1}{1+y^2}\left(\frac{1}{2}y^2+6\left(1-\frac{\tan^{-1}(y)}{y}\right)\right) 
						\geq \frac{y^2}{2(1+y^2)}.
					\end{split}
				\end{equation*}
				Similarly, 	
				\begin{equation*}
					\begin{split}
						\frac{11}{2}+\frac{3\overline{U}'}{2}+\frac{1}{1+y^2}\left(\frac{9}{2}+\frac{\overline{U}}{y}\right)&\geq \frac{1}{1+y^2}\left(4y^2+3\left(1-\frac{\tan^{-1}(y)}{y}\right)\right) 
						\geq \frac{4y^2}{1+y^2}.
					\end{split}
				\end{equation*}
				This completes the proof. 
			\end{proof}

The remaining inequalities \eqref{num_6} and \eqref{num_1} are proved in the subsequent lemmas.

\begin{lemma} 
Let $\overline{U} := \overline{U}_1$ be the self-similar profile constructed in
Proposition~\ref{Profile-construct}. Then there exists a constant $\textstyle \delta\in(\frac{1}{2},1)$ such that
\begin{equation}\label{num_6-lem}
|\overline{U}''(y)|\frac{y^2+1}{y^2} \int^{|y|}_0{\frac{y'^2}{1+y'^2 }\,dy'} \leq  \delta \left(1+\frac{\overline{U}'}{2} +\frac{2}{y(1+y^2)}\left(\frac{9y}{2}+\overline{U}\right)\right) -\frac{\delta y^2}{48(1+y^2)}
\end{equation}
for all $y \in \mathbb{R}$.
\end{lemma}

\begin{proof}
Throughout the proof, we use the following relations obtained in the proof of Proposition~\ref{Profile-construct}:
\begin{align}
-4 & \leq \overline{U}'(y) < 0, \quad y \in \mathbb{R}, \label{Wbar_neg'} \\
\overline{U} + \frac{9}{2}y &= \left( \frac{4 + \overline{U}'}{2^8(-\overline{U}')^9} \right)^{1/2}, \quad y \geq 0, \label{3.4'} \\
\overline{U}'' &= 4 (4 + \overline{U}')^{1/2} (-\overline{U}')^{11/2}, \quad y \geq 0, \label{U''_TEMP'}
\end{align}
and
\begin{equation} \label{U0FAR'}
- y^{2/9} \overline{U}' = (5040)^{-2/9} (4+\overline{U}')^{1/9} (\overline{U}'^4 - 2\overline{U}'^3 + 6\overline{U}'^2 - 20\overline{U}' + 70)^{2/9}.
\end{equation}
These correspond to \eqref{Wbar_neg}, \eqref{3.4}, \eqref{U''_TEMP}, and \eqref{U0FAR}, respectively.

By the odd symmetry of $\overline U$ and continuity at $y=0$, it is enough to prove \eqref{num_6-lem} for $y>0$. We first rewrite the desired inequality on $y>0$. From \eqref{Wbar_neg'} and \eqref{U''_TEMP'}, we observe that
\begin{equation*}
\overline{U}''(y) = 4(4+\overline{U}'(y))^{1/2}(-\overline{U}'(y))^{11/2} \geq 0,
\end{equation*}
and thus $|\overline U''(y)|=\overline U''(y)$ on $y>0$. Using this, \eqref{num_6-lem} can be rewritten as
\begin{equation} \label{num_6-lem'}
\overline{U}''(y) \frac{y^2+1}{y^2} \int_0^y \frac{y'^2}{1+y'^2} \, dy' \leq \delta \left( \frac{4+\overline U'}2 + \frac{2}{1+y^2} \left(4+\frac{\overline U}{y}\right) - \frac{49}{48}\frac{y^2}{1+y^2} \right), \quad y >0.
\end{equation}
We now split the domain into $\textstyle \{ y >0: - \frac{1}{2} \leq \overline{U}'(y) < 0 \}$ and $\textstyle \{ y > 0: -4 < \overline{U}'(y) \leq - \frac{1}{2} \}$. On each of these regions, we express both sides of \eqref{num_6-lem'} in terms of $\overline{U}'$ and estimate the resulting expressions directly, with the admissible choice $\textstyle \delta=\frac{7}{8}$.

First, we consider the region where $\textstyle -\frac{1}{2} \leq \overline{U}'(y) < 0$. Since
\begin{equation*}
\int_0^y\frac{{y'}^2}{1+{y'}^2} \, dy'\leq y \quad \text{for } y>0,
\end{equation*}
the left-hand side of \eqref{num_6-lem'} is bounded by
\begin{equation*}
\frac{1+y^2}{y}\overline U''(y).
\end{equation*}
To estimate this quantity, we use the relation \eqref{U0FAR'}. Set
\begin{equation} \label{zP}
z:=-\overline{U}'(y), \quad P(z):=z^4+2z^3+6z^2+20z+70 .
\end{equation}
Then \eqref{U0FAR'} gives
\begin{equation} \label{yPz}
y=\frac{\sqrt{4-z}}{5040z^{9/2}}P(z).
\end{equation}
Using this together with \eqref{U''_TEMP'}, we have
\begin{equation*}
\begin{split}
\frac{1+y^2}{y}\overline U''(y) = \frac{\overline U''(y)}{y} + y\overline U''(y) = \frac{20160z^{10}}{P(z)} + \frac{(4-z)zP(z)}{1260}.
\end{split}
\end{equation*}
Note that, on $\textstyle 0 < z\leq \frac12$,
\begin{equation*}
70 = P(0) \leq P(z)\leq P\!\left(\frac12\right)=\frac{1309}{16}
\end{equation*}
and
\begin{equation*}
z^{10}\le 2^{-10}, \quad (4-z)z\leq \frac74.
\end{equation*}
Therefore the left-hand side of \eqref{num_6-lem'} is bounded as
\begin{equation*}
\begin{split}
\mathrm{LHS} & \leq \frac{20160}{70}2^{-10} + \frac1{1260}\cdot\frac74\cdot\frac{1309}{16}  = \frac{4549}{11520}.
\end{split}
\end{equation*}
On the other hand, using $\textstyle 4+ \frac{\overline{U}}{y} \geq 0$ and $\textstyle \frac{y^2}{ 1+y^2} \leq 1$, we obtain
\begin{equation*}
\begin{split}
\mathrm{RHS} &\geq \delta \left( \frac{4+\overline U'}2 - \frac{49}{48} \right)  \geq \delta \left( \frac74-\frac{49}{48} \right).
\end{split}
\end{equation*}
Here the bound $\textstyle 4+ \frac{\overline{U}}{y} \geq 0$ follows from $\overline{U}(0)=0$ and $4+\overline{U}'\geq 0$:
\begin{equation*}
4+\frac{\overline{U}(y)}{y} = \frac1y\int_0^y \left(4+\overline{U}'(y')\right)\,dy' \geq 0.
\end{equation*}
For $\textstyle \delta = \frac{7}{8}$, we have
\begin{equation*}
\mathrm{LHS} \leq \frac{4549}{11520} < \frac{245}{384} \leq \mathrm{RHS},
\end{equation*}
and thus the desired inequality holds.

We next consider the region where $\textstyle -4 < \overline{U}'(y) \leq -\frac12$. Similarly as above, we set \eqref{zP}. Using \eqref{3.4'} and \eqref{yPz}, we find that
\begin{equation*}
\begin{split}
4+\frac{\overline{U}(y)}{y}-\frac{4-z}{3} &= \frac{(4-z)^2(2z^3+9z^2+30z+70)}{6P(z)} \geq 0.
\end{split}
\end{equation*}
Applying this in the right-hand side of \eqref{num_6-lem'} with $\textstyle \delta=\frac78$, we have
\begin{equation} \label{2RHSb}
\mathrm{RHS} \geq \frac78 \left( \frac{4-z}{2} + \frac{2(4-z)}{3(1+y^2)} - \frac{49}{48} \frac{y^2}{1+y^2} \right).
\end{equation}
We now estimate $y^2$ in terms of $4-z$. By \eqref{yPz}, it holds that
\begin{equation*}
\frac{y^2}{4-z} = \frac{P(z)^2}{5040^2z^9}.
\end{equation*}
A direct computation gives
\begin{equation*}
\frac{d}{dz}
\left( \frac{P(z)^2}{5040^2z^9} \right) = -\frac{P(z)(z^4+6z^3+30z^2+140z+630)}{5040^2z^{10}} <0.
\end{equation*}
Therefore, for $\textstyle \frac12 \leq z < 4$,
\begin{equation*}
\frac{y^2}{4-z} \leq \frac{P(\frac12)^2}{5040^2(\frac12)^9} = \frac{34969}{259200} < \frac17.
\end{equation*}
In particular,
\begin{equation} \label{4-z7}
y^2 \leq \frac{4-z}{7} \leq \frac12, \quad \frac1{1+y^2} \geq \frac23.
\end{equation}
Substituting these estimates into \eqref{2RHSb}, we obtain
\begin{equation} \label{2RHS_est}
\mathrm{RHS} \geq \frac78 \left( \frac{4-z}{2} + \frac49(4-z) - \frac{49}{336}(4-z) \right) = \frac{805}{1152}(4-z).
\end{equation}

It remains to estimate the left-hand side of \eqref{num_6-lem'}. Since
\begin{equation*}
\int_0^y \frac{{y'}^2}{1+{y'}^2}\,dy' \leq \int_0^y {y'}^2\,dy' = \frac{y^3}{3},
\end{equation*}
we have
\begin{equation*}
\mathrm{LHS} \leq \frac{y(1+y^2)}{3}\overline{U}''(y).
\end{equation*}
Using \eqref{U''_TEMP'} and \eqref{yPz}, we get
\begin{equation*}
\frac{y\overline{U}''(y)}{3} = \frac{(4-z)zP(z)}{3780}.
\end{equation*}
Together with $\textstyle y^2 \leq \frac{4-z}{7}$ from \eqref{4-z7}, this implies
\begin{equation*}
\begin{split}
\mathrm{LHS} &\leq \frac{(4-z)zP(z)}{3780} \left(1+\frac{4-z}{7}\right) = \frac{(4-z)zP(z)(11-z)}{26460}.
\end{split}
\end{equation*}
Finally, for $\textstyle \frac12 \leq z < 4$,
\begin{equation*}
17640-zP(z)(11-z) =(4-z) \left( z^4(5-z)+36z^3+190z^2+910z+4410 \right) \geq 0,
\end{equation*}
and hence
\begin{equation} \label{2LHS_est}
\mathrm{LHS} \leq \frac{17640}{26460}(4-z) = \frac23(4-z).
\end{equation}
Since $\textstyle \frac23 < \frac{805}{1152}$, \eqref{num_6-lem'} holds by the estimates \eqref{2RHS_est} and \eqref{2LHS_est}.

Combining the two cases proves \eqref{num_6-lem} for all $y>0$. The case $y<0$ follows from the odd symmetry of $\overline U$, since both sides of \eqref{num_6-lem} are even in $y$. The case $y=0$ follows by taking the limit $y\to0$. This completes the proof.
\end{proof}

\begin{lemma}
Let $\overline{U} := \overline{U}_1$ be the self-similar profile constructed in
Proposition~\ref{Profile-construct}. Then there exists a constant $m_0>0$ such that
\begin{multline} \label{num_1'}
	\frac{63}{62}(y^{2/9}+1)|\overline{U}''(y)|\int^{|y|}_0\frac{dy'}{1+y'^{2/9}}
\\
\leq 1-\frac{1}{4(y^{2/9}+1)}+\frac{\overline{U}'}{2}-\frac{2y^{2/9}}{9(y^{2/9}+1)}\left(\frac{9}{2}+\frac{\overline U}{y}+\frac{1}{y}\int_{0}^{y}{\frac{dy'}{y'^{2/9}+1}}\right), \quad |y|\geq m_0.
\end{multline}
\end{lemma}

\begin{proof}
By the odd symmetry of $\overline U$, both sides of \eqref{num_1'} are even in $y$. Thus it is enough to consider $y>0$. By l'H\^opital's rule, we have
\begin{equation} \label{limit79}
\lim_{y\to\infty}y^{-7/9}\int_0^y\frac{dy'}{1+y'^{2/9}}=\lim_{y\to\infty}\frac{9y^{2/9}}{7(1+y^{2/9})}=\frac97.
\end{equation}
For $y>0$, \eqref{U''_TEMP} gives $|\overline U''|=\overline U''$. Hence, by \eqref{asymp-y-infty},
\begin{equation*}
\begin{split}
&\lim_{y\to\infty}y^{2/9}\frac{63}{62}(y^{2/9}+1)|\overline U''(y)|\int_0^y\frac{dy'}{1+y'^{2/9}}\\
&\quad=\frac{63}{62}\lim_{y\to\infty}(1+y^{-2/9})(y^{11/9}\overline U''(y))\left(y^{-7/9}\int_0^y\frac{dy'}{1+y'^{2/9}}\right)=\frac{9\Theta}{31}.
\end{split}
\end{equation*}
We next consider the right-hand side of \eqref{num_1'}. First, we see that
\begin{equation*}
1-\frac{1}{4(y^{2/9}+1)}-\frac{2y^{2/9}}{9(y^{2/9}+1)}\frac92=\frac{3}{4(y^{2/9}+1)}.
\end{equation*}
Moreover, \eqref{limit79} and \eqref{asymp-y-infty} yield
\begin{equation*}
\lim_{y\to\infty}y^{2/9}\left(\frac{\overline U(y)}{y}+\frac1y\int_0^y\frac{dy'}{1+y'^{2/9}}\right)=\frac97(1-\Theta).
\end{equation*}
Therefore, after multiplying the right-hand side of \eqref{num_1} by $y^{2/9}$, we obtain
\begin{equation*}
\begin{split}
&\lim_{y\to\infty}\left\{\frac{3y^{2/9}}{4(y^{2/9}+1)}+\frac{y^{2/9}\overline U'(y)}2-\frac{2}{9}\frac{y^{2/9}}{y^{2/9}+1}y^{2/9}\left(\frac{\overline U(y)}{y}+\frac1y\int_0^y\frac{dy'}{1+y'^{2/9}}\right)\right\}\\
&\quad=\frac34-\frac{\Theta}{2}-\frac27(1-\Theta)=\frac{13-6\Theta}{28}.
\end{split}
\end{equation*}
Since $\Theta=(1/1296)^{1/9}<1/2$, it holds that
\begin{equation*}
\frac{9\Theta}{31}<\frac{9}{62}<\frac{5}{14}<\frac{13-6\Theta}{28}.
\end{equation*}
It follows that \eqref{num_1} holds for all sufficiently large $y>0$. By evenness, it holds for all $|y|\geq m_0$ after choosing $m_0>0$ sufficiently large. This completes the proof.
\end{proof}

\subsection{Maximum principles}		
		
		We present a maximum principle, a modified form of that developed in \cite{BSV}, for use in our analysis. 
		We consider the initial value problem: 
		\begin{equation}\label{IVP-f}
			\begin{split} 
				& \partial_s f(y,s)+D(y,s) f(y,s)+ {U}(y,s)\partial_yf(y,s)=F(y,s)+\int_{\mathbb{R}}f(y',s)K(y,s;y')\,dy',  \quad
				s\ge s_0,  y\in \mathbb{R},
				\\
				& f(y,s_0) = f_0(y), 
			\end{split} 
		\end{equation}
		where $D$, $U$, and $F$ are smooth functions. 
		\begin{lemma}[\cite{BKK2}] \label{max_2} 
			Let $f$ be a classical solution to IVP \eqref{IVP-f}. Let $\Omega \subseteq \mathbb{R}$ be any compact set. 
			Suppose that 
			\begin{subequations}
				\begin{align}
					& \|f(\cdot,s)\|_{L^{\infty}(\Omega)}\leq \mo, \label{max_2_1}
					\\
					&  \|f(\cdot,s_0)\|_{L^{\infty}(\mathbb{R})}\leq \mo, \label{max_2_1'}
					\\
					& \int_{\mathbb{R}}|K(y,s;y')|\,dy'\leq \delta D(y,s)\quad \text{for}\quad (y,s)\in \Omega^c\times [s_0,\infty), \label{max_2_2}
					\\
					& \inf_{(y,s)\in \Omega^c\times [s_0,\infty)}D(y,s)\geq \lambda_D >0, \label{max_2_3}
					\\
					& \|F(\cdot,s)\|_{L^{\infty}(\Omega^c)}\leq F_0, \label{max_2_4}
					\\ 
					& \limsup_{|y|\rightarrow \infty}|f(y,s)| <  2\mo \label{max_2_6}
				\end{align}
			\end{subequations} 
			for some $\mo, F_0, \lambda_D >0$, and $\delta\in(0,1)$ satisfying 
			\begin{equation}\label{D-cond} 
				\mo \lambda_D> \frac{F_0}{ 2 (1-\delta ) }. 
			\end{equation}  
			Then, it holds that  $\|f(\cdot,s)\|_{L^{\infty}(\mathbb{R})}\leq 2 \mo$ for all $s\ge s_0$. 
		\end{lemma}
			For the proof of Lemma~\ref{max_2}, we refer to that of Lemma~SM$2.1$ in the supplementary materials of \cite{BKK2}.
		
		Next, we consider the initial value problem: 
		\begin{equation}\label{f-ivp-2}
			\begin{split}
				& \partial_s f(y,s) + D(y,s)f(y,s) +  U(y,s)\partial_y f(y,s)  =F(y,s),  \quad
				s\ge s_0, \  y\in \mathbb{R},
				\\
				& f(y,s_0) = f_0(y).  
			\end{split}
		\end{equation}  
		The following lemma addresses the spatial decay properties of solutions to the transport-type equation \eqref{f-ivp-2} under suitable assumptions. 
		\begin{lemma}[\cite{BKK2}] \label{rmk2}
			Let $f$ be a classical solution to IVP \eqref{f-ivp-2}. Assume that $U$, $D$, and $F$ are smooth functions satisfying
			\begin{subequations}
				\begin{align*} 
					& \inf_{\{|y| \geq N ,\, s\in[s_0,\infty)\}}  U(y,s) \frac{y}{|y|}  > 0,
					\\
					& \inf_{\{|y| \geq N,\, s\in[s_0,\infty)\}}D(y,s)\geq \lambda_D,
					\\
					& \|F (\cdot, s) \|_{L^{\infty}(|y| \geq N)}\leq F_0e^{-s\lambda_F}
				\end{align*}
			\end{subequations}
			for some {$ \lambda_D, \lambda_F, N, F_0  \geq  0$}.   Then it holds that  
			\begin{equation*}
				\begin{array}{l l}
					\limsup_{| y |\rightarrow \infty}|f(y,s)|\leq \limsup_{|y|\rightarrow \infty}{|f(y,s_0)|}e^{-\lambda_D(s-s_0)}+\frac{F_0}{\lambda_D-\lambda_F}e^{-s\lambda_F} & \quad \text{if } \lambda_D>\lambda_F, \\
					\limsup_{| y |\rightarrow \infty}|f( y ,s)|\leq \limsup_{|y|\rightarrow \infty}{|f(y,s_0)|}e^{-\lambda_D(s-s_0)}+\frac{F_0e^{-s_0\lambda_F}}{\lambda_F-\lambda_D}e^{-\lambda_D(s-s_0)} & \quad \text{if } \lambda_F>\lambda_D.
				\end{array}
			\end{equation*}
			Here the limits are uniform in $s$. 
		\end{lemma}
		For the proof, we refer to that of Lemma~$3.5$ in \cite{BKK2}.


\begin{thebibliography}{10}

			
			\bibitem{BKK2} J. Bae, Y. Kim, and B. Kwon, \textit{Delta-shock for the pressureless Euler--Poisson system}, SIAM J. Math. Anal. 57, 3255-3296 (2025)
			
			\bibitem{BKK} J. Bae, Y. Kim, and B. Kwon, \textit{Structure of singularities for the Euler--Poisson system of ion dynamics}, preprint, arXiv:2405.02557 (2024)
			
			\bibitem{BSV} T. Buckmaster, S. Shkoller, and V. Vicol, \textit{Formation of shocks for 2{D} isentropic compressible {E}uler}, Comm. Pure Appl. Math. 75(9), 2069-2120 (2022)
		

\bibitem{CGLQ}
R.M. Chen, F. Guo, Y. Liu, and C. Qu,
\textit{Analysis on the blow-up of solutions to a class of integrable peakon equations},
J. Funct. Anal. 270(6), 2343--2374 (2016).


			\bibitem{CLWX}
			R.M. Chen, W. Lian, D. Wang, and R. Xu,
			\textit{A rigidity property for the Novikov equation and the asymptotic stability of peakons},
			Arch. Ration. Mech. Anal. 241, 497--533 (2021)
			
			\bibitem{DP}
			H. Dai and M. Pavlov,
			\textit{Transformations for the Camassa--Holm equation, its high-frequency limit and the sinh--Gordon equation},
			J. Phys. Soc. Jpn. 67, 3655--3657 (1998)
			
			\bibitem{HH}
			A.A. Himonas and C. Holliman,
			\textit{The Cauchy problem for the Novikov equation},
			Nonlinearity 25, 449 (2012)
			
			\bibitem{HHK}
A.A. Himonas, C. Holliman, and C.E. Kenig,
\textit{Construction of 2-peakon solutions and ill-posedness for the Novikov equation},
SIAM J. Math. Anal. 50(3), 2968--3006 (2018).
			
			\bibitem{HW}
			A.N.W. Hone and J.P. Wang,
			\textit{Integrable peakon equations with cubic nonlinearity},
			J. Phys. A: Math. Theor. 41, 372002 (2008)
			
			\bibitem{HZ}
			J. Hunter and Y. Zheng,
			\textit{On a completely integrable nonlinear hyperbolic variational equation},
			Physica D 79, 361--386 (1994)
			
			\bibitem{JN}
			Z. Jiang and L. Ni,
			\textit{Blow-up phenomenon for the integrable Novikov equation},
			J. Math. Anal. Appl. 385, 551--558 (2012)
			
			\bibitem{KKS}
			\newblock	Y. Kim, B. Kwon, and W. Shim,
			\newblock	\textit{Asymptotic self-similar blow-up for the regularized Saint--Venant equations},
			\newblock	preprint,
			\newblock	arXiv:2604.03188, 2026.
			
			\bibitem{KKY}
			\newblock	Y. Kim, B. Kwon, and J. Yoon,
			\newblock	\textit{Sharp regularity of gradient blow-up solutions in the Camassa--Holm equation},
			\newblock	preprint,
			\newblock	arXiv:2412.00558, 2024.
			
			\bibitem{Laf}
			S. Lafortune,
			\textit{Spectral and linear stability of peakons in the Novikov equation},
			Stud. Appl. Math. 152, 1404--1424 (2024)
			
			\bibitem{Lai}
			S. Lai,
			\textit{Global weak solutions to the Novikov equation},
			J. Funct. Anal. 265, 520--544 (2013)
			
			\bibitem{LLQ}
			X. Liu, Y. Liu, and C. Qu,
			\textit{Stability of peakons for the Novikov equation},
			J. Math. Pures Appl. 101, 172--187 (2014)
						
			\bibitem{Nov}
			V. Novikov,
			\textit{Generalizations of the Camassa--Holm equation},
			J. Phys. A: Math. Theor. 42, 342002 (2009)
			
			\bibitem{NZ}
			L. Ni and Y. Zhou, \textit{Well-posedness and persistence properties for the Novikov equation}, J. Differential Equations 250, 3002-3021
			(2011)
			
			\bibitem{OP} S.-J. Oh and F. Pasqualotto, \textit{Gradient blow-up for dispersive and dissipative perturbations of the Burgers equation}, Arch. Ration. Mech. Anal. 248(3), 54(61pp) (2024)
			
			\bibitem{WY}
			X. Wu and Z. Yin,
			\textit{Well-posedness and global existence for the Novikov equation},
			Ann. Sc. Norm. Super. Pisa Cl. Sci. (5) 11, 707--727 (2012)		
			
			\bibitem{WY2}
			X. Wu and Z. Yin,
			\textit{A note on the Cauchy problem of the Novikov equation},
			Appl. Anal. 92, 1116--1137 (2013)
			
			\bibitem{WY2011}
			X. Wu and Z. Yin,
			\textit{Global weak solutions for the Novikov equation},
			J. Phys. A: Math. Theor. 44, 055202 (2011)
			
			\bibitem{Y} R. Yang, \textit{Shock Formation of the Burgers--Hilbert Equation}. SIAM J. Math. Anal. 53, 5756--5802 (2021)
			
			
			\bibitem{YLZ}
			W. Yan, Y. Li, and Y. Zhang,
			\textit{The Cauchy problem for the Novikov equation},
			Nonlinear Differ. Equ. Appl. 20, 1157--1169 (2013)
			
			
			\bibitem{YLZ2}
			W. Yan, Y. Li, and Y. Zhang,
			\textit{The Cauchy problem for the integrable Novikov equation},
			J. Differential Equations 253, 298--318 (2012)
			
			
			\bibitem{ZY}
			R. Zheng and Z. Yin,
			\textit{Wave breaking and solitary wave solutions for a generalized Novikov equation},
			Appl. Math. Lett. 100, 106014 (2020)
			

		\end{thebibliography}
	\end{document}